\documentclass{article}

\usepackage[utf8]{inputenc}  %allow utf-8 input
\usepackage[T1]{fontenc}     %use 8-bit T1 fonts
\usepackage{hyperref}        %hyperlinks
\usepackage{url}             %simple URL typesetting
\usepackage{booktabs}        %professional-quality tables
\usepackage{amssymb}
\usepackage{amsfonts}        %blackboard math symbols
\usepackage{nicefrac}        %compact symbols for 1/2, etc.
\usepackage{microtype}       %microtypography
\usepackage{parskip}
\usepackage{enumerate}
\usepackage{enumitem}
\usepackage{graphicx}
\usepackage{xcolor}

\title{Independence of the Continuum Hypothesis: an Intuitive Introduction}
\author{Serafim Batzoglou}

\begin{document}

\maketitle
\begin{abstract}
The independence of the continuum hypothesis is a result of broad impact: it settles a basic question regarding the nature of $\mathbb{N}$ and $\mathbb{R}$, two of the most familiar mathematical structures; it introduces the method of forcing that has become the main workhorse of set theory; and it has broad implications on mathematical foundations and on the role of syntax versus semantics. Despite its broad impact, it is not broadly taught. A main reason is the lack of accessible expositions for nonspecialists, because the mathematical structures and techniques employed in the proof are unfamiliar outside of set theory. This manuscript aims to take a step in addressing this gap by providing an exposition at a level accessible to advanced undergraduate mathematicians and theoretical computer scientists, while covering all the technically challenging parts of the proof.

\end{abstract}
\section{Introduction}

Infinity comes in different sizes. Circa 1874, Cantor proved that the countable infinity $\mathbb{N}$ is strictly smaller than the continuum $\mathbb{R}$ because there is no surjection from $\mathbb{N}$ to $\mathbb{R}$. The Continuum Hypothesis ($\mathrm{CH}$) postulates that there is no infinity bigger than $\mathbb{N}$ and smaller than $\mathbb{R}$. 

We are all familiar with the natural numbers and the real line, so asking whether there is a set of intermediate size is natural and basic. In his famous 1900 speech, Hilbert named the 23 most important open mathematical problems at the time and $\mathrm{CH}$ was number 1 on the list. Astonishingly, the question is independent of the axioms of set theory---essentially, independent of common mathematics. If a proof of either $\mathrm{CH}$ or ``not $\mathrm{CH}$'' is found, the proof can be extended to a contradiction. Proving the independence of $\mathrm{CH}$ was a long journey from Cantor's 1874 formulation, to Hilbert naming it as the number 1 open problem in 1901, to G\"{o}del proving that $\mathrm{CH}$ is consistent with set theory in 1939, to finally Cohen settling the question in 1963 by proving that ``not $\mathrm{CH}$'' is also consistent with set theory.

The independence of $\mathrm{CH}$ is one of the most important mathematical results of all time, and arguably one of the great accomplishments of human intellect. However, as Scott Aaronson indicated in his blog\footnote{https://scottaaronson.blog/?p=4974}, ``perhaps no other scientific discovery of equally broad interest remains so sparsely popularized''.  Aaronson and Timothy Chow (2007), who wrote an excellent readable exposition on the topic, call this an ``open exposition problem''. My aim here,  besides having fun writing on the topic, is to take a stab at this exposition problem aiming for an audience of computer scientists and amateur mathematicians, assuming discrete math or theoretical computer science major-level background. I aimed to write a manuscript that would have been helpful to me in my self study of the topic, and decided to forego brevity and instead motivate and elaborate on the key ideas and especially on the notion of ``genericity''.

\textbf{Layers of complexity in understanding the proof.} The proof is generally considered exceptionally hard---especially in one direction, consistency of ``not $\mathrm{CH}$''. To prove this, Cohen devised the powerful technique of \textit{forcing} that has since  become a main workhorse of set theorists. There are three independent sources of complexity in this proof:
\begin{enumerate}
\item The general unfamiliarity of broad audiences with foundational mathematics: set theory, model theory and proof theory.
\item As S. Aaronson described eloquently\footnote{https://scottaaronson.blog/?p=4974}, within the proof there are three distinct layers whose boundaries should stay crystal clear but are often confused. The first layer is syntactic and constructive, and is conducted outside of ZFC. It is a terminating algorithm that takes a ZFC proof of $\mathrm{CH}$ as input---as a string of letters---and converts it into a ZFC proof of contradiction. This layer can be conducted in a formal system such as Peano Arithmetic. The second layer is set theoretic, conducted within ZFC, and involves showing that theorems of ZFC can be reflected in a model $\mathrm{M}$ so that a ZFC proof of $\mathrm{CH}$ can be converted to a proof of $\mathrm{CH}$ within $\mathrm{M}$. The third layer is again conducted within ZFC, and is the forcing construction within a model $\mathrm{M}$ ``looking outside''. The set theoretic layers are highly non-constructive and involve transfinite constructions of sets with heavy use of the axiom of choice. Even so, seen as syntax the proof is a fully specified terminating algorithm.
\item Third, the forcing construction itself is complex especially in the definition and motivation of \textbf{\textit{genericity}}. This concept is usually introduced at the front without sufficient explanation or motivation, while its hidden purpose and power are made clearer later in the proof, when genericity allows defining the forcing relation that enables key ZFC axioms like Comprehension to be transferred from the model $\mathrm{M}$ to the bigger model $\mathrm{M}[G]$. Here, I  reverse the order of exposition from the end goal to the front: first motivate the need to transfer axioms ``outside'' $\mathrm{M}$, and then introduce genericity as a natural way to accomplish that. Once this need to transfer axioms outside $\mathrm{M}$ is understood, one is \textit{forced} to introduce genericity.
\end{enumerate}

Once spelled out, the modular structure of the proof facilitates its understanding. First understand the main syntactic argument using forcing within ZFC as a black box, and separately understand the forcing construction.

\textbf{Recommended background for reading this manuscript.} The aim of this exposition is to require minimal background. Unfortunately, for this topic the minimal is still high: college--major-level mathematics or theoretical/algorithmic computer science is probably necessary. In addition, familiarity with first-order logic and especially with the classical theorems by G\"{o}del---Incompleteness, Completeness and Compactness---is recommended. These results are highly accessible compared to the independence of CH. Readers are referred to background material on these topics (Kim B.\footnote{https://web.yonsei.ac.kr/bkim/goedel.pdf}; Smith 2013; Nagel and Newman 2001; Stanford Encyclopedia of Philosophy\footnote{https://plato.stanford.edu/entries/goedel-incompleteness/}; Buildt 2014, Batzoglou 2021).

\textbf{Expositions of forcing and the continuum hypothesis.} Many excellent expositions of the independence of $\mathrm{CH}$ are available, even though most of them are aimed at the committed student of set theory or the specialist in the field. A few words on my self study of the topic may be helpful to a reader embarking in self study. After reading the respective chapters in Jech (2002)---an advanced comprehensive exposition of set theory that is probably not best for self study of forcing---I found Forcing for Mathematicians by Weaver (2014). The name says it all, and reading this brief book a couple of times made everything clear. Then I got additional insight by reading Smullyan and Fitting (2010). Additionally, Cohen's book (1966) is highly readable and clear. The highly readable high-level exposition by Chow (2007) might be a good place to start. Other references that are considered outstanding are Bell (1977) and Kunen (1980). The exposition here most closely follows the one by Weaver: the syntactic approach of introducing a symbol $\mathrm{M}$ postulated to be a countable transitive model, as well as the choices of foregoing Boolean algebras in favor of posets and using generic ideals rather than generic filters, are in my opinion excellent choices that clarify the presentation. In addition, I blend in useful elements from Smullyan and Fitting (2010) and Cohen (1966), reverse some of the order of the presentation to motivate the defined structures and especially genericity, and provide examples and discussion to help the reader absorb the concepts and proofs.

\section{Overview}

This section provides an overview of basic background material and a summary of the proof of independence of $\mathrm{CH}$. The aim is to to give the reader a lay of the land before delving into the long exposition in Sections 3 and 4.

\subsection{Levels of Infinity and the Continuum Hypothesis}

A set A  is bigger than a set B when there is no way to pair every element of $A$ with a distinct element of $B$. Formally, a function $f: A \to B$ is $1-1$ or an \textit{injection} if distinct elements of $A$ map to distinct elements of $B$; a \textit{surjection} if every element of $B$ is the image $y=f(x)$ of some element $x$ of $A$; and a \textit{bijection} if it is injective and surjective. A bijection between $A$ and $B$ establishes that they are of the same size, while if no surjection exists then $|A| < |B|$. This is obvious for finite sets, and the same notion applies to infinite sets. 
\begin{quote} 
\textit{Infinity comes in different sizes, and set theory can be thought of largely as the study of infinity. }
\end{quote}

Importantly, a set is \textbf{countable} if there is a surjection from $\mathbb{N}$ to the set. Let's see some examples of countable sets:
\begin{itemize}

\item $\mathbb{N}^2 = \{(i, j) \mid i,j \in \mathbb{N} \}$. One simple way to enumerate $\mathbb{N}^2$ is $(0,0)$, $(0,1)$, $(1,0)$, $(0,2)$, $(1,1)$, $(2,0)$, $(0,3)$, $...$ by going through all pairs of sum $0, 1, 2, ...$ successively. Similarly, $\mathbb{N}^k$ is countable for any $k \in \mathbb{N}$: just enumerate the $k$-tuples of sum $0, 1, 2, 3, \ldots$ in order.
\item A countable union of countable sets is countable. Let $A_1, A_2, A_3, \ldots$ be countable sets. Each $A_i$ can be enumerated, $a_{i1}, a_{i2}, \ldots$. Enumerate all of them in stages: at stage $k$ enumerate the $k$th element of sets $A_1, \ldots, A_{k-1}$ followed by the first $k$ elements of $A_k$: $a_{1k}, a_{2k}, \cdots, a_{(k-1)k}, a_{k1}, \ldots, a_{kk}$. Every element in $\bigcup_{i=1}^{\infty} A_i$ is $a_{kl} \in A_k$ for some $k, l$ and it is thus enumerated at stage $max\{k,l\}$.\footnote{Note that to enumerate the sets, the axiom of countable choice is used.}
\item $\bigcup_{k = 1} ^\infty \mathbb{N}^k$ is countable because it is a countable union of countable sets.\footnote{Similarly, $\bigcup_{l=1}^{\infty} \: (\bigcup_{k = 1} ^\infty \mathbb{N}^k)^l$ is countable because it is a countable union of countable sets. Each of its member is an $l$-tuple $(w_1, \ldots, w_l)$ where each $w_i$ is a $k_i$-tuple of elements $(n_1, \ldots, n_{k_i})$ where $n_j \in \mathbb{N}$, and $k_i: \: 1\leq i \leq l$ and $l$ are all arbitrary $\in \mathbb{N}$. The above process can continue for $m$ steps $\bigcup_{l_m = 1}^{\infty} \: (\ldots ((\bigcup_{l_1 = 1}^{\infty} \mathbb{N}^k)^{l_1})^{\cdots})^{l_m}$ for any $m$, yielding a countable set, call it $N_m$. The union of all these has a special name, $\varepsilon_0 := \bigcup_m N_m$ and is a countable union of countable sets, therefore countable.}
\end{itemize}

However, $\mathbb{R}$ is bigger than $\mathbb{N}$. This is Cantor's famous theorem and has a simple proof. Let $x_1, x_2, x_3, \ldots$ be an enumeration of all members of $(0,1)$. Represent $x_i$ by $0.x_{i1} x_{i2} \ldots$, in decimal, following the convention that $0. \ldots x_{ik}99999...$ is picked over $0. \ldots (x_{ik}+1)00000...$. Now construct $y = 0.y_1 y_2 \ldots$ as follows: if $x_{ii} = 1$ let $y_i = 2$ and if $x_{ii} \neq 1$ let $y_i = 1$. Then, $y$ cannot be in the list because it differs from $x_i$ in its $i^{th}$ digit, contradicting the assumption that the $x_i$s completely enumerate $(0,1)$. Therefore $(0,1)$ and $\mathbb{R}\supset (0,1)$ are not countable. The size of infinity represented by the real line is called the \textbf{continuum}, denoted by $c$. The following are examples of sets of size continuum:

\begin{itemize}
\item $|(0,1)| = |\mathbb{R}|$: $f = (2x-1) / (x - x^2)$ is a  bijection between $(0,1)$ and $\mathbb{R}$.
\item $|\mathbb{R}| = |\mathbb{R}^k|$. Can you think of a way to show that? Here is a sketch: first, map them to $(0,1)$ and $(0,1)^k$, respectively. Then, write $x \in (0,1)$ in decimal: $x = 0.x_1 x_2 x_3 \ldots$ where $x_i \in \{0, 1, ..., 9\}$. Again follow the convention that among $0.x_1 \ldots x_k 99999.... = 0.x_1 \ldots (x_k +1) 00000......$, with $x_k \leq 8$, the former is picked. Finally, interleave digits. For example, $(0,1)^2$ maps to $(0,1)$ by mapping $(x, y) = (0.x_1 x_2 \ldots, 0.y_1 y_2 \ldots)$ to $0.x_1 y_1 x_2 y_2 \ldots$.\footnote{Interleaving infinite decimals is a bit tricky. For example, de-interleaving $0.94040404...$ leads to $0.90000....$, which is not allowed. It is a bit cumbersome but not hard to deal with this subtlety.} 
\item $|\mathbb{R}| = |\mathbb{R}^\mathbb{N}|$. There is a bijection between $\mathbb{R}$ and the set of tuples $(x_1, x_2, \ldots)$ with each $x_i \in \mathbb{R}$: first map $\mathbb{R}^{\mathbb{N}}$ to $(0,1)^\mathbb{N}$. Then, each $(x_1, x_2, \ldots)$ can be thought of as a $\mathbb{N} \times \mathbb{N}$ table of single digits $[x_{ij}]$ with $i,j \in \mathbb{N}$. By using the map from $\mathbb{N}$ to $\mathbb{N}^2$ just wrap those digits to a single sequence of digits $y_1, y_2, \ldots$ that represents a single real number in $(0,1)$. 
\item $\mathcal{N} := \omega \to \omega$, the set of functions from $\omega$ to $\omega$, also known as the Baire space. It is easy to see that $|\mathcal{N}| = |\mathbb{R}|$; the members $r \in \mathcal{N}$ are infinite sequence over $\mathbb{N}$ and are often called reals.
\end{itemize}

Two concepts are key for measuring infinity: the concepts of ordinal and cardinal. Those concepts are covered in Section 3.2. For now, it suffices to say that ordinals are a way of counting past $\mathbb{N}$, in the following sense. Denote $\mathbb{N}$ with $\omega$ as is customary when we think of it as an ordinal. Its successor is defined to be  $\omega +1 := \omega \cup \{\omega \}$. This successor operator $x+1 := x \cup \{x\}$ can be applied repeatedly to yield $\omega+2, \ldots, \omega\cdot 2 := \{0,1,\ldots, \omega, \omega+1, \ldots\}$, and continue to $\omega\cdot 3, \ldots \omega^2, \ldots, \omega^3, \ldots, \omega^{\omega}, \ldots$, each of these sets being distinct from each other and well ordered. \textbf{Cardinals} are a way of sizing infinity; for cardinals $\kappa, \lambda$ we say $\kappa < \lambda$ when there is no surjection from $\kappa$ to $\lambda$ (and consequently, there is a surjection but no bijection from $\lambda$ to $\kappa$). 

\textbf{The power set.} Given any set $S$, the set of all its subsets is called the \textit{power set} of $S$, denoted by $\mathcal{P}(S)$ or by $2^S$, the latter notation emphasizing the equivalence between a subset and a function $g: S  \to\{0,1\}$ that chooses which elements of $S$ to include in the subset. $\mathcal{P}(S)$ is always larger than $S$: assume for contradiction that there is a bijection $f: S \to\mathcal{P}(S)$. Consider the subset $D$ of $S$ of elements that do not belong to their image under $f$: $D = \{s \in S \mid s \not \in f(s)\}$.\footnote{$D$ is always a subset of $S$ and thus a member of $\mathcal{P}(S)$. If there is no $s$ in $D$, then $D$ is $\emptyset \in \mathcal{P}(S)$, and therefore $f$ has to cover it.} Since $f$ is a bijection, let $d \in S$ be the unique element such that $f(d) = D$. Is $d$ a member of $D$? Well, if $d \not \in D$, then $d \not \in f(d)$ so it satisfies the condition for being in $D$. But if $d \in D$, then $d \in f(d)$ so it does not satisfy the condition of being in $D$. This is a contradiction, therefore $f$ cannot be a bijection.

Therefore, the power set is always bigger than the original set. Infinite cardinals are ordered: $\aleph_0 = \omega, \aleph_1, \aleph_2, \ldots, \aleph_{\omega}, \ldots$. The \textbf{\textit{continuum hypothesis} ($CH)$} states that there is no infinity strictly between $\mathbb{N}$ and $\mathbb{R}$, more commonly stated as $2^{\aleph_0} = \aleph_1$.The \textbf{\textit{generalized continuum hypothesis ($\mathrm{GCH}$)}} states that there is no infinity strictly between an infinite set and its power.

So, is there a set of size between $\mathbb{N}$ and $\mathbb{R}$? This is a simple question that ought to have a clear answer.  Astonishingly, the question is independent of the rest of mathematics: if either case is proven, \textit{the proof can be converted to a proof of contradiction by a complex but well-defined algorithm}. One thing becomes immediately clear: no set $S$ can be described mathematically and demonstrated to be between $\mathbb{N}$ and $\mathbb{R}$. In that sense, such a set does not exist. However, the existence or non-existence of such a set can be made into a new axiom that introduces no new contradiction to existing axioms, and which alters the properties of $\mathbb{R}$.

\subsection{Outline of this manuscript}

An outline/summary of the manuscript is provided below. \textit{It is not self-contained} so some readers may wish to skip to Section 3.

Starting in Section 3, the manuscript covers preliminary background on axiomatic set theory, including the formulation of infinity in the concepts of ordinals and cardinals and the formulation of the universe(s) of sets, and model theory, or how to define mini-universes of elements (such as numbers, fractions, and in our case sets) that obey a collection of axioms. Readers familiar with the ZFC axioms, ordinals, cardinals, sets and classes, relativization and reflection, are welcome to skip these sections.

\textbf{Set theory.} Section 3.1 covers basic axiomatic set theory. All ordinary mathematics can be formalized using the versatile concept of a set. Familiar entities such as numbers, vector spaces, groups, graphs, and functions and relations over them, can all be expressed as sets. A collection of axioms that formalize the notion of a set can serve as a foundation of all mathematics. One such axiomatization is the Zermelo-Fraenkel set theory with the Axiom of Choice (ZFC). It includes the standard first-order logical axioms and inference rules, plus set--theory-specific axioms: seven axioms and two axiom templates that take as input a formula $\phi$ and generate an axiom where $\phi$ is a substring.\footnote{The need for a separate axiom for every formula stems from a desire for all quantifiers to quantify over sets. A finite axiomatization would require quantifying over properties over sets, i.e., formulas. ZFC cannot ``talk'' about formulas.} In this way, ZFC is an infinite but compactly described (recursive) collection of axioms that specify properties of sets, upon which all mathematical concepts can be built.

Section 3.2 covers the basic concepts of ordinals and cardinals. The brief Section 3.3 covers the concept of proper class: a collection of sets that satisfy a ZFC formula that is "too large" to be a set. The collection of all ordinals is a canonical example of a proper class. In that section, the universe of sets as a ``cumulative hierarchy" and the rank of each set within this hierarchy are also defined.  

\textbf{Models.} Section 3.4 covers basic model theory. A model for a collection of axioms is a set together with an interpretation of the constants, functions and relations mentioned in the axioms within the set: constants are mapped to elements of the model, and functions and relations are mapped to functions and relations across the elements of the model. For example,  consider the axioms of group theory: (1) associativity: $\forall x \forall y \forall z \: (x \cdot y) \cdot z = x \cdot (y \cdot z)$; (2) identity: $\forall x \: e\cdot x = x \cdot e = x$; (3) inverse: $\forall x \: \exists x^{-1} \: x \cdot x^{-1} = x^{-1} \cdot x = e$. $\mathbb{Z}$ is a model of group theory if the identity $e$ is mapped to $0$, the function $\cdot(x,y)$ is mapped to $x+y$, and the inverse of $n$ is $-n$. The finite groups $\mathbb{Z}_n$ of addition modulo $n$ provide infinitely many additional models. Countless other models exist such as $I_h = A_5\times\mathbb{Z}_2$, $SO(2)$ and $GL(n, \mathbb{R})$. Importantly, \textit{every logically consistent collection of axioms has a model}. This is the famous G\"{o}del completeness theorem. In this manuscript, a more specific result will play a key role: any finite consistent collection of axioms involving the $\in$ relation\footnote{The specific property needed from $\in$ here is \textit{extensionality}: $x$ and $y$ are equal if and only if for all $z$, $z \in x$ if and only if $z\in y$.} admits a countable transitive set model $M$: a countable set $M$ such that whenever $x \in y \in M$ then $x \in M$. This result is called a \textbf{reflection principle} because the finite collection of axioms is \textbf{reflected} in the set $M$. As we will see, the reflection principle is not hard to prove but its use can be counter-intuitive and confusing. In the case of  independence of $\mathrm{CH}$, the subject matter is existence of uncountable sets and thanks to the reflection principle uncountability will be reflected/modeled within a countable set $M$! This subtle point will be made clear in Section 3.4. Briefly, the definition of a set $S$ being ``uncountable'' within $M$ is that there is no function $f \in M$ that is a bijection from $\mathbb{N}$ to $S$.\footnote{If $M$ is a model of ZFC, then $\mathbb{N}\in M$.} Therefore, a set $S\in M$ can be countable ``in reality'' but uncountable with respect to $M$. This property will be key to a devious construction that explores it to ``break'' $\mathrm{CH}$.

The statement that a given $S$ is uncountable is a ZFC formula $\mathrm{Uncountable}(x)$ taking a single argument. To say that a set $S$ is uncountable \textit{within} $M$, we \textbf{relativize} the formula $\mathrm{Uncountable}(x)$ to $M$ by changing the "there is no surjection $f$ from $\mathbb{N}$ to $S$" to "there is no surjection $f\in M$ from $\mathbb{N}$ to $S$". 

\begin{quote}
\textbf{Notation.} For every formula $\phi$, we denote by $\phi^M$ its version \textbf{relativized} to $M$ where all variables take values from $\mathrm{M}$. Then $\aleph_1^M$, $\aleph_2^M$ and  $(2^\mathbb{N})^M$ are the sets within $M$ that satisfy the formulas $\psi^M$ where $\psi$ is the  formula that defines the properties of the corresponding set.
\end{quote}

In the long Section 4, the manuscript outlines the proof of the independence of $\mathrm{CH}$, omitting some details but covering all the technically and conceptually challenging parts, with commentary that is aimed at the non-specialist reader. A brief non self-contained summary follows here; concepts mentioned in this summary will be fully exposed in Section 4. 

\textbf{Concluding sections.} Section 5 provides a bullet-point review of the proof so that the reader anchors the big picture in a summary format. Section 6 covers two key proofs that were omitted in the main exposition---the model-theoretic preliminary of Reflection including the Mostowski Collapse, and the proof of the consistency of $\mathrm{CH}$ using forcing. Section 7 covers alternative ways to prove independence of $\mathrm{CH}$; all these alternatives use a form of forcing, so in a sense there is only one grand proof. Section 8 concludes with my personal remarks on lessons learned.

\subsection{Non self-contained summary of the proof} 

The following summary provides a high-level map of the proof. Some readers may find this a useful introduction to the main structure and ideas of the proof. However, it is not self contained and many concepts will only become clear in the main exposition, so the reader can optionally skip to Section 3.

Assume for contradiction that there is a proof of $\mathrm{CH}$. A proof is a list of formulas $\phi_1, \ldots, \phi_n = \mathrm{CH}$ that are either axioms or logical consequences of previous formulas in the list. The axioms in this list, plus finitely many additional axioms that will be needed, are all collected in a list that we may denote again by $\phi_1, \ldots, \phi_k$. Using the Reflection Principle some countable transitive set $M$ reflects them. Therefore, $\mathrm{CH}$ is provable within $M$. $\mathrm{CH}$ says that $2^{\aleph_0} = \aleph_1$. Within $M$, these uncountable sets are represented by countable ordinals that are uncountable in $M$. Therefore there is a $(2^{\aleph_0})^M$ set, a $\aleph_1^M$ set, and so on. Each of these sets is a countable set that is uncountable in $M$ in the sense that there is no bijection in $M$ from $\aleph_0$ to that set. (As it turns out, $\aleph_0^M = \aleph_0$.) Now comes the fun part. Because $M$ is countable, we can expand it in some drastic ways: because $\aleph_2^M$ is countable, there exist $\aleph_2^M$-many distinct subsets of $\aleph_0$ that are not in $M$. Let $G$ denote a table of $\aleph_2^M$-many distinct rows, each of which is a subset of $\aleph_0$ not in $M$. Attaching $G$ to $M$ will not result in a model immediately, because $M\cup \{G\}$ will be missing new sets that are definable by the set theory axioms (covered in Section 3.1) using $G$ as a parameter. The idea is to attach $G$ as well as additional subsets so that the new expanded set, call it $M[G]$, is a model of ZFC.\footnote{Actually, it suffices that it is a model of the axioms required in the hypothesized proof of $\mathrm{CH}$.} Moreover, it should be ensured that $\aleph_2^{M[G]}$ is equal to $\aleph_2^M$ in the sense that the same exact countable ordinal plays the role of the second smallest uncountable infinity in both $M$ and $M[G]$. Extending $M$ with $G$ while making sure that the finite list of axioms we started with still provably hold in $M[G]$ is super delicate, and is done by a technique called \textbf{forcing}. Once this is accomplished, within $M[G]$ it is simultaneously true that $(2^{\mathbb{N}})^{M[G]} \geq $ $\aleph_2^{M[G]}$ and also $CH^{M[G]}$, which says $(2^{\mathbb{N}})^{M[G]} = \aleph_1^{M[G]}<\aleph_2^{M[G]}$, a contradiction.

\begin{quote}
\textit{The independence proof is an algorithm whose input is a ZFC proof of $\mathrm{CH}$ or "not $\mathrm{CH}$" and whose output is a ZFC proof of contradiction.}
\end{quote}

\textbf{The proof comes in two main parts:} (1) the logical part whereby a proof of $\mathrm{CH}$ (or $\lnot \mathrm{CH}$, which we will cover after the main exposition) is relativized to a countable transitive set $M$; (2) the forcing part of extending $M$ by adding the table $G$ of $\aleph_2^M$-many distinct rows, defining $M[G]$ that contains $M$ and $G$ as well as additional sets, and then showing that the resulting expanded set $M[G]$ is a model of the axioms used in the hypothesized proof of $\mathrm{CH}$. The logical argument is not technically hard. However, it is a nasty source of confusion. The argument is carried out outside ZFC. It is an argument \textbf{about} the axioms of ZFC that can be expressed in informal mathematical language or carried out formally in Peano Arithmetic (PA). Importantly, it defines a terminating algorithm whose input and output are ZFC proof syntax. The forcing argument (2) is technically much harder. It is conducted within ZFC. It has different layers of ``truth'': because of the models $M$ and $M[G]$, any sentence $\phi$ has in principle different truth values as $\phi$, $\phi^{\mathrm{M}}$, and $\phi^{\mathrm{M}[G]}$. For example $(2^{\mathbb{N}})^M$ is countable and at the same time uncountable$^M$. The proof breaks $CH^M$ not  $\mathrm{CH}$. This last point can be confusing and even infuriating. Nothing is concluded about the ``true'' truth value of $\mathrm{CH}$. All that is shown is that if a proof of $\mathrm{CH}$ existed, the axioms used in the proof would cogently hold within a countable transitive set $M$, then an extension $M[G]$ would be demonstrated where the same axioms hold, and proofs of $CH^{M[G]}$ as well as $\lnot \mathrm{CH}^{M[G]}$ would go through resulting in contradiction. Does this mean that $\mathrm{CH}$ cannot hold? Not really. For example, $\mathrm{CH}$ can be appended as an axiom. Then the axiom can hold as $CH^M$ in a countable transitive $M$, and $M$ can be extended to $M[G]$ as before, but then the $CH^{M}$ axiom will \textit{not} carry to $CH^{M[G]}$ so the ``proof'' of $\mathrm{CH}$---a one-liner just stating the axiom---will not carry to $M[G]$ and no contradiction will ensue.

\textbf{Expanding the ground model with new subsets of $\mathbb{N}$.} In the second part, the ground model $M$ is expanded with a table $G$ that encodes $\aleph_2^M$-many distinct subsets of $\aleph_0$ that are not in $M$. For the resulting set $M[G]$ to be a model of ZFC it has to include all sets that are constructible from sets in $M$ combined with $G$. One way to construct new subsets is by taking preimages under $G$ as follows. $G$ is a table of $\aleph_2^M \times \aleph_0$ entries $0/1$, each $(\alpha, k)$ entry indicating whether the $\alpha$th subset of $\aleph_0$ contains (1) or doesn't contain (0) element $k$. This table is not in $M$, however every finite subset of it is in $M$ because \textit{every model contains every finite set} (and incidentally it also contains every Turing machine and every recursive function). Those finite subsets of $G$ are part of a larger collection that is in $M$, the collection of all finite subsets of any $\aleph_2^M\times \aleph_0$ $0/1$ table, and $M$ does not know which ones are the former.\footnote{The expression ``$M$ does not know'' is anthropomorphism for the syntactic fact that no statement pertaining to subsets in $G$ relativizes to $M$ and no set that isolates these subsets is a member of $M$.} Let's denote by $\mathbb{P}$ the collection of finite fragments of any 0/1 table of size $\aleph_2^M\times \aleph_0$. Now, consider any element $x$ in $M$, and any relation $r$ in $M$ that is a subset of $x\times \mathbb{P}$. The preimage of $r$ by $G$, i.e., all elements $(y, p)$ of $x \times \mathbb{P}$ that are in $r$ and such that $p \subset G$, should also be in $M[G]$ because $M[G]$ has the $G$ ``key''. This preimage operation can be repeated infinitely many times and all the results need to be in $M[G]$ if $M[G]$ has any hope of being a model. This need motivates the definition of $M[G]$: it is precisely the collection of all such preimages. At first glance it would be wishful thinking that such a simple definition of $M[G]$ should work, and miraculously it does.

\textbf{Forcing the expanded set to be a model.} The next task is to show that $M[G]$ defined thus is a model. The desired property for $M[G]$ is that for every axiom of ZFC relativized to $M$, $\phi^M$, there is a ZFC proof of $\phi^{M[G]}$. To achieve that for certain strong axioms $\phi$ such as Comprehension, $\phi^{M}$ should be invoked to derive $\phi^{M[G]}$. However, when invoking $\phi^{\mathrm{M}}$, there is no access to $G$; only to finite parts $p\in \mathbb{P}$ of G. So it turns out that to invoke $\phi^{\mathrm{M}}$ to prove $\phi^{\mathrm{M}[G]}$, the truth of  $\phi^{\mathrm{M}[G]}$ has to be decided on whether some finite $p\in \mathbb{P}$ is in $G$. Such a $p$ is said to \textit{\textbf{force}} $\phi^{\mathrm{M}[G]}$. This requirement immediately leads to the concept of genericity of $G$: because each property of $G$ is conferred by some finite $p \in M$, all properties of $G$ are generic with respect to $M$: no specific infinite ``AND'' of new properties can be formed. An infinite ``AND'' encoding a specific new property in $G$ not present in $M$, such as for example a specific new subset (as opposed to an unspecified, generic new subset), or the equality between two specific rows of $G$, would correspond to a sentence $\phi^{\mathrm{M}[G]}$ that cannot be made true or false based on a finite $p \in \mathrm{M}$. Genericity of $G$ is defined elegantly in terms of the concept of density of subsets of $\mathbb{P}$. $G$ is defined to be generic if it meets every dense subset $D\subset \mathbb{P}$, $D\in M$ with some $p\in D$ that is a subset of $G$, because to avoid $D$ would require the truth of a new infinite ``AND''. It is then proven that a generic $G$ exists as long as $M$ is countable. And that's why we insist on countable models. The relation between a sentence $\phi^{\mathrm{M}}$ and the parts $p$ whose presence in $G$ ensures the corresponding sentence $\phi^{\mathrm{M}[G]}$ is true is called the \textit{forcing} relation in $M$. A key theorem shows that every sentence of $M[G]$ is forced true or false by the forcing relation in $M$. Based on that theorem, with some technical but not difficult proofs, $M[G]$ is shown to reflect the same list of axioms $\phi_1, \ldots, \phi_k$, so that the $CH^{M[G]}$ proof goes through. At the same time it is shown that there are at least $\aleph_2^{M[G]}$ subsets of $\aleph_0$ in $M[G]$, which is a contradiction. This part of the proof utilizes an elegant argument involving combinatorics of infinite sets. The conclusion is that there cannot be a proof $\phi_1, \ldots, \phi_n = \mathrm{CH}$.

\section{Preliminaries---Set Theory, Ordinals, Cardinals and Models}

This section covers concepts that are prerequisite to understanding the proof of independence of $\mathrm{CH}$ in a manner meant to be intuitive and assuming minimal background: (1) ZFC, which axiomatizes the notion of a set, a versatile concept upon which all familiar mathematical structures can be built; (2) ordinals and cardinals, notions that enable ``counting'' beyond $\mathbb{N}$ and comparing different sizes of infinity; (3) models and their ability to reflect the axioms of set theory.

\subsection{Axiomatic Set Theory}

Sets are a versatile concept that can represent any familiar mathematical structure. For example, numbers can be represented as sets by letting each number be the set of all smaller numbers, with $0$ being the set $\emptyset = \{\}$. Ordered pairs $(x, y)$ can be represented as sets $\{ \{x\}, \{x,y\}\}$ and sequences $(x_1, x_2, \ldots, x_k)$ as nested pairs $(x_1, (x_2, (\ldots, x_k)\ldots))$. Functions $f: X \to Y$ can be represented as sets of ordered pairs $(x,y)$ with $x\in X$, $y\in Y$ such that for every $x$ in $X$ there is precisely one $(x,y)$ in $f$; $k$-ary relations over $X$ can be represented as subsets of $X^k$. Because of the versatility of sets, set theory has been the go-to notion as a foundation of mathematics.

However, defining what a set is is tricky. In 1901, Bertrand Russell posed his famous paradox, ``the set of all sets that are not members of themselves'' (and the ``barber who shaves all those, and those only, who do not shave themselves''). Call this set $x$. Is $x$ a member of itself? A moment's reflection reveals that $x$ is a member of itself if and only if $x$ is not a member of itself, a contradiction. Clearly, the notion of what a  set is cannot be unconstrained.

To remedy such paradoxes, in 1908 Earnst Zermelo proposed the first axiomatic set theory. The axioms laid out constrained the notion of a set enough so as to avoid self-referential paradoxes like Russell's barber. However, as was discovered later by Abraham Fraenkel and others, the axioms were not powerful enough to yield certain sets whose existence was commonly assumed by mathematicians. After additions of axioms by Fraenkel, John von Neumann and others, the Zermelo-Fraenkel theory (ZF) with the Axiom of Choice (AC), commonly known as ZFC, was settled on. ZFC is  designed to be powerful enough to prove the existence of all kinds of sets commonly assumed by mathematicians, but to refrain from constructs that lead to contradiction.\footnote{Many alternative axiomatizations exist. The von Neumann-Bernays-G\"{o}del set theory (NBG) is an extension of ZFC that allows classes in addition to sets. Classes are collections of elements that are defined by a property---a formula---but which may be too big to form a set. Morse-Kelley set theory goes further, allowing classes to be defined by a property that has quantifiers that range over classes. The topic is rich and extensive, and beyond the scope here. For a review see the Stanford Encyclopedia of Philosophy (https://plato.stanford.edu/entries/settheory-alternative/).}

In this section we will go over the formal axioms of ZFC, their meaning, and how they can be utilized to express familiar mathematical concepts. But first a few words on formal axiomatic systems.

\subsubsection{Formal languages and proofs}

Mathematicians typically write proofs in informal English interleaved with mathematical formulas. Informal proofs are understood to be shorthand for formal proofs that can be conducted entirely within a formal axiomatic system. Every proof is in principle machine-verifiable, otherwise it is not a proof.

A formal axiomatic system consists of a \textbf{language} of symbols and syntactic rules with which two types of objects can be composed: \textbf{terms}, which refer to objects in the domain and can be thought of as words of the language,  and \textbf{formulas}, which are mathematical assertions and can be thought of as phrases of the language. Terms and formulas may contain free (input) variables, so their value can vary: terms with no free variables evaluate to a single item in the domain and formulas with no free variables---the \textbf{sentences} (also known as \textit{statements})---evaluate to true or false.\footnote{For example, with $\mathbb{R}$ as the domain  $+, \cdot$ as binary operators and $>, <, =$ as binary predicates, $2.718$, $2\pi r$, and $(x+y)^2$ are examples of terms; $1+1 = 2$, $x > y$, $\forall x \: \forall y \: (x+y)^2 = x^2 + 2xy+ y^2, x < x+1$ are examples of formulas. The first and third are true sentences, the second is a formula whose truth depends on $x$ and $y$, and the fourth is a true formula with free variable x.} A finite collection of \textbf{axioms} or axiom schemas (templates that capture a family of axioms) specify all the formulas that can be taken for granted. 

\textbf{Logical axioms and inference rules.} Standard first-order logical axioms  accompany the domain-specific axioms of a formal language. Many different axiomatizations of first-order logic are essentially equivalent.\footnote{The following one is from the exposition of forcing by Weaver (2014):
\begin{enumerate}
\item $\phi \rightarrow (\psi \rightarrow \phi)$
\item $(\phi \rightarrow(\psi \rightarrow \chi)) \rightarrow ((\phi \rightarrow \psi) \rightarrow (\phi \rightarrow \chi))$ 
\item $(\lnot \phi \rightarrow \lnot \psi) \rightarrow (\psi \rightarrow \phi)$
\item ($\forall x(\phi \rightarrow \psi)) \rightarrow (\phi \rightarrow \forall x \psi)$, where $\phi$ is a formula that does not contain any free occurrences of $x$.
\item $\forall x (\phi(x)) \rightarrow \phi(t)$ where t is a term that does not share any variables with $\phi$. 
\end{enumerate}} 
Then, \textbf{logical inference rules} specify how to derive a new formula from previous formulas. Again, there are many equivalent sets of inference rules, although some care needs to be taken that the inference rules interact well with the logical axioms and no contradictions are created.\footnote{The following three inference rules are in Weaver (2014): (1) from $\phi$ and $\phi \rightarrow \psi$ we can derive $\psi$; this principle is known as \textit{modus ponens}. (2) From $\phi$ we can derive $\forall x \phi$ for any variable $x$, whether it appears in $\phi$ or not; this is \textit{universal generalization}. In some systems it needs to be restricted to variables not free in $\phi$, depending on how it combines with other axioms and inference rules to potentially lead to contradictions. (3) From $\phi$ we can derive any $\phi'$ that constitutes a renaming of the variables of $\phi$. A renaming $\rho$ of variables is one where $\rho(x) \not = \rho(y)$ iff $x$ and $y$ are different variables. This avoids disasters like renaming $x < y \rightarrow \exists z \: (x+z = y)$ to $a < b \rightarrow \exists b \: (a + b = b$). } All such inference systems are complete in making all logically valid inferences (in one or more steps), a fact that we will take for granted.\footnote{For example, if we know $\phi$ and $\psi$ we can derive $\phi \land \psi$, if we know $\phi \rightarrow (\psi \lor \chi)$ then $\phi\rightarrow \psi$, and $\forall x \: \phi(x)$ implies $\phi(t)$ for any term $t$ that does not contain variables appearing in $\phi$.}

\textbf{Proofs within a formal system.} A \textbf{proof} of a formula $\phi_k$ is simply a list of formulas, $\phi_1, ..., \phi_k$ such that each $\phi_i$ either is an axiom or is logically derived from previous formulas in the sequence. A formula $\phi$ is a \textbf{theorem} of the language $\mathcal{L}$, denoted $\mathcal{L} \vdash \phi$ when there is a proof $\phi_1, ..., \phi_k = \phi$ of $\phi$. The validity of a proof can be syntactically checked by a mechanical procedure. Equivalently, the language and the proofs can be \textbf{arithmetized}, a procedure devised by G\"{o}del that turns the question of validity of proofs into an arithmetic statement: ``$\phi_1, ..., \phi_k$ is a valid proof of $\phi_k$'' becomes ``the number $\lceil (\phi_1, ..., \phi_k)\rceil$ is an encoding of a syntactically valid sequence of formulas that constitute a valid proof of the formula encoded by the number $\lceil \phi_k \rceil$''. Here, $\lceil s \rceil$ is meant to be a function that turns a string $s$ into a number. Readers not familiar with G\"{o}del's arithmetization and the incompleteness theorems may want to read on the topic and understand the theorem statement and scope as well as the proof idea before reading this manuscript (Kim B.\footnote{ https://web.yonsei.ac.kr/bkim/goedel.pdf}; Smith 2013; Nagel and Newman 2001; Stanford Encyclopedia of Philosophy\footnote{https://plato.stanford.edu/entries/goedel-incompleteness/}; Buildt 2014; Batzoglou 2021). $\mathrm{CH}$ is a prime example of the incompleteness phenomenon. Moreover, expression of recursive functions in Peano Arithmetic, which is a key component of G\"{o}del's incompleteness theorem's proof, has great similarity to transfinite recursive definitions at the heart of ZFC. Such definitions are used heavily in the subsequent sections, including defining the universe of sets $V$, the class of names $N$, and the forcing-equality relation $Force^=$. Studying the proof of  G\"{o}del's first incompleteness theorem is a good warm-up and a soft prerequisite for reading this manuscript.

\subsubsection{The axioms of ZFC}

In the subsequent exposition we will assume that every item in the mathematical world is a set. This is not a necessary assumption in ZFC: for instance we could allow the existence of atoms (apple, orange, banana, Agamemnon, Julius Caesar, Cleopatra) that are not sets but can belong to a set. This generalization adds no power because all such atoms can be tokenized with numbers IDs and numbers can be represented as sets. So we will simplify and assume everything is a set.

The language of ZFC consists of the logical symbols $\land$, $\lnot$, and $\forall$ from which additional familiar symbols $\lor, \rightarrow, \leftrightarrow, \exists, \exists!$, can be defined as shorthand\footnote{For example $\phi \lor \psi$ is $\lnot(\lnot \phi \land \lnot \psi)$, $\exists x \phi$ is $\lnot \forall x \lnot \phi$, and $\exists ! x \: \phi$ is $\exists x (\phi \land \forall y (\phi \rightarrow y = x))$.}, the binary predicates $=$ and $\in$ from which other predicates such as $\not \in, \neq, \subset, \subseteq, \subsetneq$, connectors such as $\cap, \cup$ and abbreviations $\forall x \in y, \exists x \in y$ can be defined\footnote{For example, $x \subseteq y$ is defined as $\forall z \: z\in x \rightarrow z \in y$, $z = x \cap y$ is shorthand for $\forall w \:(w \in z \leftrightarrow w\in x \land w \in y)$, $\forall x \in y \: \phi$ stands for $\forall x \: x\in y \rightarrow \phi$ and $\exists x \in y\: \phi$ stands for $\exists x \: x\in y \land \phi$.}, and a countably infinite collection of variables $x_0, x_1, x_2, \ldots$.When we use ``$x, y, z, ...$'' within formulas, those are meant to represent one of the valid variables $x_i$ of the language. Terms of ZFC are very simple: they are just the variables. Formulas are defined recursively to be either a ground formulas $x \in y$ and $x=y$ where $x, y$ are terms, or composites $\phi \land \psi$, $\lnot \phi$, $\forall x \: \phi$ and so on, where $\phi, \psi$ are formulas and $x$ is a variable. 

The axioms of ZFC include: (1) the standard first order logical axioms; (2) two axioms for equality: $x = y \rightarrow (x \in z \rightarrow y \in z)$ and  $x = y \rightarrow (z \in x \rightarrow z \in y)$; (3) the domain-specific axioms presented below. 

There are many equivalent formulations of ZFC. The following one is similar to the one by Kunen (1980). Notice that most of the axioms postulate the existence of a certain set. The aim is to make the notion of a set as expansive as possible while avoiding  self-referential contradictions.

\begin{enumerate}
\item \textbf{Extensionality.} Two sets are equal if and only if they have the same elements. $\forall x \: \forall y \: (\forall z \: (z\in x \leftrightarrow z \in y) \leftrightarrow x = y)$. 
\item \textbf{Pairing.} For any two sets, there exists a set that contains them as members. $\forall x \: \forall y \: \exists z \: (x \in z \land y \in z)$.
\item \textbf{Union.} For any set, the union of all its member sets is a set. $\forall x \: \exists u \: \forall y \: (y \in u \leftrightarrow (\exists z \: (y \in z \land z \in x)))$. 
\item \textbf{Power.} For any set, there is a set of all its subsets. $\forall x \: \exists y \: \forall z \: (z \subseteq x \leftrightarrow z \in y)$. 
\item \textbf{Infinity.} There is a set with infinitely many members. $\exists x \: (\exists e \: (\forall y \: y \not \in e) \land e \in x \land (\forall z \in x \: z \cup \{z\} \in x))$.
\item \textbf{Foundation.} Any nonempty set has a member with which it shares no members. $\forall x \: (x \not = \emptyset \rightarrow (\exists y \: (y \in x \land y \cap x = \emptyset)))$. 
\item \textbf{Schema of Comprehension.} Let $\phi(x, y, w_1, ..., w_k)$ be a formula. There is a set that contains all elements $x \in y$ that satisfy $\phi$, for any choices of set $y$ and parameters (sets) $w_i$. $\forall y \: \forall w_1 \cdots \forall w_k \: \exists z \: (\forall x \: x \in z  \leftrightarrow (x \in y \land \phi(x, y, w_1, ..., w_k)))$. 
\item \textbf{Schema of Replacement.} The image of a set under a class function falls within a set. Let $\phi(x, y, A)$ be any formula with free parameters $x, y, A$, such that for every $x \in A$ there is a unique $y$ making $\phi$ true. The formula $\phi$ can be thought of as a function that takes $A$ as input and specifies some set $B$ replacing all elements $x$ of $A$ with the $y$ that makes $\phi(x,y,A)$ true. However, $\phi$ is not always a function because the domain over which $A$ ranges is not specified. Instead, $\phi$ is a \textit{class function}, and the reason for the terminology will become clear in Section 3.3. For now, it suffices to consider $\phi$ as any formula of the requisite property. Then, the Replacement Schema asserts that $B$ is a set that contains all the elements $y$ that are the image of set $A$ under $\phi$. $\forall A \:[(\forall x\in A\: \exists ! y \: \phi(x, y, A)) \rightarrow \exists B \: (\forall x \in A \: \exists y \in B \: \phi(x, y, A)]$. 
\item \textbf{Axiom of Choice (AC).} Given a set of nonempty sets, there exists a function that picks one element from each member set. To write down the axiom, it is convenient to use some shorthand: $\mathrm{Fun}(f)$ for the condition of $f$ being a function, $\mathrm{Dom}(f), \mathrm{Ran}(f)$ for the domain and range of a function, respectively,  $\bigcup x$ for the union of all sets in $x$, and $f(y)$ assuming $f$ is a function, for the unique $z$ s.t. $(y,z) \in f$. $\forall x \: (\emptyset \not \in x \rightarrow \exists f \: (\mathrm{Fun}(f) \land \mathrm{Dom}(f) = x \land \mathrm{Ran}(f) = \bigcup x \land \forall y \in x \: f(y) \in y)$.
\end{enumerate}

Some formulations of ZFC introduce another axiom, that of the \textbf{Existence of the Empty Set}: $\exists e \: \forall w \: w \not \in e$ or something equivalent. This is convenient but not needed because the existence of the empty set can be derived from Infinity and Comprehension. The proof sketch is included here as an example of how ZFC ``works''.

\begin{quote}
\textbf{Existence of $\emptyset$---proof sketch.} Let us write Infinity in shorthand as $\exists x \: \mathrm{Inf}(x)$, and Comprehension as $\forall y \: \mathrm{Comp}(\phi, y)$ noting that this is a separate axiom for every formula $\phi$ that contains no free occurrences of variable $y$. Let $\phi(w) := w \not = w$. A ZFC proof is just a list of formulas. For this proof, let the first two formulas be $\exists x \: \mathrm{Inf}(x), \forall y \: \mathrm{Comp}(w \not  =w, y)$. $A$ and $B$ implies $A\land B$ so let the third formula of the proof be  $\exists x \: \mathrm{Inf}(x) \land \forall y \: \mathrm{Comp}(w \not = w, y)$. Because $x$ does not appear to the right of $\land$, this implies $\exists x \: (\mathrm{Inf}(x) \land \forall y \: \mathrm{Comp}(w \not = w, y))$. Instantiating the $\forall y$ quantifier to $x$ yields $\exists x \: (\mathrm{Inf}(x) \land \mathrm{Comp}(w \not = w, x))$. Let's unroll the ``$\mathrm{Comp}$'' shorthand: $\exists x \: \exists z \: (\mathrm{Inf}(x) \land \forall w \: (w \in z \leftrightarrow (w \in x \land w \not = w)))$ (here, we also moved $\exists z$ outside the parentheses as $z$ does not appear within $\mathrm{Inf}(x)$). Because $A \land B$ implies $B$, we now get $\exists x \:\exists z  \:\forall w \:(w \in z \leftrightarrow (w \in x \land w \not = w))$\footnote{Perhaps in more than one derivation step, depending on the precise set of deduction rules.}. Now $x$ appears nowhere in the formula, so we can drop it: $\exists z \: \forall w \: (w \in z \rightarrow w \not = w)$ which is easily shown to yield $\exists z \: \forall w \: w \not \in z$. $\square$
\end{quote}

We can shorthand this to $\exists z \: \mathrm{Empty}(z)$ or $\exists z \: z = \emptyset$. Whenever $\emptyset$ is used as a shorthand, it is understood that a proof such as the above precedes it. All sets are compositions of the empty set, as will become clear. 

ZFC can be thought of as a rather low-level programming language. The syntax is kept to a minimum in order to reason about it precisely. Actually writing formal proofs in ZFC is extremely tedious but not needed in practice. Some remarks on the axioms:

\textbf{Extensionality} simply says that two sets are equal iff they have the same elements. That simplifies discourse. More generally, one could allow multiple distinct entities having the same elements, and then define a formula $\mathrm{EqualSets}(x,y)$ to be true whenever $x, y$ have the same elements. The resulting theory would be usually more cumbersome, sometimes more convenient, and overall equivalent.

Notice that \textbf{Pairing} allows other elements beyond $x,y$ to be in $z$. That is fine because using Comprehension in conjunction with Pairing one can define precisely the pair $\{x,y\}$.

The \textbf{Power} axiom is powerful: it postulates the existence of the set of subsets of a set. This allows jumping from $\mathbb{N}$ to continuum, and from there to much much larger infinities. It has been said that Power is way too powerful, and arguably it is the source of all trouble with infinities. The power set of a set $x$ will be denoted by $\mathcal{P}(x)$ or $2^x$.

\textbf{Comprehension} and \textbf{Replacement} require some discussion. These are not individual axioms, but rather axiom schemas. For any formula $\phi$ of the language, a separate axiom is instantiated. They correspond to mathematician's statements such as ``let $X$ be the set of all elements $x$ of $\mathbb{R}$ that satisfy property $P(x)$'' (Comprehension) and ``let $Y$ be the image of $X$ under function $f$'' (Replacement). First-order logic does not allow variables ranging over properties, and this is why axiom schemas are needed. Ranging over properties is allowed in second-order logic and finite second-order axiomatizations of set theory are possible. Discussing the merits of first order logic and (mostly de-)merits of second-order logic is beyond the scope here. It suffices to understand that there is a separate Comprehension and a separate Replacement axiom for every formula $\phi$. Great care needs to be taken to avoid creating self-referential contradictions such as a ``barber who shaves all those and those only who do not shave themselves''. This is why Comprehension creates a set of all elements $x$ \textbf{within a set $y$} that satisfy $\phi$. This simple condition prohibits circularity: assuming $y$ is a set and involves no circularity, $x$ involves no circularity. Similarly, Replacement creates a set out of the image of \textbf{an existing set}.

What about \textbf{Infinity}?  Infinity postulates a set that contains $\emptyset$ and is closed under $x \mapsto \{x\}$. Given a set $J$ with this property as given by Infinity, we can then use Comprehension to take the intersection of all subsets of $J$ that also satisfy this property, and this intersection will be isomorphic to $\omega$. Without Infinity it would not be possible to prove that infinite sets exist.

\textbf{Foundation} is mostly cosmetic. It ensures that there are no downward chains of membership $\cdots \in x_k \in x_{k-1} \in \cdots \in x_2 \in x_1 \in x_0$. It grants no additional power to the theory, but makes certain proofs way easier.

Finally, the \textbf{Axiom of Choice} (AC) is a non-obvious axiom on which enormous literature exists. As it turns out, AC is itself independent of the other axioms of ZF, and a proof can be provided  based on the forcing technique described in this manuscript. In fact there are axioms such as Determinacy that are contradictory with AC and could be appended to ZF as alternatives that lead to different mathematical results. Mathematicians often pay attention to the use of AC in proving a theorem.

\subsection{Levels of Infinity: Ordinals and Cardinals}

A set is \textit{well ordered} if its elements are ordered and, importantly, if every nonempty subset contains a smallest element:

\textbf{Definition: well-ordered set.} A set $S$ is well ordered by a binary relation $<$ if: (1) for any $x,y \in S$, either $x<y$ or $y<x$ or $x = y$; (2) if $x<y$ and $y<z$ then $x<z$; (3) every subset of $S$ contains a smallest element.

The natural numbers, $\mathbb{N}$ are well ordered. However, the real numbers $\mathbb{R}$ are \textit{not} well ordered by the natural relation $<$ because there are subsets such as $(0,1)$ that do not contain a smallest element. Requiring that every subset contains a smallest element (an \textit{infimum}) is a powerful condition. Every set can be well ordered as can be shown with the axiom of choice. Conversely, if every set can be well ordered then it is easy to prove the axiom of choice by defining a function that picks the infimum of each set in a collection of sets. 

Given a set $S$ well ordered by $<$, the order type of $S$ is the graph induced by $<$. As an example, let $S = \mathbb{N} \cup \{\mathrm{Julius}, \mathrm{Augustus}, \mathrm{Cleopatra}\}$ and let $n < \mathrm{Julius} < \mathrm{Augustus} < \mathrm{Cleopatra}$ for all $n \in \mathbb{N}$. The resulting graph is $0<1<2<...<\mathrm{Julius}<\mathrm{Augustus} < \mathrm{Cleopatra}$. This is a different order type from the graph $0<1<2<...$ even though both sets are countably infinite. 

As another example, what is the order type of $\mathbb{N}^2$? This depends on the induced order. For example, we can order $\mathbb{N}^2$ lexicographically: $(a,b) < (c,d)$ iff $a<c \lor (a = c \land b\leq d)$. Then, the order type is $(0,0)<(0,1)< \cdots < (1,0)<(1,1) < \cdots < ... < (n, 0) < (n,1)< \cdots$. If instead we order $\mathbb{N}^2$ by $(a,b) < (c,d)$ iff $ a+b < c+d \lor (a+b = c+d \land a\leq b)$ then the order type is same as that of $\mathbb{N}$, as a moment's reflection will reveal.

Ordinals are a way to express the notion of order type. They can be thought of as a way of counting past infinity. The cleanest definition of ordinals is by von Neumann:

\begin{quote}
A set $S$ is \textbf{transitive} if every member $x \in S$ is a subset $x \subset S$. A set $\alpha$ is an \textbf{ordinal} if it is transitive and well-ordered by $\in$.
\end{quote}

This definition leads to a way of representing numbers in $\mathbb{N}$ as ordinals. Each number is the set of smaller numbers. The ordinal $\mathbb{N}$ is usually denoted by $\omega$. The same scheme allows ``counting'' past $\mathbb{N}$.

\begin{itemize}
\item $0 := \emptyset$; $1 := \{ 0\} = \{ \emptyset \}$; $2 := \{0,1 \} = \{\emptyset, \{ \emptyset \}\}$; $3 := \{0,1,2\} = \{ \emptyset, \{ \emptyset \}, \{\emptyset, \{ \emptyset \}\} \}$,...
\item $\omega$ := $\{0,1,2,...\}$; $\omega + 1 := \omega \cup \{ \omega \}$, ..., $\omega\cdot 2 = \{0, 1, ..., \omega, \omega+1, \omega+2, ...\}$,... 
\item $\omega^2 = \{0, 1, ..., \omega, ..., \omega\cdot 2, ..., \omega\cdot 3, ...\}$
\item ...
\end{itemize}

As can be seen above, ordinals larger than $\emptyset$ come in two flavors: (1) \textbf{successor ordinals} $\alpha +1 := \alpha \cup \{ \alpha \}$; (2) \textbf{limit ordinals} such as $\omega = \bigcup_{\alpha < \omega} \alpha$, $\omega\cdot 2$, $\omega^2$ and so on.

\textbf{Ordinal addition: an example of transfinite recursion.} The addition $\alpha + \beta$ of two ordinals can be defined with \textbf{transfinite recursion}, a general technique for writing down a formula applicable to a well-ordered class of objects. To define a function $f$ for all ordinals, it suffices to define it for the base case $f(\emptyset)$, and then define $f(\alpha)$ assuming $f(\beta)$ is defined for any $\beta < \alpha$. Two cases need to be distinguished: the successor case $\alpha = \gamma+1$ for some $\gamma$, and the limit case $\alpha = \bigcup_{\gamma<\alpha} \gamma$. For ordinal addition, the definition of $\alpha + \beta$ uses transfinite recursion on $\beta$:

\begin{itemize}
\item Base case, $\beta = 0\:$: then, $\alpha + 0 := \alpha$; 
\item Successor case,  $\beta = \gamma+1\:$: then, $\alpha + \beta := (\alpha + \gamma) + 1$; 
\item Limit case, $\beta$ is a limit ordinal: then, $\alpha + \beta := \bigcup_{\gamma < \beta} \alpha + \gamma$. 
\end{itemize}

In the above definition, note that $\alpha, \beta$ can be infinite and even uncountable. When seeing transfinite recursion for the first time, a natural question is how does it ``work'' in ZFC. The way it is implemented is by a formula, in this case $\mathrm{Plus}(x, y, z)$, which is true iff $x$ and $y$ are ordinals and $z$ is their sum. It is almost like coding. Let's define $\mathrm{Plus}(x,y,z)$ in full detail to exemplify how this is done.

Let $\mathrm{Ord}(x)$ be shorthand for the formula that asserts $x$ is an ordinal,  $\mathrm{Succ}(x, y)$ for $y = x \cup \{x\}$, and $\mathrm{LimOrd}(x)$ for $x$ being a limit ordinal. $\mathrm{Plus}(x, y, z)$ has to ensure $x$ and $y$ are ordinals and $z$ is their ordinal sum. This is accomplished with a helper sequence $s$ of length $y+1$ (a function whose domain is $y+1$\footnote{Recall that $y+1$ is the set of ordinals $\leq y$.}) such that at every ``position'' $w, 0 \leq w \leq y$ the value $s[w]$ is $x + w$. Let $\mathrm{Fun}(x)$ be shorthand for $x$ being a function\footnote{A set of pairs $(a, b)$ such that for any $a$ there is at most a single $(a, b)$.} and $\mathrm{Dom}(x, y)$ for ``the domain of $x$ is $y$''. Also, let the shorthand $s[i] = y$ stand for a formula $\mathrm{Eval}(s, i, y)$ that is true iff $s(i) = y$, i.e., if the pair $(i, y)$ is a member of $s$. The reader may want to write down their $\mathrm{Plus}$ version before reading the one below. For readability, $\alpha, \beta, ...$ are used as variables taking ordinal values:

$$\mathrm{Plus}(\alpha, \beta, \gamma) := \mathrm{Ord}(\alpha) \land \mathrm{Ord}(\beta) \land \exists \beta' (\mathrm{Succ}(\beta, \beta') \land \exists s \: (\mathrm{Fun}(s) \land \mathrm{Dom}(s, \beta') \land $$
$$\forall \delta \in \beta' \: ((\delta = \emptyset \land s[\delta] = \alpha) \lor (\exists \zeta \: \mathrm{Succ}(\zeta, \delta) \land \mathrm{Succ}(s[\zeta], s[\delta])) \lor (\mathrm{LimOrd}(\delta) \land s[\delta] = \bigcup_{\zeta \in \delta} s[\zeta]))  \land \gamma = s[\beta])$$

\begin{quote}
Note the use of shorthand $\bigcup_{\zeta \in \delta} s[\zeta]$. The range of $s$ for $\zeta \in \delta$ exists by Replacement and its union exists by Union. Composing a formula for $\bigcup$ is an exercise left for the reader.
\end{quote}

Readers familiar with G\"{o}del's incompleteness theorem and the expression of recursive functions with Peano Arithmetic formulas will notice a similarity here, where intermediate values are ``saved'' in $s$. The upshot is that \textbf{every transfinite recursive definition can be expressed as a formula}. One way to do so is the technique given above, where a helper sequence $s$ stores all intermediate values. From now on, we will liberally make transfinite recursive definitions without the need to provide formulas.

Notice that $1+\omega = \omega$ while $\omega + 1 = \{0,1,...,\omega\} \neq \omega$. As another example, $\omega + \omega$ is a limit ordinal equal to $\bigcup_{\alpha < \omega} \omega + \alpha$ and it has the same order type as $0<1<2<\cdots < 0'<1'<2'<\cdots$. Figure 1 displays ordinals up to $\omega^{\omega}$.

\begin{figure}
\centering
\includegraphics[width=300pt]{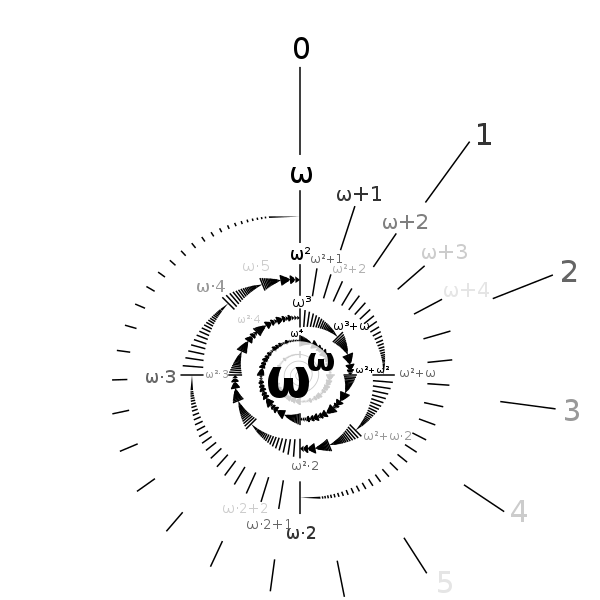}
\caption{\textbf{Representation of ordinal numbers up to $\omega^{\omega}$.}} Picture found in Wikipedia (https://en.wikipedia.org/wiki/Ordinal\_number) without attribution. 
\end{figure}

\textbf{Ordinal multiplication.} The ordinal $\alpha \cdot \beta$ can be defined as the Cartesian product of $\alpha$ and $\beta$ with elements sorted lexicographically. By convention, the \textit{least} significant position is placed first. As an example, $\omega \cdot 2 = \{(0,0), (0,1), (0,2), \ldots, (1,0), (1,1), (1,2), \ldots \} = \omega + \omega$ whereas $2 \cdot \omega = \{(0,0), (0,1), (1,0), (1,1), (2,0), (2,1), \ldots \} = \omega$. As another example, $\omega^2$ has the same order as $\mathbb{N}^2$ under lexicographic order $(a, b) < (c, d) \leftrightarrow a < c \lor (a = c \land b<d)$. 

No finite ordinal can be mapped bijectively to a different ordinal. However, this is no longer true for infinite ordinals. For example, $\omega^2$ can be mapped bijectively to $\omega$ because both sets are countable. 

There are certain ordinals that cannot be put in a bijection with a smaller ordinal. All $n \in \omega$ are examples. Also, $\omega$ is an example. A \textbf{cardinal} is an ordinal that cannot be put in a bijection with a smaller ordinal. Every member of $\omega$ as well as $\omega$ itself are cardinals and $\omega$ is the smallest infinite cardinal. Let $\mathrm{Card}(\kappa)$ be shorthand for the formula asserting $\kappa$ is a cardinal. There is special notation for infinite cardinals: $\aleph_0 := \omega$ and every subsequent infinite cardinal is given an ordinal subscript: $\aleph_0 <\aleph_1 < \aleph_2 < \cdots < \aleph_k < \cdots < \aleph_{\omega} < \cdots$. Each of these cardinals is a different order of infinity, ``bigger'' than the previous infinity in the sense that there cannot be a bijection between them. The notation $\aleph_\alpha$ can be thought of as shorthand for a formula $\mathrm{Aleph}(\kappa, \alpha)$ asserting that $\alpha$ is an ordinal and $\kappa$ is precisely the $\alpha$th infinite cardinal. Given any set $S$, we denote by $|S|$ the unique cardinal that maps bijectively to $S$.

The \textbf{cardinal successor} of a cardinal $\kappa$ is the cardinal immediately larger than $\kappa$, denoted $\kappa^+$. This should not be confused with the ordinal successor of $\kappa$, which is $\kappa \cup \{\kappa\}$.

\textbf{Cardinal exponentiation.} Given infinite cardinals $\kappa$ and  $\lambda$, $\kappa^{\lambda}$ is defined to be the cardinality of the set of functions from $\lambda$ to $\kappa$. The usual properties of exponentiation such as $\kappa^{\lambda + \xi} = \kappa^{\lambda}\kappa^{\xi}$ and $\kappa^{\lambda^{\xi}} = \kappa^{\lambda \xi}$ follow easily from the definition. With the Axiom of Choice (AC), additional properties simplify cardinal exponentiation:
\begin{itemize}
\item For infinite $\kappa \leq \lambda$, we have $\kappa^{\lambda} = 2^{\lambda}$.
\item If there exists cardinal $\mu < \kappa$ with $\mu^{\lambda} \geq \kappa$, then $\kappa^{\lambda} = \mu^{\lambda}$.
\end{itemize}

\begin{quote}
\textbf{Confusing abuse of notation: ordinal vs. cardinal exponentiation.} There is an unfortunate notational discrepancy between ordinal and cardinal exponentiation.  Exponentiating to a limit ordinal such as  $\omega^{\omega}$ denotes taking the limit $\sup\{\omega^n | n \in \omega\}$, while if $\omega$ is viewed as the cardinal $\aleph_0$, $\aleph_0^{\aleph_0}$ denotes the cardinal of same size as the set of functions $\aleph_0 \to \aleph_0$. Therefore, $\omega^{\omega} = \bigcup_{n\in \omega}\omega^n$ is a countable union of countable sets and therefore countable, while $\aleph_0^{\aleph_0}$ is of cardinality continuum. To see why this is so, observe that $2^{\aleph_0} \leq \aleph_0^{\aleph_0} \leq (2^{\aleph_0})^{\aleph_0} = 2^{\aleph_0^2} = 2^{\aleph_0}$. 
\end{quote}

With this terminology, $\mathrm{CH}$ states that $2^{\aleph_0} = \aleph_1$. $\mathrm{GCH}$ states that $2^{\aleph_{\alpha}} = \aleph_{\alpha+1}$ for any ordinal $\alpha$.

\begin{quote}
\textbf{How big are countable ordinals?} Countable ordinals are way smaller than continuum. Since we are all familiar with the real line or interval $(0,1)$, which we don't consider all that ``big'', we may conclude that countable ordinals can't be too big. In fact, countable ordinals keep advancing to stages that are absurdly high. As a first flavor, Cantor's ordinal $\varepsilon_0$ is defined as the first ordinal that is a fixed point of $\varepsilon = \omega^{\varepsilon}$, and consequently $\varepsilon_0 = \omega^{\omega^{\omega^{\cdots}}} = \sup\{\omega, \omega^{\omega}, \omega^{\omega^{\omega}}, \ldots \}$. It is still countable. Next, $\varepsilon_1$ is the next fixed point, $\varepsilon_2$ the next, and so on: $\varepsilon_0, \varepsilon_1, \ldots, \varepsilon_{\omega}, \ldots, \varepsilon_{\varepsilon_0}, \ldots, \varepsilon_{\varepsilon_{\varepsilon_{\ldots}}}, \ldots$ are all fixed points of $\varepsilon = \omega^{\varepsilon}$ and each such $\varepsilon_{\zeta}$ is countable as long as $\zeta$ is countable. This topic is complex and beyond the scope here, but worth keeping in mind in the subsequent sections when wondering how a countable model $M$ could cogently reflect properties of uncountable cardinals.
\end{quote}

\subsubsection{Expressing mathematics in ZFC}
ZFC is powerful enough to express all ordinary mathematics. Some examples: 
\begin{itemize}
\item \underline{Definition of $\mathbb{N}$.} $\mathbb{N}$ is defined using Infinity and Comprehension: using abbreviations $\emptyset$ and $S(x) := \ x \cup \{x\}$, we can define $I$ as an infinite set satisfying the axiom of Infinity and $\mathbb{N}$ as $\{x \in I: \forall J \: \mathrm{Inf}(J) \rightarrow x \in J\}$, which exists by the axiom of Comprehension. Now we can define $0 := \emptyset, 1 := S(0)$ and addition inductively for any $m,n \in \mathbb{N}$ expressed as a formula $\phi_{+}(m,n,s)$ that is true iff $s = m+n$ by $0+0 := 0; m+S(n) = S(m)+n := S(m+n)$ and with some effort prove all the properties of addition such as associativity and commutativity. Similarly a formula $\phi_{\times}$ can be defined for multiplication with $m\times S(n)  := m\times n + n$. These recursive definitions cam be turned into formulas using precisely the same technique employed above to define ordinal addition. Then, all properties of these operations can be proven using ZFC axioms.
\item \underline{Definition of $\mathbb{Q}$.} $\mathbb{Q}$ can be defined in a multitude of ways; a clean one is representing $m/n$ as a pair $(m,n)$ with no common divisor\footnote{A formula to achieve that would be $\lnot \exists r,s,t\in \mathbb{N} \:  (r > 1 \land s>1 \land t>1 \land \phi_{\times}(r,s,m) \land \phi_{\times}(r,t,n)$.} and defining addition $(a,b) + (c,d) := (r, s)$ whereby $r/s = (ad+bc)/cd$ but $r,s$ share no common divisor.
\item \underline{Definition of $\mathbb{R}$.} $\mathbb{R}$ can also be defined in many ways, such as using \textit{Dedekind cuts}: a cut is a partition of $\mathbb{Q}$ into two subsets $A,B$ such that $A$ is nonempty and a proper subset of $\mathbb{Q}$, if $x<y \in A$ then $x \in A$, and additionally  $\forall x \in A \: \exists y \in A \: y>x$. $\mathbb{R}$ can be defined as the the set of all cuts $(A,B)$ with '$<$' defined by $(A,B) < (C,D)$ iff $A$ is a proper subset of $C$. For example, $\sqrt{2}$ becomes the cut $(\{x\in \mathbb{Q} : x<0 \lor x^2 < 2 \}, \{y \in \mathbb{Q} : y>0 \land y^2 > 2\})$. Addition of $(A,B), (C, D)$ can be defined as $(A', B')$ with $A' = \bigcup \{x+y : x \in A \land y \in C\}$ and $B' = \mathbb{Q} - A'$. And so on. All familiar properties of $\mathbb{R}$ can be proven. 
\end{itemize}

With some pain, all familiar concepts in mathematics can be defined, and all mathematical proofs can be expanded to formal proofs within ZFC. ZFC can be thought of as a low-level programming language for mathematics. 

\subsection{Sets, Classes and the Universe}

The ZFC axioms postulate the existence of sets built out of other sets. However, not everything expressed syntactically in ZFC is a set. The famous Russell's paradox of a set of all sets that are not members of themselves is avoided here by simply not allowing $x \in x$. Specifically, $\{x | x \in x\} = \emptyset$ and $\{x | x \not \in x\}$ is the collection of all sets, which is not a set in ZFC.

\textbf{Proper classes.} A formula $P(x)$ of a single free variable $x$ (a property $P$), sometimes defines a set $\{x \:|\: P(x)\}$ of all elements satisfying the property. For example, the collection of all $x$ that satisfy $x \not = x$ is $\emptyset$, which is a set. The collection of $x$ that satisfy $\mathrm{Inf}(x) \land \forall z \: \mathrm{Inf}(z) \rightarrow x \subseteq z$ is precisely $\{\omega\}$. However, $P(x)$ need not define a set. Sometimes the set status of $x : P(x)$ is unknown, and sometimes it is \textit{known} not to be a set, in which case $x :  P(x)$ is called a \textbf{proper class}. For example, $x : \mathrm{Ord}(x)$ cannot be a set: if it were, it would be (1) transitive because given any $\alpha$ such that $\mathrm{Ord}(\alpha)$, all $\alpha$'s members $\beta\in \alpha$ are also ordinals and therefore $\mathrm{Ord}(\beta)$; and (2) well ordered by $\in$ because given any two ordinals $\alpha, \beta$ either $\alpha \in \beta$ or $\beta \in \alpha$ or $\alpha = \beta$, and given any subset there would be a smallest element because of Foundation\footnote{Given a subset $s$ consisting of ordinals, let $\alpha_0 \in s$. If $\alpha$ is the smallest element, we are done. If not, pick $\alpha_1 <  \alpha_0$, i.e., $\alpha_1\in \alpha_0$. Continue on with $\alpha_{i+1} \in \alpha_i$ until the smallest element is hit. The sequence has to terminate because by Foundation there is no infinite descending chain $\cdots \in \alpha_{i+1}\in \alpha_i \cdots \alpha_0$.}. Therefore, if it were a set, $\mathrm{Ord}(x)$ would be an ordinal! But then $\mathrm{Ord}(x)$ and indeed $\mathrm{Ord}(x)+1$ would need to be ``in'' $\mathrm{Ord}(x)$, an absurdity. Therefore, the ordinals are a proper class and not a set. We will denote them by $ON$: when we say $x \in ON$ this is just another way of writing $\mathrm{Ord}(x)$.\footnote{Incidentally, there are countably many formulas $P(x$) and therefore there are countably many proper classes definable by a formula and countably many sets definable by a formula. For instance, when considering the continuum-many subsets of $\omega$, only countably many of them can ever be specified.}

\textbf{The universe $V$ and the cumulative hierarchy.} The most important proper class of sets is the entire universe $V$ of all sets. $V$ can be defined with transfinite recursion similarly to how $\mathrm{Plu}s(x,y,z)$ was defined in the preceding section.  The set $V_{\alpha}$ is defined for every ordinal $\alpha$. Note that $V_{\alpha}$ is indeed a set and is called the \textit{universe at stage} $\alpha$:
\begin{itemize}
\item $V_0 := \emptyset$
\item $V_{\alpha + 1} := \mathcal{P}(V_{\alpha})$, the power set of $V_{\alpha}$ for any successor ordinal $\alpha$.
\item $V_{\alpha} := \bigcup_{\beta < \alpha} V_{\beta}$ for any limit ordinal $\alpha$.
\item $V := \bigcup_{\alpha \in ON} V_{\alpha}$.
\end{itemize}

This definition stands for a ZFC formula, call it $V(\alpha, x)$ which says that ``$x$ is a member of $V_{\alpha}$'', and a formula $V(x)$ which says that ``$x$ is a member of some $V_{\alpha}$''. Before proceeding, the reader may want to write down versions of these formulas.

Let's write down these formulas. $Pow(x,y)$ is an abbreviation for the formula ``$y$ is the power set of $x$''. The other abbreviations have been introduced. Notice the use of $V_{aux}$, which is a helper sequence that contains all stage-$\beta$ universes for $0 \leq \beta \leq \alpha$.\footnote{$V_{aux}$ is a set that is a function with domain $\alpha+1$, similarly to $s$ in the definition of $\mathrm{Plus}$. $V_{aux}[\beta] = \gamma$ stands for $\mathrm{Eval}(V_{aux}, \beta, \gamma)$, or equivalently for $(\beta, \gamma) \in V_{aux}$.} Also note that it is easy to show that if $x\in V_{\alpha}$ then $x \in V_{\gamma}$ for all $\gamma > \alpha$. This is called a \textbf{cumulative} property, and $V$ is a \textbf{cumulative hierarchy}.
$$V(\alpha, x) := \mathrm{Ord}(\alpha) \land \exists V_{aux} \: \mathrm{Fun}(V_{aux}) \land \mathrm{Dom}(V_{aux}, \alpha+1) \land$$ $$\forall \beta \in \alpha + 1 \: (\beta = \emptyset \land V_{aux}[\beta] = \emptyset) \lor (\exists \gamma \in \beta \: \beta = \gamma + 1 \land \mathrm{Pow}(V_{aux}[\gamma], V_{aux}[\beta])) \lor$$
$$ (\mathrm{LimOrd}(\beta) \land V_{aux}[\beta] = \bigcup_{\gamma \in \beta} V_{aux}[\gamma])) \land x \in V_{aux}[\alpha]$$
We can define the proper class $V$ with the formula $V(x) := \exists \alpha \: V(\alpha, x)$. It is not hard to prove that $\forall x \: V(x)$. 

Having defined $V(\alpha, x)$, a useful property to define for any set $x$ is its \textbf{rank}. Let $\mathrm{rank}(x)$ be  the first ordinal $\alpha$ such that $x \in V_{\alpha}$, i.e., the least ordinal $\alpha$ such that $V(\alpha, x)$. It is easy to see that $\mathrm{rank}(x)$ is the smallest ordinal that exceeds the rank of any member of $x$. If $\beta$ is the supremum of ranks of elements in $x$, then $\mathrm{rank}(x) =\beta+1$ if there is $y\in x$ with $\mathrm{rank}(y) = \beta$, and otherwise $\mathrm{rank}(x) = \beta$.

\textbf{Class functions.} A class function is a unary operation defined on a class of sets. Consider a class of sets $A  :  P(A)$. A class function $\phi$ over the class is a unary operation that takes any $A$ such that $P(A)$ as argument and returns some set $B$. For example, the formula $\phi(x, y, A)$ in the Replacement Schema takes as input any $A$ in the universe $V$ and defines a set $B$ that consists of replacing every $x$ in $A$ with the unique $y$ such that $\phi(x, y, A)$. The Replacement Schema postulates that $B$ is a set.

\begin{quote}
\textbf{Summary: sets and classes.} A set is any entity proven to exist by ZFC. Syntactically, think of a set as a variable $x$ bound within a formula that states that $x$ exists and is the unique entity with so and so properties. A good example is $\omega$, which should be thought of as the bound variable $x$ within the formula $\exists x \: \mathrm{Inf}(x) \land \forall y \: \mathrm{Inf}(y) \rightarrow x \subseteq y$. A class is just a formula---a property---$P(x)$ of a single free variable. It defines any $x$ that makes the resulting sentence true. So, when we define a class such as $ON$ (the ordinals), we really mean anything that satisfies the formula $\mathrm{Ord}(x)$. And when we abuse notation and say $\alpha \in ON$, we really mean $\mathrm{Ord}(\alpha)$. A proper class is simply a class---a single-variable formula---which is satisfied by a collection of items that cannot possibly be a set. $ON$ and $V$ are examples of proper classes. When we use transfinite recursion to define a set, a class, a function or a relation, such as in the case of $V_{\alpha}$ or of ordinal addition $\mathrm{Plus}(\alpha, \beta, \gamma)$, we simply are writing down a \textit{single} ZFC formula in readable form, by exposing the base case, the successor ordinal case, and the limit ordinal case.
\end{quote}

\subsection{Models and Reflection}

Given a collection $\Phi$ of formulas of a language $\mathcal{L}$, a model $M$ is basically a set of elements in which the formulas are assigned meaning and are true. More precisely, a model $M$ is a set together with assignments of constants of $\mathcal{L}$ to elements of $M$, functions of $k$ variables of $\mathcal{L}$ to functions $M^k \to M$, and relations of $\mathcal{L}$ to subsets of $M\times M$. $M$ is a \textbf{model} of $\Phi$ if each formula in $\Phi$ is true under this assignment. 

As an example, consider again the axioms of group theory: (1) associativity: $\forall x \: \forall y \: \forall z \: (x \cdot y) \cdot z = x \cdot (y \cdot z)$; (2) identity: $\forall x \: e\cdot x = x \cdot e = x$; (3) inverse: $\forall x \: \exists x^{-1} \: x \cdot x^{-1} = x^{-1} \cdot x = e$. One model of these axioms is $M = \{0, 1, 2, 3, 4\}$ with $e := 0$, $x \cdot y := x+y$ mod $5$. Another model is $\mathbb{Z} = \{\ldots, -2, -1, 0, 1, 2, \ldots\}$ with $e := 0$, $x \cdot y := x+y$ and $x^{-1} := -x$. Yet another model is $SO(3)$, the group of rotations in $3$ dimensions: identity is the rotation by zero degrees in any axis; associativity holds between any three rotations; and each rotation has the reverse rotation as its unique inverse.

The notion of a model is useful in determining whether a given statement is the logical consequence of a collection of axioms. As an example, does the commutative law $\forall x \: \forall y \: x \cdot y = y \cdot x$ follow from the laws of group theory? It doesn't, because it holds for some models of group theory, such as $G$ and $\mathbb{Z}$, but doesn't hold for others such as $SO(3)$: performing a rotation by $90^o$ by the $z$-axis followed by one by the $y$-axis takes $(1,0,0)$ to $(0,0,1)$ while reversing the order results in $(0,1,0)$. An important theorem is behind this: 

\begin{quote}
\textbf{G\"{o}del's completeness theorem:} any consistent collection of formulas of first order logic has a model. 
\end{quote}

Consequently, if neither $\phi$ nor $\lnot \phi$ is a logical consequence of $\Phi$, then there exists models both for $\Phi \cup \{\phi\}$ and $\Phi \cup \{\lnot \phi\}$. The commutative law is an example of a formula independent of the axioms of group theory. Consequently, there are commutative models such as $\mathbb{Z}$ and non-commutative models such as $SO(3)$.

G\"{o}del's completeness theorem is outside the scope of this manuscript. The reader can take it as given, or preferably refer to an exposition of the Henkin version of the proof (Arjona and Alonso 2014). The proof is not hard, and is illuminating on how models arise from syntax.

\textbf{Models of ZFC.} In ZFC there are no constants (alternatively, sometimes $\emptyset$ is defined as a constant), there are countably many variables $x_0, x_1, x_2, \ldots$, no functions, and two relations: $\in, =$. A model of ZFC is a set $M$ together with assignments of the relations $\in$ and $=$ to subsets of $M \times M$ so that each of the countably many axioms of ZFC is logically valid under these assignments.

\textbf{Standard transitive models.} The language of ZFC devoid of a model is just syntax and $\in$ is just a letter that by itself carries no meaning. The axioms do give it structure, but in principle a model could be devised where $\in$ is a relation other than set membership. A \textbf{standard} model $M$ is one where $\in$ is set membership. There is a bit of meta reasoning here. Who decides what is set membership? We do! Examining ZFC as syntax, and thinking of  $M$ semantically, we decide that indeed $\in$ is set membership within $M$. Formally, if we assume the class $V$ to be the ``true'' set universe, a standard model $M$ is a set in $V$ (i.e., $\exists \alpha \: M \in V_{\alpha}$) and we check that $\in_M$, viewed as a subset of $M\times M$, is precisely $\in_V \cap \: M\times M$, where $\in_V$ is the subset of $V$ that stands for ``true'' $\in$. If in addition $M$ is a transitive set, then $M$ is a \textbf{standard transitive} model and has nice properties that will be covered below.

In the proof of independence of CH, the notion of a model of ZFC will be used extensively. The models will be standard countable transitive sets. 

Each model will be a mini-universe within which a sentence $\phi$ may be true or false, and may have a truth value different from its truth value across $V$. To keep track of truth of formulas within models, two concepts are key: \textbf{relativization} and \textbf{reflection}. Informally, a formula $\phi$ is \textit{relativized} to a set $M$ when $\phi$ talks about objects in $M$ only. For a trivial example, if $\phi := \forall x \: x = \emptyset$ then $\phi$ is false; however, relative to $M := \{\emptyset\}$ $\phi$ is true. As another example, $\psi := \forall x \exists y \: y = x+1$ is true relativized to $\mathbb{N}$ but false relativized to $\{0,1,\ldots, 2022\}$. We denote $\phi$ relativized to $M$ by $\phi^M$. A formula $\phi$ is \textit{reflected} in $M$ when $\phi$ is equivalent to its relativization: $\phi \leftrightarrow \phi^M$. Let's make these notions more precise:

\textbf{Relativization.} A formula $\phi$ is \textbf{relativized} to a set $M$ if all its variables range over $M$. This is achieved by a simple syntactic manipulation to make all existential and universal quantifiers range over $M$: all $\forall x \: \psi(x)$ and $\exists x \: \psi(x)$ are simply converted to $\forall x \in M \: \psi(x)$ and $\exists x \in M \: \psi(x)$. More precisely, any set $M$ can be thought of as a formula of a single variable $M(x)$ that defines a set of elements satisfying it. Then, $\forall x \in M \: \psi(x)$ and $\exists x \in M \: \psi(x)$ are simply shorthand for  $\forall x \:(M(x) \rightarrow \psi(x))$ and $\exists x \: (M(x) \land \psi(x))$. We denote this syntactically modified formula by $\phi^M$. Some examples:
\begin{itemize}
\item $\phi := \forall x \: x = \emptyset$ is false. However, if we let $M := \{\emptyset\}$, $\phi^M$ becomes $\forall x \in M \: x = \emptyset$, which is true.
\item $\phi := \forall x \: x = \emptyset \lor \exists y \: x = y \cup \{y\}$ is also false. However, relativized to $\omega$ it is true.
\item $\phi := \forall x \: x = \emptyset \lor \exists y \: y \in x$ is true. However, if we let $M := \{1\}$\footnote{Recall, $1 := \{\emptyset\}$.} then $\forall x \in M \: x = \emptyset \lor \exists y \in M \: y \in x$ is false.
\end{itemize}
 In general, formulas relativized to a set $M$ state something \textit{about} $M$. 
 
\textbf{Relativization to a class/property.} The notion of a set model can be generalized to a class model. Syntactically, a set in ZFC is just a formula $\phi(y)$ for which the sentence $\exists x \: \forall y \: (y\in x \leftrightarrow \phi(y))$ is derivable in ZFC. For example, $\emptyset$ is captured by $\exists x \: \forall y \: (y \in x \leftrightarrow y \not = y)$.\footnote{When a ZFC proof uses $\emptyset$, this should be understood as shorthand for the following: (1) append to the beginning of that proof the derivation of the sentence $\exists x \: \forall y \: (y \in x \leftrightarrow y \not = y)$; (2) understand whatever formula uses $\emptyset$ to be a $\land$ with $\exists x \: \forall w \: w \not \in x$ (where the $x$ variable is treated properly in the rest of the formula). For example, $1 = \{\emptyset\}$ should be understood as shorthand for something like $\exists w \: \exists x \: (\forall y \: (y \in x \leftrightarrow y \not = y)\land \forall z \: z \in w \rightarrow z = x)$. Writing down proofs in full mechanically verifiable detail is incredibly painful but it works. The same is true for assembly code, which is simpler to reason about and implement in hardware, but cumbersome to program on. In our case, $\emptyset$, $1$, $2$, and so on are examples of shorthand that unrolls into longer formulas and proofs.} More generally, a formula $P(y)$ defines the class of all $y$ that satisfy property $P$. This class becomes a set if it is provable in ZFC that $\exists x \: \forall y \: (y \in x \leftrightarrow P(y))$. Even when such a proof is not possible, however, we can relativize a formula $\phi$ to $P$: simply require that all variables of $\phi$ range over values for which $P$ is true. Syntactically, all subformulas of the form $\forall x \: \psi(x)$ and $\exists x \: \psi(x)$ are converted to $\forall x (P(x) \rightarrow \psi(x))$ and $\exists x (P(x) \land\psi(x))$, respectively. For example:
\begin{itemize}
\item $\forall x\: \forall y \: x \in y \lor x = y \lor y \in x$ is clearly false. However, we can relativize it to the proper class $ON$ of ordinals, which is \textit{not} a set. $\forall x \in ON \: \forall y \in ON \: x \in y \lor x = y \lor y \in x$ is true. This latter formula should be understood as shorthand for $\forall x \: (\mathrm{Ord}(x) \rightarrow \forall y \: (\mathrm{Ord}(y) \rightarrow x\in y \lor x = y \lor y \in x))$. 
\end{itemize}

\textbf{Reflection.} Reflection is the property of a set $M$ whereby a  formula $\phi$ is equivalent to $\phi^M$: we say that $M$ \textbf{reflects} $\phi$ just in case $\phi \leftrightarrow \phi^M$ is true whenever all free variables of $\phi$ range over $M$. The latter point can be a source of confusion. Say that $\phi(x,y)$ has free variables $x$ and $y$. ``$M$ reflects $\phi$'' means that whenever $a, b \in M$, we have $\phi(a,b) \leftrightarrow \phi^M(a,b)$. For example, let $\phi(x, y) := x \in y \lor y\in x \lor x = y$. Then, $\omega$ reflects $\phi$ trivially because $\phi = \phi^{\omega}$ and therefore for any $m, n \in \omega$, $\phi(m, n) \leftrightarrow \phi^{\omega}(m, n)$. But now let $\psi(x) := \forall y \: \phi(x, y)$. In this case, $\psi^{\omega} := \forall y \in \omega \: x \in y \lor y \in x \lor x = y$, and $\omega$ does not reflect $\psi$ because we can find $n \in \omega$ (in fact, any $n \in \omega$) for which $\psi(n)$ is false and $\psi^{\omega}(n)$ is true.

Recall the definition of transitivity. A set $M$ is \textbf{transitive} if every member's member is a member: $M$ is transitive if $\forall x \in M \: y\in x \rightarrow y \in M$. A transitive set $M$ reflects a large class of formulas: those that involve only looking at elements of $M$ and their contents. Some examples before this concept is made more precise.

Consider $\mathrm{Empty}(x) := \forall w\: w \in x  \rightarrow w \not = w$. Relativized to $M$, $\mathrm{Empty}^M(x) := \forall \: w \in M \: w \in x \rightarrow \: w \not = w$. Is $\mathrm{Empty}(x)$ reflected in $M$? Equivalently, does $\mathrm{Empty}(x) \leftrightarrow \mathrm{Empty}^M(x)$ hold whenever $x$ ranges over $M$? Well, it is easy to show this does not hold in general. Let $M := \{1, 2\}$ and let $x = 1$. By our definition of $1$ $:= \{0\}$, $\mathrm{Empty}(1)$ is false because $0 \in 1$ but $\mathrm{Empty}^M(1)$ is true because $\forall w \in M \: w \not \in 1$ holds. What if $M$ were transitive? Then, because $0 \in 1$ whenever $1 \in M$, it would necessarily be the case that $0 \in M$. It is easy to show that $\mathrm{Empty}(x)$ is reflected in any transitive set $M$.

As a more involved example, consider the definition of an ordinal as a transitive set totally ordered by $\in$:  $\mathrm{Ord}(x) := \forall y \: \forall z \: y \in x \rightarrow (z \in y \rightarrow z \in x) \land (z \in x \rightarrow (z \in y \lor y \in z \lor y = z))$. Let's relativize $\mathrm{Ord}(x)$ to $M$: $\mathrm{Ord}^M(x) := \forall y \in M \: \forall z \in M \: y \in x \rightarrow (z \in y \rightarrow z \in x) \land (z \in x \rightarrow (z \in y \lor y \in z \lor y = z))$. Let $M$ be transitive and let some ordinal $\alpha \in M$. $\mathrm{Ord}(\alpha)$ holds by assumption. Does $\mathrm{Ord}^M(\alpha)$ hold? Well, if $\forall y \: \forall z \: y\in \alpha \rightarrow \cdots$ holds, then certainly $\forall y \: \forall z\in M \: y\in \alpha \rightarrow \cdots$ holds as well. Therefore, $\mathrm{Ord}(x)$ is reflected in any transitive set.

The property of $\mathrm{Empty}(x)$ and $\mathrm{Ord}(x)$ of being reflected in every transitive set, is a general property of a formula called absoluteness. Informally, $\mathrm{Empty}(x)$ and $\mathrm{Ord}(x)$ are formulas that only talk about sets within $x$, and therefore if the model is transitive and has $x$, it has all the relevant elements to decide the truth of the formula:

\textbf{Absoluteness.} A formula $\phi$ is called \textbf{absolute} across transitive models if $\phi \leftrightarrow \phi^M$ for any transitive set $M$. Absoluteness is a useful concept in reasoning about formulas. There is an important class of formulas of ``low complexity'' that are absolute in every transitive set. 

\textbf{Logical complexity.} The quantifiers of a formula, $\forall x$ and $\exists x$, determine its logical complexity. Quantifiers are \textbf{bounded} when they are of the form $\forall x \in y\: \psi$ and $\exists x \in y \: \psi$ (formally, $\forall x\: x\in y \rightarrow \psi$ and $\exists x \: x\in y \land \psi$), and otherwise they are \textbf{free}. Bounded quantifiers let the variable range over a given set, while free quantifiers let the variable range over the universe. Bounded quantifiers do not make the formula any more complex. Free quantifiers increase the logical complexity of the formula. A formula is of the lowest complexity $\Delta_0$ if all its quantifiers are bounded. A formula is $\Sigma_1$ if it is $\exists x \: \psi$ where $\psi$ is $\Delta_0$. A formula is $\Pi_1$ if it is $\forall x \: \psi$ where $\phi$ is $\Delta_0$. Here, $x$ could represent a tuple so these are equivalent to allowing the formula to be $\exists x_1 \: ...\exists x_k \: \phi$ or $\forall x_1 \: ...\forall x_k \: \phi$, respectively. Higher complexity formulas involve alternating quantifiers, but we will not need to discuss this topic. The following important properties are easy to prove:

\begin{quote}
Whenever $M$ is a transitive set, $\phi$ is $\Delta_0$ and free variables range over $M$, 
\begin{itemize}
\item \textbf{Absoluteness of bounded formulas.} $\phi$ is reflected in $M$:  $\phi \leftrightarrow \phi^M$. \textbf{Proof sketch.} The idea is simple: $\phi^M$ is simply $\phi$ with the syntactic change of converting every quantifier $\exists x \: \psi(x)$ to $\exists x \: x \in M \land \psi^M(x)$. Because all quantifiers in $\phi$ are bounded, they are of the form $\exists x \: x\in y \land \psi(x)$, for some $y$. When interpreting $\phi^M$ by definition all variables range over $M$ and so does $y$. Therefore, $x\in y$ implies $x \in M$. Therefore, $\exists x \: x\in y \land x\in M \land \psi^M(x)$ is equivalent to $\exists x \: x\in y \land \psi^M(x)$, which in turn is equivalent to $\exists x \: x\in y \land \psi(x)$  by inductive assumption on the structure of $\phi$. $\square$
\item Moreover, for any $\Sigma_1$ formula $\psi := \exists x \:\phi$, it is true that $\exists x \in M \: \phi \rightarrow \exists x \: \phi$, therefore $\psi^M \rightarrow \psi$. This is obvious: if there is $x$ within $M$ with a given property, then there is $x$ with the property.
\item Conversely, for any $\Pi_1$ formula $\psi := \forall x \:\phi$, it is true that $\psi \rightarrow \psi^M$. This is also obvious: if for all $x$ a property holds, then the property holds for all $x$ in $M$.
\end{itemize}
\end{quote}

Many useful formulas are \textbf{$\Delta_0$ and therefore absolute across transitive sets}, and can be transferred freely across transitive models. Here is a useful list:
\begin{itemize}

\item ``$x$ is empty''; ``$x$ is the union of $y$ and $z$''.\footnote{The formula for ``$x$ is empty'' is $\lnot \exists y \in x \: y = y$ (for instance) and the formula for ``$x$ is the union of $y$ and $z$'' is $(\forall w \in x \: w \in y \lor w \in z) \land \forall w \in y \: w \in x  \land \forall w \in z \: w \in x$. Both formulas are $\Delta_0$.}
\item ``$x$ is transitive''.\footnote{$\forall y \in x \: \forall z \in y \: z \in x$ is $\Delta_0$.}
\item ``$x$ is an ordinal": $\mathrm{Ord}(x)$; ``$x$ is a successor ordinal''; ``$x$ is a limit ordinal''.\footnote{For $x$ being a limit ordinal, $\mathrm{Ord}(x) \land \forall y \in x \: \exists z \in x \: y \in z$ is $\Delta_0$. Successor ordinal is left as an exercise.}
\item ``$ x$ is a natural number''.\footnote{The following formula is $\Delta_0$: $\mathrm{Ord}(x) \land$ ($x$ is not a limit ordinal $\lor x = \emptyset)$ $\land \forall y \in x \: (y = \emptyset \lor y$ is not a limit ordinal$)$.}
\item $\omega$: we only need to look at elements within a set $x$ to tell whether $x$ is indeed equal to $\omega$.\footnote{The formula ``$x = \omega \leftrightarrow (\mathrm{Ord}(x) \land x \not = \emptyset \land x$ is a limit ordinal $\land \forall y \in x \: y$ is a natural number'' is $\Delta_0$.} Therefore, for any transitive model $M$, $\omega^M = \omega$, or in the notation that we will use most frequently, $\aleph_0^M = \aleph_0$. 
\item ``$f$ is a function from $x$ to $y$''; similarly for $f$ is a bijection, surjection, injection etc. Also $\mathrm{Dom}(f)$, $\mathrm{Ran}(f)$ are absolute: we only need to look inside $f$ to see which elements of the pairs $(x,y)\in f$ form the domain and range of $f$. Showing that these formulas are $\Delta_0$ is left as an easy exercise.
\item ``$x$ is finite''.\footnote{"$\exists n \in \omega \: \exists f \in \omega \times \omega \: f$ is a function from $n$ to $x$''. This definition uses the fact that $\omega \times \omega$ is absolute, which is left to the reader.}
\end{itemize}

\textbf{Cardinals are not absolute.} $\mathrm{Card}(\kappa) :=$ ``$\kappa$ is a cardinal'' is not absolute. Unrolling this formula, it says that ``$\mathrm{Ord}(\kappa)\land$ there is no surjection $f$ between $\alpha \in \kappa$ and $\kappa$''. Now let $M$ be a transitive set and let $\kappa \in M$. The part $\mathrm{Ord}(\kappa)$ is absolute, so $\kappa$ is an ordinal with respect to $M$ iff $\kappa$ is an ordinal. However, there could be a surjection $f : \alpha \to \kappa$ outside $M$ but none in $M$. Therefore, it is possible that $\mathrm{Card}^M(\kappa)$ is true but $\mathrm{Card}(\kappa)$ is false. In general, if we have models $M \subset N$, cardinals of $M$ can stop being cardinals in $N$, but any cardinals of $N$ are guaranteed to be cardinals of $M$.

Finally we are ready to express the Reflection Principle:

\begin{quote}
\textbf{Reflection Principle.} Given any formula $\phi(x_1, \ldots, x_k)$, it is a ZFC theorem that there exists a countable transitive set $M$ such that whenever the variables of $\phi$ take values $a_1, \ldots, a_k \in M$, 
$$\phi(a_1, \ldots, a_k) \leftrightarrow \phi^M(a_1, \ldots, a_k)$$We say that $\phi$ is \textbf{reflected} in $M$.
\end{quote}

\begin{quote}
\textit{The Reflection Principle is an algorithm that takes any $\phi$ as input and returns a ZFC proof of $\exists M \: \mathrm{CntTrans}(M) \land \phi \leftrightarrow \phi^M$.} 
\end{quote}

Importantly, $\phi$ can be the conjunction of any finite collection of theorems of ZFC:

\begin{quote}
\textbf{Reflection Principle and Finite Fragments of ZFC.} Given a finite collection of ZFC axioms $\phi_1, ..., \phi_k$, it is a ZFC theorem that there exists a countable transitive set $M$ where the axioms are reflected: $\exists M \: \mathrm{CntTrans}(M) \land \phi_1^M \land \cdots \land \phi_k^M$.
\end{quote}

This is a stunningly powerful result: there exist countable transitive set models for arbitrarily large finite fragments of ZFC. Consequently, any theorem of ZFC---basically, any mathematical theorem---can be converted into a theorem holding within a countable transitive set model.

This also implies that there cannot be a finite first-order logical axiomatization of ZFC. If there was, ZFC could prove its own consistency through proving the existence of a model of all its axioms, which is impossible by G\"{o}del's second incompleteness theorem.

The Reflection Principle is at the heart of the proof of independence of CH. The proof is not hard; we will skip it and come back to it in Section 6. For now, it is important to keep in mind the following point: ZFC does not have any way to ``talk'' about formulas as objects. Only after fixing formula $\phi$ can the Reflection Principle algorithm be run to output a ZFC proof that there is a countable transitive set $M$ such that $\phi \leftrightarrow \phi^M$. However, we can arithmetize  ZFC because the collection of axioms is recursive, and prove the Reflection Principle in Peano Arithmetic.

\section{Outline of the proof of independence of $\mathrm{CH}$}

\subsection{Statement}

To recap, the natural question we will be addressing is whether there is an infinity between $|\mathbb{N}| = |\omega| = \aleph_0$ and $|\mathbb{R}| = c = 2^{\aleph_0}$. The formula $\mathrm{CH}$ is the statement $2^{\aleph_0} = \aleph_1$. We will show that a proof of either $\mathrm{CH}$ or $\lnot \mathrm{CH}$ leads to contradiction. 

It is clear that a set $S$ such that $|\mathbb{N}| < |S| < |\mathbb{R}|$ can never be demonstrated the way $\mathbb{N}$ and $\mathbb{R}$ are: if it ever did, a contradiction would ensue. But neither can inexistence be proven. In fact, the inability to properly size $c$ in the pecking order of infinities is extreme: there can be 2, 3, 2022, countably many, or uncountably many infinities in-between $\aleph_0$ and $c$.\footnote{However, Konig's theorem restricts the precise position of $c$: it cannot be the union of a collection of fewer than $c$-many sets each of which is smaller than $c$. This topic is addressed in the section ``What values can $2^{\aleph_0}$ take?'' below.}

G\"{o}del showed that $\mathrm{CH}$ is consistent by the \textit{inner model} method (1939), a topic beyond the scope here; a simple argument shows that this method cannot show that $\lnot \mathrm{CH}$ is consistent.\footnote{The topic is beyond the scope of this manuscript and for excellent further reading material I refer the reader to Smullyan and Fitting (2010), but a brief summary is in order. In the inner model method, a formula $L(x)$ is devised for the universe of constructible sets, analogous to $V(x)$ for the universe of all sets. The difference is that at every stage $\alpha$, $L_{\alpha +1}$ is defined to contain the sets that are definable by formulas with parameters from $L_{\alpha}$ instead of $V_{\alpha+1}$ being the power set of $V_{\alpha}$. The collection $x : L(x)$ is a proper class containing all the ordinals and such that every axiom of ZFC $\phi$ is proven in its relativized version $\phi^L$, i.e., replacing all $\forall x \ldots$ with $\forall x \: L(x) \rightarrow \ldots$ and all $\exists x \ldots$ with $\exists x \: L(x) \land \ldots$. Importantly, a formula is defined $\mathrm{Constr}(x)$ that states that $x$ is constructible, which is true iff $x$ is in $L$, i.e., if $L(x)$. And not only that, $\mathrm{Constr}(x)$ is true over $L$, meaning that the relativized version $\mathrm{Constr}(x)^L$ is provable. This implies that $\forall x \: L(x) \rightarrow \mathrm{Constr}(x)^L$ is proven. Therefore the class of constructible sets $x : L(x)$ obeys all ZFC and is a model of the \textbf{axiom of constructibility}, $\forall x \: \mathrm{Constr}(x)$ or as it is better known, $V = L$. Therefore $V = L$ is consistent with ZFC. Moreover, it is shown that within $L$ the continuum hypothesis holds, i.e., $CH^L$ is a theorem of ZFC. Finally, it is shown that $L$ is the smallest model of ZFC containing all the ordinals, and any model of ZFC containing all the ordinals has to include all of $L$.

Can the same technique be applied in the other direction to show that $\lnot \mathrm{CH}$ is consistent? Let's assume that we could devise a formula $M(x)$ such that for every ZFC axiom $\phi$ the version $\phi^{M}$ holds and such that $(\lnot \mathrm{CH})^M$ holds. The class construction $L$ can be relativized to $M$, as $L^M$. This is just a formula $L^M(x)$ that defines the class of constructible $x$ within the class $M(x)$, therefore $L^M \subseteq L$. However, $L$ is the smallest proper class model of ZFC, and therefore $L^M = L$, which implies that $L \subset M$. Simply, for any $x$, $L(x) \rightarrow M(x)$. Now because $\mathrm{CH}$ holds in $L$ and $\lnot \mathrm{CH}$ holds in $M$, it follows that $M \neq L$---$M$ cannot be identical to $L$ so it has to be properly larger. But then it follows that $\lnot \forall x \: \mathrm{Constr}(x)$ or $V \neq L$ because $L$ is precisely the class of sets $x$ such that $\mathrm{Constr}(x)$. This is a contradiction with consistency of $V=L$. In conclusion, there cannot be such a formula $M(x)$ where $\lnot \mathrm{CH}^M$ holds.} Instead, Cohen's method of forcing starts with a ground model $M$  and extends to a larger, \textit{outer model} $N$ in which $\mathrm{CH}$ fails.

\subsection{The Ground Model $\mathrm{M}$}

\textbf{Assume $\mathrm{CH}$ is provable for contradiction.} The proof is an algorithm that transforms a hypothetical proof of $\mathrm{CH}$ to an inconsistency proof of ZFC. At a high level, a model $M$ is established for the finitely many axioms of ZFC that are used in the hypothetical proof of $\mathrm{CH}$. The proof of $\mathrm{CH}$ is copied within $M$, and then $M$ is expanded to a model $N\supset M$ in which there exist more than $\aleph_2$ subsets of $\mathbb{N}$ and therefore $\lnot \mathrm{CH}^N$ is proven. At the same time the proof of $\mathrm{CH}$ is shown to still go through in $N$ to derive  $CH^N$, leading to a contradiction. The same technique, building an \textit{outer model} of a ground model $M$ through the forcing method, is general and can be used to also prove consistency of $\mathrm{CH}$. This exposition will focus on the proof of consistency of $\lnot \mathrm{CH}$, which is the harder part. Consistency of $\mathrm{CH}$ is covered in Section 6.2, after the main exposition.

\textbf{Relativize the proof to a set $M$.} Any proof of $\mathrm{CH}$ would be  an ordered list of formulas, $\phi_1, ..., \phi_n$ such that each $\phi_i$ is either an axiom or a logical consequence of preceding formulas in the list and $\phi_n = \mathrm{CH}$. We might as well assume that $\phi_1, ..., \phi_k$ ($k<n$) are the axioms used in the proof. By the Reflection Principle, a ZFC proof can be produced that there is a countable transitive set $M$ where all $\phi_i$s hold, and therefore the proof goes through relativized to $M$: first prove within ZFC that $\exists M \: \mathrm{CntTrans}(M) \land \phi_1^M \land \cdots \land \phi_k^M$ and then continue on with the rest of the proof: $\phi_{k+1}^M, \ldots, \phi_n^{M}$. An additional subtlety is that the list $\phi_1, \ldots, \phi_k$ needs to contain (1) all axioms needed for the proof of $\mathrm{CH}$; (2) all axioms needed for the Reflection Principle to be applied to $\phi_1 \land \cdots \land \phi_k$; (3) all axioms needed for subsequent constructions below. Still this is a finite list of axioms, so without loss of generality assume that $\phi_1, \ldots, \phi_k$ includes all of them.

So far, the algorithm has produced a ZFC proof of $CH^{M}$.

\textbf{The Skolem paradox.} This is a moment to pause and reflect. How is it possible that ${M}$ is a \textit{countable} set reflecting $\mathrm{CH}$? The whole question of $\mathrm{CH}$ is whether an \textit{uncountable} infinity exists between $\aleph_0$ and continuum. How can this question be relativized to a countable model ${M}$? The seeming conflict between uncountability and the existence of countable models is known as the Skolem paradox.

This is a part of the proof that initially looks like cheating. First, note that ${M}$ contains an intact copy of $\aleph_0$. Then, $M$ does \textit{not} contain all of $2^{\aleph_0}$ but instead a countable subset of $2^{\aleph_0}$, call it $(2^{\aleph_0})^{M}$, which acts as the continuum within ${M}$. $(2^{\aleph_0})^{M}$ is a countable ordinal that is uncountable with respect to $M$---call it uncountable$^{M}$---because within $M$ there is no bijection between $\aleph_0$ and $(2^{\aleph_0})^{M}$. Moreover, any property of $2^{\aleph_0}$ expressed by a formula, such as containing all limits of Cauchy sequences (if viewed as $\mathbb{R}$), is obeyed: $(2^{\aleph_0})^{M}$ contains all limits of \textit{those Cauchy sequences }\textit{that exist in ${M}$}. This may all sound bizarre, but the Reflection Principle and the more general Lowenheim-Skolem reflection theorems were known well before 1963 when Cohen proved the consistency  of $\lnot \mathrm{CH}$, and are way simpler to prove.

\subsubsection{Extend ZFC to postulate a countable transitive model $\mathrm{M}$}

\begin{quote}
\textit{Since every proof of a formula $\phi$ can be transformed to $\exists M \: \mathrm{CntTrans}(M) \land \phi^M$, we may as well assume that there exists a countable transitive model $\mathrm{M}$.}
\end{quote}

\textbf{It is safe to extend ZFC to ZFC+.} Following Weaver (2014), this observation suggests an extension of ZFC to ZFC+ that simplifies the proof: introduce a new symbol '$\mathrm{M}$', which acts as a term, together with the axiom $CntTrans(\mathrm{M})$ and all axiom schemas of ZFC relativized to $\mathrm{M}$. These additions cannot introduce a contradiction where none existed: suppose that a proof of contradiction follows from new axioms $\psi_1^{\mathrm{M}}, ..., \psi_l^{\mathrm{M}}, \mathrm{CntTrans}(\mathrm{M})$. Now $\psi_1 \land \cdots \land \psi_l$ is trivially a theorem, therefore $\exists M \: \mathrm{CntTrans}(M) \land \psi_1^{M} \land \cdots \land \psi_l^{M}$ is a theorem by the Reflection Principle as explained above. A perfect syntactic copy of the proof of contradiction, with $\mathrm{M}$ replaced with $M$, brings a contradiction within ZFC. Therefore, if ZFC is consistent, so is ZFC+.

In ZFC+, $\mathrm{M}$ is a symbol that represents \textit{by fiat} a countable transitive model of ZFC. All ZFC theorems are theorems of ZFC+ relativized to $\mathrm{M}$. The formula $\mathrm{Ord}(x)$ is $\Delta_0$, therefore it is derivable in ZFC+ that being an ordinal in $\mathrm{M}$ means being an ordinal: $\forall x \in \mathrm{M} \: \mathrm{Ord}(x)^{\mathrm{M}} \leftrightarrow \mathrm{Ord}(x)$. Also, $\aleph_0^{\mathrm{M}} = \aleph_0$, i.e., ``($x$ is the smallest infinite cardinal)$^{\mathrm{M}}$'' $\leftrightarrow$ ``$x$ is the smallest infinite cardinal''.

It is derivable in ZFC+ that there exists an ordinal $\alpha_{\mathrm{M}} = \sup\{\alpha \: | \: \mathrm{Ord}(\alpha) \land \alpha \in \mathrm{M}\}$ that is countable,  $\not \in \mathrm{M}$, and is the supremum of ordinals in $\mathrm{M}$. Also that the sets $\aleph_0^{\mathrm{M}} < \aleph_1^{\mathrm{M}} < \aleph_2^{\mathrm{M}} < \ldots < \aleph_{\alpha}^{\mathrm{M}} < \ldots$  are ordinals $\in \mathrm{M}$ as long as $\alpha \leq \alpha_{\mathrm{M}}$, are countable, are uncountable$^\mathrm{M}$ and cardinals$^\mathrm{M}$. This means that it is a theorem of ZFC+ that there is no surjection $f\in \mathrm{M}$ between any of these $\mathrm{M}$-cardinals: the formula $\lnot \exists f \in  \mathrm{M} \: \mathrm{Surj}(f) \land \mathrm{Dom}(f, \aleph_{\alpha}^{\mathrm{M}}) \land \mathrm{Ran}(f, \aleph_{\alpha+1}^{\mathrm{M}})$ is a theorem of ZFC+ for any ordinal $\alpha \in \mathrm{M}$. 

\textbf{The uncountable cardinals in $\mathrm{M}$.} All infinite sets $\aleph_1^{\mathrm{M}} < \aleph_2^{\mathrm{M}} <\ldots$ are countable and uncountable$^\mathrm{M}$. This is a simple point that needs to be made clear, otherwise it can be enormously confusing. $\aleph_1^{\mathrm{M}}$ is nothing other than a formula of ZFC+. It says ``there exists $x\in \mathrm{M}\land x$ is a countable ordinal $\land$ there is no surjection $f: \omega \to x$ such that $f\in \mathrm{M}$ $\land$ for any ordinal $\alpha \in \mathrm{M}$ with $\alpha < x$, there is a surjection $f\in \mathrm{M}$, $f: \omega \to \alpha$''. And $\mathrm{M}$ in turn is just a symbol of ZFC+. Beyond that, $\mathrm{M}$ does not ``exist'' other than syntactically---\textit{no proof that $\mathrm{M}$ exists is  possible in ZFC by G\"{o}del's second incompleteness theorem}---and neither does $\aleph_1^{\mathrm{M}}$ exist other than as a shorthand for the above formula that is provable in ZFC+.

\textbf{Relativize the assumed proof of $\mathrm{CH}$ to $\mathrm{M}$.} Let $\phi_1, ..., \phi_k$ be the axioms appearing in the proof of $\mathrm{CH}$. The same proof, using $\phi_1^{\mathrm{M}}, ..., \phi_k^{\mathrm{M}}$ instead yields $CH^{\mathrm{M}}$, or equivalently $(2^{\aleph_0})^{\mathrm{M}} = \aleph_1^{\mathrm{M}}$. 

\subsection{Breaking $\mathrm{CH}$ by extending $\mathrm{M}$ to $\mathrm{M}[G]$}

\textbf{A countable collection of $\aleph_2^\mathrm{M}$-many subsets of $\aleph_0$.} Now comes a devious trick of leveraging the countability of $\mathrm{M}$ to construct an extension of $\mathrm{M}$ that breaks $\mathrm{CH}$. In particular, precisely because $\aleph_1^{\mathrm{M}}, \aleph_2^{\mathrm{M}}, \ldots$ as well as $(2^{\aleph_0})^{\mathrm{M}}$ are countable, it is easy to prove in ZFC+ the existence of a countable collection of $\aleph_2^{\mathrm{M}}$-many distinct subsets of $\aleph_0$ (remember, $\aleph_0^{\mathrm{M}} = \aleph_0$), none of which are in $\mathrm{M}$. If these subsets can be added to $\mathrm{M}$ in a way so as to produce a new model, $\mathrm{CH}$ will fail in this new model. And this is precisely what the forcing construction does.

\begin{quote}
\textbf{Representing $\aleph_2^{\mathrm{M}}$-many subsets of $\aleph_0$.} These subsets can be represented as a function $g: \aleph_2^{\mathrm{M}} \times \aleph_0 \to \{0,1\}$:  every ``row'' $\alpha \in \aleph_2^{\mathrm{M}}$ of the ``table'' $g$ is a function $\aleph_0 \to \{0,1\}$, or equivalently a subset of $\omega$. As will become clear, it is convenient to represent $g$ with the collection of \textit{finite} partial functions $p: \aleph_2^{\mathrm{M}} \times \aleph_0 \to \{0,1\}$ such that $p \subset g$, i.e., all finite partial functions that agree with $g$ on all elements of their domain. All such $p$ are in $\mathrm{M}$. (All hereditarily finite sets---sets with a finite number of members that are also hereditarily finite, recursively---are in $\mathrm{M}$. Given a hereditarily finite set $f$, a finite number of Pairing applications establishes that $f\in \mathrm{M}$.)  We will denote this collection with $G := \{p \:|\: p$ is a finite partial function  $\aleph_2^{\mathrm{M}}\times \aleph\to \{0,1\}$ and $\forall (\alpha, n) \in \mathrm{Dom}(p) \: p(\alpha, n) = g(\alpha, n)\}$. The ``rows'' of $g$, which by abuse of terminology will also be referred to as rows of $G$, will appear below in several examples and will be denoted by $s_{\alpha}$ for the $\alpha$th row.
\end{quote}

The goal is now to add $G$ to $\mathrm{M}$ in a manner that results in a larger model: extend $\mathrm{M}$ to a set $\mathrm{M}[G]$ that contains every set in $\mathrm{M}$ as well as $G$ and all sets that are implied by the presence of $G$ through ZFC axioms. $\mathrm{M}[G]$ should be a model: for every axiom $\phi^{\mathrm{M}}$ there should be a ZFC+ proof of $\phi^{\mathrm{M}[G]}$. Then, the proof of $CH^{\mathrm{M}}$ can be copied to yield a ZFC+ proof of $CH^{\mathrm{M}[G]}$. If in addition $\aleph_2^{\mathrm{M}} = \aleph_2^{\mathrm{M}[G]}$, i.e., the precise same (countable) ordinal that is the third smallest infinite cardinal in $\mathrm{M}$ is the third smallest infinite cardinal in $\mathrm{M}[G]$, then a contradiction ensues.

\subsubsection{New sets that need to be in $\mathrm{M}[G]$}

What kind of new sets are required to be in $\mathrm{M}[G]$ for all ZFC axioms $\phi$ to hold as $\phi^{\mathrm{M}[G]}$? All kinds of sets will be required by Pairing, Power, Union, Comprehension, Replacement and AC.  The Comprehension axiom presents the clearest challenge: for any formula $\phi$ of $k+1$ free variables, for any $y\in \mathrm{M}[G]$ and any sets $w_1, \ldots, w_k \in \mathrm{M}[G]$, the set $z := \{x \in y \:|\: \phi^{\mathrm{M}[G]}(x, w_1, \ldots, w_k)\}$ should be in $\mathrm{M}[G]$. 

Consider the set $\mathbb{P}$ of all finite partial functions $\aleph_2^{\mathrm{M}}\times \aleph\to \{0,1\}$. Each $p\in \mathbb{P}$ is in $\mathrm{M}$ as discussed and therefore $\mathbb{P}\in \mathrm{M}$. $\mathbb{P}$ will be an important structure: it is a partial order with sets $p\in \mathbb{P}$ ordered by $\subseteq$, and $G$ is a subset of $\mathbb{P}$ that is not in $\mathrm{M}$. In the next section, we will define $\mathbb{P}$ as a \textit{forcing notion} and $G$ as a \textit{maximal ideal} over $\mathbb{P}$. For now, let's just focus on the question of what new sets in $\mathrm{M}[G]$ are implied by $G$.

Take any set $x \in \mathrm{M}$. Having access to the ``key'' $G\subsetneq \mathbb{P}$, can we generate any new subsets $y \subset x$ that may not already be in $\mathrm{M}$?

One obvious way is to find any relation $r\subset x \times \mathbb{P}$ such that $r\in \mathrm{M}$ and take the preimage of $x$ under $G$: $\{z \in x \:|\: (z, p) \in r$ and $p \in G\}$ is in general a set not in $\mathrm{M}$, and if $\mathrm{M}[G]$ is a model then it contains this set by Comprehension$^{\mathrm{M}[G]}$. As an aside, every $x \in \mathrm{M}$ is the preimage of $r := x \times \mathbb{P}$ under $G$. Also, $G$ is the preimage of $r := \{(p,p) \:|\: p\in \mathbb{P}\}$ under $G$. Therefore, the collection of preimages under $G$ of $x\times r$ for $x, r \in \mathrm{M}$ and $r \subseteq x\times \mathbb{P}$ includes all of $\mathrm{M}$ as well as $G$. 

\textbf{Taking preimages transfinitely.} Is this a general form of all sets that are required to be in $\mathrm{M}[G]$? Well, what about preimages under $G$ of $y \times r'$ where $y$ is the preimage of $x$ under $r$? These sets are also required to be in $\mathrm{M}[G]$ if $\mathrm{M}[G]$ is to be a model, and specifically if Comprehension$^{\mathrm{M}[G]}$ holds. One stage of taking preimages under $G$ does not suffice; the construction needs to be applied transfinitely. The details of how to do this will be given in the next section but the idea is simple. Consider a relation $r \subset r' \times \mathbb{P}$ and $r'$ in turn being a relation $\subset r'' \times \mathbb{P}$, with $r'' \in \mathrm{M}$ and $r''$ a relation of the same form, and so on. The preimage under $G$ of $r$, call it $r^G$, can be taken recursively as $r^G := \{r'^G \: |\: (r',p) \in r$ and $p\in G\}$. From now on, ``preimage under $G$'' will mean taking preimages recursively.

\begin{quote}
\textbf{Sets that could lead to contradiction.} Certain sets that contain information about $\mathrm{M}$ lead to contradiction if included in $\mathrm{M}[G]$. For example, consider the countable ordinal $\alpha_{\mathrm{M}}$ that is the supremum of all ordinals in $\mathrm{M}$. As it turns out, $\alpha_{\mathrm{M}} \in \mathrm{M}[G]$ would lead to contradiction:\footnote{\textbf{Proof sketch.} In ZFC+ the existence of $\mathrm{CntTrans}(\mathrm{M})$ implies the existence of a minimal model $\mathrm{M_{MIN}}$ that is a subset of any transitive model. This result is covered by Cohen (1966) and is beyond the scope here. It is worth mentioning that every element $a$ of $\mathrm{M_{MIN}}$ is associated with a formula $\phi_a(x)$ that is true iff $x = a$---every element is named. Of course, $\mathrm{M_{MIN}}$ is countable. Let $\alpha_0$ be the supremum of all ordinals in $\mathrm{M_{MIN}}$. It is consistent with ZFC+ that there is no model $\mathrm{N}$ with supremum of ordinals $\alpha_1 > \alpha_0$. To see why, let $\mathrm{N}$ be a model with supremum $\alpha_1$ such that $\alpha_1 > \alpha_0$ is minimal. If no such $\mathrm{N}$ exists, we are done. Otherwise, $\mathrm{M_{MIN}} \in \mathrm{N}$ and ``there is no model with supremum of ordinals $> \alpha_0$"$^{\mathrm{N}}$ is a ZFC+ theorem. Therefore, it is consistent that there is no model with supremum $\alpha_1 > \alpha_0$, because the statement holds relativized to $\mathrm{N}$ and $\mathrm{N}$ is a model. Now consider $\alpha_{\mathrm{M}} \in \mathrm{M}[G]$. Because the construction of $\mathrm{M}[G]$ can be carried starting with any $\mathrm{M}$, even with $\mathrm{M} = \mathrm{M_{MIN}}$, it is inconsistent for $\mathrm{M}[G]$ to have a supremum of ordinals greater than that of $\mathrm{M}$. Therefore, $\alpha_{\mathrm{M}} \not \in \mathrm{M}[G]$. $\square$} Because $\alpha_{\mathrm{M}}$ is countable, it is just a shuffling of $\omega$ that can be expressed as a function $f: \omega \to \omega$, which can be expressed as a subset $s_f$ of $\omega \times \omega$, which can be expressed as a subset of $\omega$. Call that subset $s_{\alpha_{\mathrm{M}}}$. The maps $a_{\mathrm{M}} \to f \to s_f \to  s_{\alpha_{\mathrm{M}}}$ are definable by rather simple formulas of ZFC---they are absolute---and therefore the maps and their inverses are present in any model. So if $s_{\alpha_{\mathrm{M}}} \in \mathrm{M}[G]$ then $\alpha_{\mathrm{M}} \in \mathrm{M}[G]$ and we get a contradiction. 
\end{quote}

Therefore, not any $G$ works. If $G$ contains $s_{\alpha_{\mathrm{M}}}$ for instance, it won't produce a model. To show how $\mathrm{M}[G]$ can be made into a ZFC model, we will: (1) describe the new sets implied by $G$, and define $\mathrm{M}[G]$ as a collection of these sets; (2) devise a strategy of how to prove $\phi^{\mathrm{M}[G]}$ for every ZFC axiom $\phi$, which will illuminate the conditions needed so that the proofs go through, and (3) based on that reasoning, define the property that $G$ needs to have---\textit{genericity}---so that we can produce a ZFC proof of $\phi^{\mathrm{M}[G]}$ for every axiom $\phi$  and especially for the Comprehension axiom schema. Once these tasks are accomplished, we will finally prove that $\mathrm{M}[G]$ is indeed a model as long as $G$ is \textit{generic}.

Let's organize the above discussion with some definitions. First, let's explore the structure of $\mathbb{P}$, which is in $\mathrm{M}$ and will be leveraged to convert axioms $\phi^{\mathrm{M}}$ to their counterparts $\phi^{\mathrm{M}[G]}$. $\mathbb{P}$ will be called a \textit{forcing notion} for reasons that will become clear. 

\subsubsection{A forcing notion $\mathbb{P}$ and ideals over $\mathbb{P}$}

\textbf{Definition.} A \textbf{forcing notion} is a partial order set (poset) $\mathbb{P} \in \mathrm{M}$ consisting of elements $p \in \mathbb{P}$ that in the context of consistency of $\lnot \mathrm{CH}$ are \textit{finite} partial functions $\aleph_2^{\mathrm{M}} \times \aleph_0 \to \{0,1\}$.\footnote{Forcing is a general technique for proving independence results beyond $\mathrm{CH}$, and in some contexts elements $p$ are infinite; however, they have to be $\in \mathrm{M}$.} $\mathbb{P}$ is ordered by inclusion: $p \leq q$ whenever $p \subseteq q$ (notably, $p \leq p$). Two elements $p, q$ are \textbf{compatible} if $\exists r \: p \leq r \land q \leq r$. In other words $p, q$ are compatible if they agree on the values they assign on the intersection of their domains and thus can be patched together to form a bigger partial function $r$. If they ``clash'' at a given entry $(\alpha, n)$, then they are incompatible, denoted by $p \perp q$. 

The elements $p\in \mathbb{P}$ are  called \textbf{forcing conditions} because each of them forces additional sentences to be true or false in $\mathrm{M}[G]$. For example, take a simple $p = \{((\omega+3, 7), 1), ((11, 3), 0)\}$. This $p$ maps the element $(\omega+3, 7) \in \aleph_2^{\mathrm{M}}\times \aleph_0$ to $1$ and the element $(11,3)$ to $0$, and these choices in turn will force certain sentences in $\mathrm{M}[G]$ to be true or false. 

\textbf{Construction of $G$ from partial functions.} With the poset $\mathbb{P}$ thus defined, one way to think of  $G$ is as a \textit{subset} of $\mathbb{P}$ consisting of all $p$ that can be patched together to form some complete function $\bigcup G:  \aleph_2^{\mathrm{M}} \times \aleph_0 \to \{0,1\}$. In order to make this work, we need to mathematically define what it means for $G$ to be a subset of $\mathbb{P}$ that is patchable together to form a complete function. The structure that works is for $G$ to be a \textit{maximal ideal} of the poset $\mathbb{P}$, defined as follows:

\begin{quote}
\textbf{Definition.} $G$ is an \textbf{ideal} over $\mathbb{P}$ if: 
\begin{enumerate}
\item Elements are pairwise compatible: for any $p,q \in G$ there is $r\in G$ such that $p\leq r$, $q\leq r$. 
\item $G$ is downward closed: if $p\leq q$ and $q \in G$ then $p \in G$.
\end{enumerate}
$G$ is a \textbf{maximal ideal} if for any $p\in \mathbb{P}$ either $p\in G$ or $p$ is incompatible with some $q\in G$.
\end{quote}

In other words, $G$ is an ideal over poset $\mathbb{P}$ if (1) every pair of elements of $G$ are compatible---patchable together to form a bigger function, and (2) it is downward closed---all partial functions that are subsets of $\bigcup G$ are included in $G$. $G$ is a maximal ideal if no other element of $\mathbb{P}$ can be added to it. Then, $\bigcup G$ is a complete function because if some element $(\alpha, n)$ was not mapped to $\{0,1\}$, then $\{((\alpha, n), 0)\} \in \mathbb{P}$ would be an element $\not \in G$ compatible with every element of $G$.

The notion of ideal over a poset is related to ideals over Boolean algebras, a topic that is beyond the scope here.\footnote{Ideals of a ring form the kernel of a ring homomorphism, which consists of all the elements that a homomorphism maps to $0$. Each homomorphism of a Boolean algebra to the binary algebra $\{0,1\}$ defines a maximal ideal of all elements that map to $0$. A poset is not a Boolean algebra but it can be embedded in a Boolean algebra. In most literature on $\mathrm{CH}$, the dual concept of a \textit{filter} is used and the poset is reversed for this to work: $p\leq q$ whenever $p \supseteq q$. The complement of a maximal ideal is an ultrafilter that maps elements to 1. Here I follow the approach of Weaver (2014) that I find more natural because $p \subseteq q$ resembles $p \leq q$.}

Figure 2 displays a small portion of the poset $\mathbb{P}$ and complete function $\bigcup G$.

\begin{quote}
\textit{We want pairs of elements. to be compatible so that they are patchable together, we want downward closure so that $G$ is as simple as possible, and we want maximality so that $\bigcup G$ forms a complete function.}
\end{quote}

\begin{figure}
\centering
\includegraphics[width=400pt]{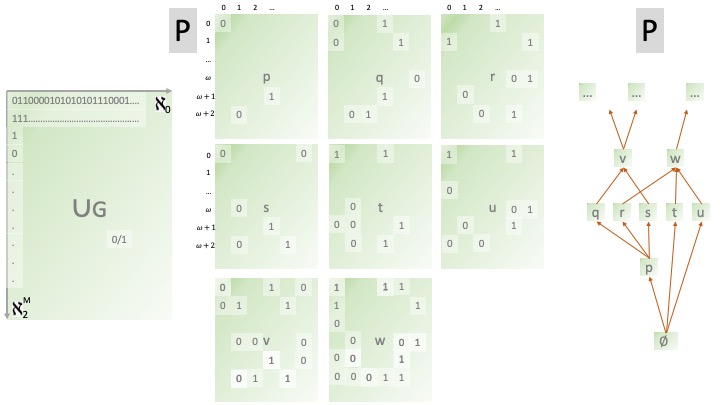}
\caption{\textbf{The forcing notion $\mathbb{P}$ and table $\bigcup G$.} $\mathbb{P}$ consists of finite partial functions $\aleph_2^{\mathrm{M}}\times \aleph_0 \to \{0,1\}$, and all of $\mathbb{P}$ is in $\mathrm{M}$. The members of $\mathbb{P}$ are partially ordered by inclusion, and the ones that are compatible can be patched together to form bigger functions. For example, $p, q$, and $s$ are patchable to form $v$. $p$ is compatible with $r$ as well, but clashes on the top left entry with $t$ and $u$. Also $r, t$ and $u$ are patchable together to form $w$. $\mathbb{P}$ continues upward to contain larger and larger functions, but always finite, whose limits are infinite (but not necessarily complete) functions $\aleph_2^{\mathrm{M}} \times \aleph_0 \to \{0,1\}$. Every infinite partial function $\aleph_2^{\mathrm{M}} \times \aleph_0 \to \{0,1\}$ is a limit $\not \in \mathbb{P}$. The function $\bigcup G$ that will be constructed in this proof is a \textit{complete} function, a very special kind of limit.} 
\end{figure}

\subsubsection{The extended model $\mathrm{M}[G]$}

Now we can go back to defining what new sets need to be in $\mathrm{M}[G]$ due to the presence of $G$. As discussed, the collection of preimages under $G$ of $x\times r$ for any $x, r \in \mathrm{M}$ with $r \subseteq x \times \mathbb{P}$ needs to be in $\mathrm{M}[G]$, and the process of taking preimages under $G$ needs to be extended transfinitely. Any model that is a superset of $\mathrm{M}$ and contains $G$ has to have all these elements.

To conduct the transfinite construction, let's start by considering what elements are in $\mathrm{M}$.

\textbf{What are the elements of $\mathrm{M}$?} $\mathrm{M}$ is just a symbol of ZFC+, and the elements that provably belong to $\mathrm{M}$ are those elements that provably belong to any countable transitive set that reflects the ZFC axioms. Recall the definition of $V$ as the universe of sets. $V$ is just a ZFC formula: $V(x) := \exists \alpha \: V(\alpha, x)$. The formula $V(\alpha, x)$ given explicitly above can be used to make $V_{\alpha}$ into a set through Comprehension as long as $\alpha$ is given. In other words, $\forall \alpha \: \exists v \: \forall x \: V(\alpha, x) \leftrightarrow x \in v$ is a ZFC theorem that can be proven with transfinite induction.

So, how about $\mathrm{M}$'s elements? Well, because $\mathrm{M}$ is a ZFC model by ZFC+ fiat, all provable ZFC sentences are reflected in $\mathrm{M}$. Every ZFC theorem $\phi$ is automatically a ZFC+ theorem $\phi^{\mathrm{M}}$. In other words, $\forall \alpha\in \mathrm{M} \: \exists v\in \mathrm{M} \: \forall x\in \mathrm{M} \: V^{\mathrm{M}}(\alpha, x) \leftrightarrow x \in v$ is a ZFC+ theorem. Therefore, the $V$ hierarchy can be defined within $\mathrm{M}$ except that it will be shorter and thinner.

\textbf{$\mathrm{M}$ is short and thin.} It is provable in ZFC+ that $V_{\alpha}^{\mathrm{M}}$ is a set in $\mathrm{M}$ as long as $\alpha \in \mathrm{M}$. Also that $\mathrm{M}$ is ``short'' because it contains only the ordinals up to and not including a countable ordinal $\alpha_{\mathrm{M}}$, and ``thin'' because it is missing a bunch of subsets at each stage of $V_{\alpha}^{\mathrm{M}}$. It has to be so: $\mathrm{M}$ is countable and $\omega \in \mathrm{M}$---in fact, $\omega \in V_{\omega+1}^{\mathrm{M}}$ (\textit{why?})---and therefore $V_{\omega + 2}^{\mathrm{M}}$ misses all but a countable collection of subsets of $\omega$. To understand what makes $\mathrm{M}$ ``thin'', the reader is encouraged to go back to the definition of formula $V(\alpha, x)$ and relativize it into $V^{\mathrm{M}}(\alpha, x)$. In particular pay attention to $\mathrm{Pow}(x,y): = \forall z \: z \subseteq x \leftrightarrow z \in y$. This is a key subformula that changes meaning when relativized to $\mathrm{M}$: $\forall z \in \mathrm{M} \: z \subseteq x \leftrightarrow z \in y$. With this ``small'' change, at each successor stage $\beta+1$ of the transfinite construction of $V_{\alpha}$, not all subsets of $V_{\beta}$ are added to $V_{\beta+1}$ \textit{but only those subsets that are in $\mathrm{M}$}. So, which subsets are these? At the very minimum, $\mathrm{M}$ includes all of the minimal ZFC model $\mathrm{M_{MIN}}$ that is implied by the existence of a model, as discussed by Cohen (1966).\footnote{However, the existence of a model cannot be proven within ZFC because then ZFC would be proving its own consistency, which contradicts G\"{o}del's second incompleteness theorem.} Because the only axioms pertaining to $\mathrm{M}$ in ZFC+ are the ZFC axioms$^\mathrm{M}$ and $\mathrm{CntTrans}(\mathrm{M})$, it is neither specified, nor excluded that $\mathrm{M}$ is $\mathrm{M_{MIN}}$ or some strictly larger countable model. So, which subsets of $V_{\beta}$ are included in $V_{\beta+1}$ is largely left unspecified, as long as they are denumerable.

\textbf{How to define new subsets in $\mathrm{M}[G]$ using $G$.} In defining $\mathrm{M}[G]$, each $V_{\alpha}^{\mathrm{M}}$ needs to be modified to $V_{\alpha}^{\mathrm{M}[G]}$. $G$ is available to $\mathrm{M}[G]$ and the extra information in $G$ is available when $V_{\beta+1}^{\mathrm{M}[G]}$ is constructed. The Cartesian product trick outlined above needs to be extended transfinitely. Analogously to sets of the form $x \times r$ above, sets $\tau \in \mathrm{M}$ will be defined recursively to consist of ordered pairs $(\sigma, p)$ with $p \in \mathbb{P}$ and $\sigma$ being similarly defined to $\tau$ but of lower rank. The sets thus defined are all in $\mathrm{M}$ and will be called the \textbf{$\mathbb{P}$-names} because each one of them \textit{names} an element in $\mathrm{M}[G]$, by taking its preimage under $G$.

The definition of $\mathbb{P}$-names that follows may seem counterintuitive at first. 

\textbf{Definition of $\mathbb{P}$-name hierarchy.} Given $\mathbb{P}\in \mathrm{M}$, the set $\mathrm{M}^{\mathbb{P}}$ of $\mathbb{P}$-names is defined. A $\mathbb{P}$-name is a set $\tau \in \mathrm{M}$ of the form $\tau \subseteq \{(\sigma, p) \:|\: \sigma$ is a $\mathbb{P}$-name and $p\in \mathbb{P}\}$. The definition is by transfinite recursion, and is turned to a formula $N(\alpha, x)$ in a similar manner to the formulas $V(x)$ and $V(\alpha, x)$, left as an exercise for the reader:
\begin{itemize}
\item $N_0 := \emptyset$.
\item $N_{\alpha+1} := \mathcal{P}(N_{\alpha} \times \mathbb{P}) \cap \mathrm{M}$, for any successor ordinal $\alpha$.
\item $N_{\alpha} = \bigcup_{\beta < \alpha} N_{\alpha}$, for any limit ordinal $\alpha$.
\item $\mathrm{M}^{\mathbb{P}} := \bigcup_{\alpha \in ON} N_{\alpha}$.
\end{itemize}

The collection of $\mathbb{P}$-names is already in $\mathrm{M}$: $\mathrm{M}^{\mathbb{P}} \subset \mathrm{M}$. Note that $\mathrm{M}^{\mathbb{P}}$ is a ZFC+ set but a proper class and not a set in $\mathrm{M}$.\footnote{A proper class is always of the same height---has the same ordinals---as the universe. If not, then let $\alpha_{MAX}$ be the maximal ordinal of the proper class $P(x)$. Then $\forall x \: (P(x) \rightarrow x \in V_{\alpha_{MAX}+1})$ therefore use Comprehension to form the set $\{x \in V_{\alpha_{MAX}+1} \:|\: P(x)\}$. A proper class with respect to the universe of sets $V$---i.e., a proper class---can be put in a bijection with all ordinals in $ON$. A proper class with respect to a set model such as $\mathrm{M}$ is always a set because $\mathrm{M}$ is a set and moreover it is ``short''. Specifically $\mathrm{M}^{\mathbb{P}}$ has height $\alpha_{\mathrm{M}}$, which is a countable ordinal. This makes $\mathrm{M}^{\mathbb{P}}$ a set with respect to $V$, but of the same height as $\mathrm{M}$ and therefore a proper class with respect to $\mathrm{M}$.} Call it an $\mathrm{M}$-class. The \textbf{name rank} of a $\mathbb{P}$-name $\tau$ is defined the same way as the rank of an element of $V$ to be the least ordinal $\alpha$ such that $\tau \in N_{\alpha}$ and is denoted by $\mathrm{nr}(\tau)$. Each $\mathbb{P}$-name $\tau$ can be thought of as a set of ``potential'' elements. The choice of $G$ determines which of these elements materialize and which ones vanish. For example say $\tau = \{(\sigma, p), (\rho, q), (\pi, r), (\pi, q) \}$, where $\sigma = \{(\pi, r), (\rho, q)\}, \rho = \{(o, q)\}, \pi = \{(o, q), (o, p)\}, o = \{(\emptyset, p) \}$. A couple of examples of how things could turn out depending of $G$:
\begin{itemize}
\item Let's say that $G$ contains $p$ and $r$ but not $q$. Then, $\xi = \{\emptyset\}$, $\pi = \{\xi\} = \{\{\emptyset\}\}$, $\rho = \emptyset$, $\sigma = \{\pi\} = \{\{\{\emptyset \}\}\}$, and $\tau = \{\sigma, \pi\} = \{\{\{\emptyset\}\}, \{\{\{\emptyset\}\}\}\}$.
\item Instead, say that $G$ contains $q$ but not $p, r$. Then $\xi = \emptyset$, $\pi = \{\emptyset \}$, $\rho = \{\emptyset\} = \pi$, $\sigma = \{\{\emptyset\}\}$, and $\tau = \{\rho, \pi \} = \{\pi\} = \{\{\emptyset\}\}$. Notice how $\rho, \pi$ collapse here but not in the previous case.
\end{itemize}

Given any $\tau \in \mathrm{M}^{\mathbb{P}}$, the preimage of $\tau$ under $G$ has to be in $\mathrm{M}[G]$. What is not clear at all is whether those are all the elements $\mathrm{M}[G]$ needs to be a model. Indeed we will show that as long as $G$ is \textit{generic} (a property that we will define soon), those are all the elements required.

The preimages of $\mathbb{P}$-names under $G$ taken with transfinite recursion will be called $G$-evaluations:

\textbf{Definition: the evaluation of $\tau$ by $G$.} Given $\mathbb{P}$-name $\tau$, its evaluation $\tau^G \in \mathrm{M}[G]$ is defined recursively to be $\{\sigma^G \:|\: (\sigma,p) \in \tau$ for some $p \in G \}$. 

\begin{quote}
\textbf{Cosmetic change to $\mathbb{P}$-names.} Notice that if $p\leq q$ and $(\sigma, p) \in \tau$, then whenever $q \in G$, by downward closure $p \in G$ and therefore $\sigma^G \in \tau^G$ whether or not $(\sigma, q) \in \tau$; the presence or absence of  $(\sigma, q)\in\tau$ makes no difference. For that reason, if $(\sigma, p) \in \tau$ we may as well assume that $(\sigma, q) \in \tau$ for all $q \geq p$. This is a cosmetic change to the definition of $\mathrm{P}$-names that simplifies their structure, as described by Weaver (2014). From now on we adopt this convention, and we modify $N_{\alpha}$ to be the subset of elements in the above definition that additionally obey this condition. This will simplify some subsequent proofs.
\end{quote}

Now we are ready to define the extended set $\mathrm{M}[G]$:
$$\mathrm{M}[G] := \{\tau^G \: | \: \tau \in \mathrm{M}^{\mathbb{P}}\}$$That's it! Elegantly, $\mathrm{M}[G]$ is just the set of preimages under $G$ of $\mathbb{P}$-names. 

$\mathrm{M}$ is a subset of $\mathrm{M}[G]$: to show that,  define a \textbf{canonical} $\mathbb{P}$-name $\check{x}\in \mathrm{M}$ for each $x\in \mathrm{M}$ recursively on the rank of $x$ by $\check{x} = \{(\check{y}, p) \:|\: y \in x$ and $p \in \mathbb{P}\}$. The idea is that all the elements of $x$ are tagged by every element of $\mathbb{P}$. No matter which elements $G$ contains, $\check{x}^G$ always equals $\{\check{y}^G \: |\: y \in x\}$ which equals $\{y \:| \: y\in x\} = x$ by transfinite induction on $\mathrm{rank}(x)$, assuming $\check{y}^G = y$ holds for all sets of smaller rank than $x$. Therefore, for all $x \in \mathrm{M}$ we have $\check{x}^G = x$ therefore $x \in \mathrm{M}[G]$.

\textbf{Proposition.} $\mathrm{M} \subset \mathrm{M}[G]$. $\square$

Does $G \in \mathrm{M}[G]$? To show that, we need to demonstrate a $\mathbb{P}$-name evaluating to $G$. The following simple one works:
$$\Gamma := \{(\check{p}, p) \:|\: p \in \mathbb{P}\}$$Then verify that $\Gamma^G = \{\check{p}^G \: | \: p \in G\} = \{p \: |\: p\in G\} = G$. Elegantly, $G\in \mathrm{M}[G]$.

Therefore $\mathrm{M}[G]$ is indeed bigger than $\mathrm{M}$: it contains $G$, which is not in $\mathrm{M}$. 

How about any new ordinals---might $\mathrm{M}[G]$ be ``taller'' than $\mathrm{M}$? No, as can be shown by a simple induction on $\mathrm{rank}(\tau^G)$, which is always $\leq \mathrm{nr}(\tau) < \alpha_{\mathrm{M}}$. Therefore $\mathrm{M}[G]$ contains exactly the same ordinals as $\mathrm{M}$. However, $\mathrm{M}[G]$ contains additional sets at infinite stages $V_{\alpha}^{\mathrm{M}[G]}$ so it is a bit ``thicker'' than $\mathrm{M}$.

To recap, the $\mathbb{P}$-names are already in $\mathrm{M}$ and are hierarchical nested relations between elements of $\mathrm{M}$ and $\mathbb{P}$. Given that $G$ is in $\mathrm{M}[G]$, the preimage of each of these relations has to be in $\mathrm{M}[G]$ if $\mathrm{M}[G]$ has any hope of being a model: Comprehension$^{\mathrm{M}[G]}$ requires that. As it turns out, $\mathrm{M}[G]$ is the minimal model that is an extension of $\mathrm{M}$ and contains $G$: if $N$ is such a model, then $\mathrm{M}[G] \subseteq N$ because  Comprehension requires that every $G$-evaluation of a $\mathbb{P}$-name is a member of $N$. Figure 3 illustrates the construction of $\mathrm{M}[G]$ from $\mathrm{M}$. Of course, we are yet to show how to define $G$ so that that $\mathrm{M}[G]$ is a model of ZFC. Indeed $\mathrm{M}[G]$ will not be a model in general if $G$ is not generic.

\begin{figure}
\centering
\includegraphics[width=400pt]{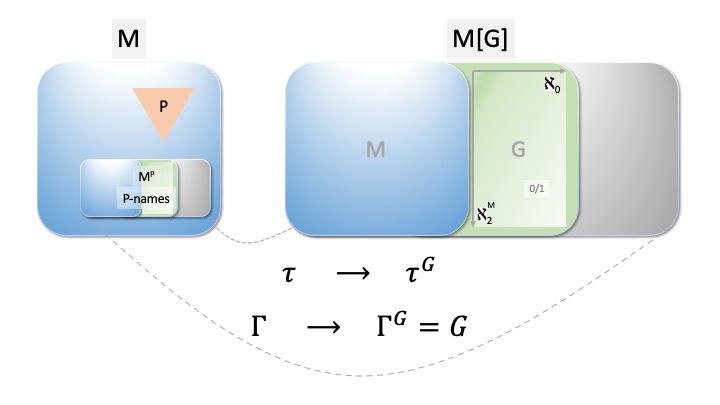}
\caption{\textbf{Definition of extended model $\mathrm{M}[G]$.} The ground model $\mathrm{M}$ contains $\mathrm{M}^{\mathbb{P}}$, the $\mathbb{P}$-names. $\mathrm{M}^{\mathbb{P}}$ is a ZFC+ set and within $\mathrm{M}$ it is a proper class and not a set. $\mathbb{P}$-names are transfinitely defined within $\mathrm{M}$. The preimage of each $\mathbb{P}$-name $\tau$ under $G$ is defined through transfinite recursion: $\tau^G := \{\sigma^G \:|\: (\sigma, p)\in \tau$ and $p \in G\}$. Then, $\mathrm{M}[G] := \{\tau^G \:|\: \tau \in \mathrm{M}^{\mathbb{P}}\}$. Every set  $x\in \mathrm{M}$ is a member of $\mathrm{M}[G]$ by defining the canonical $\check{x} = \{(\check{y}, p) \:|\: y \in x$ and $p \in \mathbb{P}\}$ such that $\check{x}^G = x$. $G$ is the preimage of $\Gamma := \{(p, p) \:|\: p \in \mathbb{P}\} \in \mathrm{M}$.} 
\end{figure}
\textbf{How $\in$ works within $\mathrm{M}[G]$.} In the next sections we will describe the property of genericity of $G$ that ensures $\mathrm{M}[G]$ is a model of ZFC. Before proceeding, some intuition is useful on how $\mathrm{M}[G]$ ``works'' and in particular how ``$\in$'' is affected by the choice of $G$. Think of two members $\sigma^G, \tau^G$ of $\mathrm{M}[G]$. What are the conditions on $\sigma, \tau \in \mathrm{M}^{\mathbb{P}}$ and $G$ that make $\sigma^G \in \tau^G$ true? One may be tempted to say $(\sigma, p) \in \tau$ for some $p \in G$. This seems perfectly reasonable, but unfortunately it has a bug. To detect the bug, let's recall what $\tau^G$ is defined to be: $\tau^G := \{\sigma^G \:|\: (\sigma, p) \in \tau$ for some $p \in G\}$. Before reading further, can you detect the bug? ... If you did, congratulations! If not, it does seem fine to require $(\sigma, p) \in \tau$ for some $p \in G$ verbatim from the definition. Well, what if $(\rho, p) \in \tau$ for some $p \in G$ and moreover $\rho^G = \sigma^G$? This much weaker condition also has the effect of making $\sigma^G \in \tau^G$ true. As we saw in an example above, it is the rule rather than the exception for $\sigma \not = \rho$ while $\sigma^G = \rho^G$. Therefore, to define membership in $\mathrm{M}[G]$ we first have to define equality. Moreover, equality better be extensional: we want $\tau^G = \upsilon^G$ iff $\sigma^G \in \tau^G \leftrightarrow \sigma^G \in \upsilon^G$. Therefore, equality needs to be defined in terms of membership, which needs to be defined in terms of equality, leading to a potential pitfall of circular reasoning. Fortunately, there is a way out: equality and membership are defined simultaneously in a transfinite recursive manner.

Keeping this in mind, let's proceed to ask what conditions $G$ and $\mathrm{M}[G]$ need to have so that $\mathrm{M}[G]$ is a model.

\subsection{Genericity and Forcing}

\subsubsection{How to transfer Comprehension to $\mathrm{M}[G]$}

$\mathrm{M}[G]$ is the set of all preimages under $G$ of $\mathbb{P}$-names in $\mathrm{M}^{\mathbb{P}}$--- the $G$-\textit{interpretations}. For $\mathrm{M}[G]$ to be a model, all these preimages are necessary because of Comprehension, but it is not clear that they are sufficient. To show that $\mathrm{M}[G]$ is a model, we need a metatheorem/algorithm that produces a ZFC+ proof of $\phi^{\mathrm{M}[G]}$ for every input axiom $\phi$.\footnote{As usual, this is a metatheorem. It is impossible to prove within ZFC that a set is a model of ZFC. Instead, we will provide an algorithm that given any axiom $\phi$ as input, produces a ZFC proof of $\phi^{\mathrm{M}[G]}$. The entire algorithm can be shown as a single arithmetic theorem in PA.}

Let's start with the most versatile ZFC axiom for creating new sets: Comprehension. For every set $y$ and formula $\phi(x)$\footnote{Suppressing extra parameters $w_i$.}, all $x\in y$ that satisfy $\phi(x)$ form a new set $z := \{x \in y \:|\: \phi(x)\}$. How do we ensure that there is a ZFC+ proof schema for Comprehension$^{\mathrm{M}[G]}$? For any $\tau^G \in \mathrm{M}[G]$, we need a ZFC+ proof that $z := \{\sigma^G \in \tau^G \:|\: \phi^{\mathrm{M}[G]}(\sigma^G)\}$ is a set in  $\mathrm{M}[G]$. Equivalently, we need to show that $z = \upsilon^G$ for some $\mathbb{P}$-name $\upsilon \in \mathrm{M}^{\mathbb{P}}$.

What axioms can ZFC+ leverage to prove this? An obvious guess is Comprehension$^\mathrm{M}$. We want to use Comprehension$^\mathrm{M}$ to prove Comprehension$^{\mathrm{M}[G]}$. So let's write down how Comprehension$^\mathrm{M}$ could yield $\upsilon$: it would look like $\upsilon := \{x \in y \:|\: \psi_{\phi}^{\mathrm{M}}(x)\}$ for some formula $\psi_{\phi}$ relativized to $\mathrm{M}$ that depends on $\phi$, and for some $y \in \mathrm{M}$ that depends on $\tau^G$.

Every element of $\mathrm{M}[G]$ has a counterpart  in $\mathrm{M}$: its $\mathbb{P}$-name. Every $\sigma^G \in \tau^G$ corresponds to $(\sigma, p) \in \tau$ for some $p\in G$.\footnote{The discussion in the previous section of how ``$\in$'' works in $\mathrm{M}[G]$ emphasized that given $\sigma, \tau \in \mathrm{M}^{\mathbb{P}}$, we can have $\sigma^G \in \tau^G$ without having $(\sigma, p) \in \tau$: we just require some $(\rho, p) \in \tau$ with $\rho^G = \sigma^G$. Are we contradicting ourselves here by saying that $\sigma^G \in \tau^G$ always corresponds to $(\sigma, p) \in \tau$? A moment's reflection will reveal that there is no contradiction. Any $\sigma^G \in \tau^G$ corresponds to \textit{some} $(\sigma, p) \in \tau$. In addition to that, there could be other $\sigma', \sigma''$ etc such that $\sigma'^G = \sigma''^G = \sigma^G$, and these additional $\mathbb{P}$-names may or may not appear inside $\tau$.} So we are starting to see the form that Comprehension$^{\mathrm{M}}$ ought to take:
$$\upsilon := \{(\sigma, p) \in \tau \:|\: \psi_{\phi}^{\mathrm{M}}(\sigma, p)\}.$$Notice that any such $\upsilon \in \mathrm{M}$ would automatically be in $\mathrm{M}^{\mathbb{P}}$ by Comprehension$^{\mathrm{M}}$. The question now is how to device a $\psi_{\phi}^{\mathrm{M}}$ that will achieve the feat of $\upsilon^G$ being equal to $\{\sigma^G \in \tau^G \:|\: \phi^{\mathrm{M}[G]}(\sigma^G)\}$.

In $\mathrm{M}[G]$, we want to collect all $\sigma^G \in \tau^G$ that satisfy $\phi$. If we could collect within $\mathrm{M}$ all $(\sigma, p) \in \tau$ such that $p \in G$ implies $\phi^{\mathrm{M}[G]}(\sigma^G)$, we would be done. So let's employ some wishful thinking and imagine we can formulate some $\psi_{\phi}^{\mathrm{M}}(\sigma, p)$ that is true in $\mathrm{M}$ iff the presence of $p$ in $G$ implies that $\phi^{\mathrm{M}[G]}(\sigma^G)$ is true.\footnote{To be perfectly clear, we want $\psi_{\phi}^{\mathrm{M}}(\sigma, p) \leftrightarrow (p \in G \rightarrow \phi^{\mathrm{M}[G]}(\sigma^G))$.}

In other words, we need to form $\psi_{\phi}^{\mathrm{M}}$ that says ``$p\in G$ \textbf{forces} $\phi^{\mathrm{M}[G]}(\sigma^G)$ to be true''.

And that is the idea at the heart of forcing! For every formula $\phi$ of ZFC, we will define a formula $\psi_{\phi}^{\mathrm{M}}$ that given $\sigma \in \mathrm{M}^{\mathbb{P}}$ specifies the set of \textbf{forcing conditions} $p\in \mathbb{P}$ whose presence in $G$ forces $\phi^{\mathrm{M}[G]}(\sigma^G)$ to be true in $\mathrm{M}[G]$. The whole point is to be able to deploy Comprehension$^\mathrm{M}$ to prove Comprehension$^{\mathrm{M}[G]}$.

\begin{quote}
\textit{The main purpose of forcing is to use Comprehension in the ground model to prove Comprehension in the outer model.}
\end{quote}

It is not at all obvious that this is possible. As it turns out, it can be done as long as $G$ is of a special kind---a \textit{generic} $G$. And generic ideals $G$ exist as long as $\mathrm{M}$ is countable, hence the need for $\mathrm{M}$ to be countable.

\subsubsection{Properties fixed by a single condition $p$ have to be generic}

Before exploring how to define formulas $\psi_{\phi}^{\mathrm{M}}(\sigma, p)$ that say ``$p$ forces $\phi^{\mathrm{M}[G]}(\sigma^G)$ to be true'', it is illuminating to discuss the implications of the ability to do so on the properties of $G$. So, suspending disbelief for a moment, assume that for any $\phi(x_1, \ldots, x_k)$, the collection of $p\in \mathbb{P}$ whose presence in $G$ guarantees the truth of sentence $\phi^{\mathrm{M}[G]}(\sigma_1^G, \ldots, \sigma_k^G)$ is a set in $\mathrm{M}$ definable with a formula $\psi_{\phi}^{\mathrm{M}}(p, \sigma_1, \ldots, \sigma_k)$. Let us for now call \textit{forcing} the ability to define such a formula $\psi_{\phi}$ for every formula $\phi$.

\begin{quote}
\textit{Let's explore this assumption down the rabbit hole of what it implies for $G$.}
\end{quote}

\textbf{Only \textit{generic} new sets can be forced into the extended model.} The subset $s_{\alpha_{\mathrm{M}}}$ discussed above that encodes the supremum $\alpha_{\mathrm{M}}$ of ordinals of $\mathrm{M}$, is specific. Can  $s_{\alpha_{\mathrm{M}}}$ possibly be in $\mathrm{M}[G]$? Well, then the formula $\phi :=$ ``there is a row $\alpha$ in $G$ that equals  $s_{\alpha_{\mathrm{M}}}$'' or something similar would need to be forced by some $p\in G$. However, $p$ is just a finite partial function. It can always be extended in each of the finite number of rows of its domain to make sure the row differs from  $s_{\alpha_{\mathrm{M}}}$. Therefore $p$ cannot ensure that any row is equal to  $s_{\alpha_{\mathrm{M}}}$.\footnote{Because $s_{\alpha_{\mathrm{M}}} \not \in \mathrm{M}[G]$, $s_{\alpha_{\mathrm{M}}}$ does not have a $\mathbb{P}$-name in $\mathrm{M}$. Therefore, $\phi :=$ ``there is a row $\alpha$ in $G$ that equals  $s_{\alpha_{\mathrm{M}}}$'' is not possible to relativize to $\mathrm{M}[G]$ and ``there is a $p$ such that $\psi_{\phi}(p, \ldots, s_{\alpha_{\mathrm{M}}})$'' saying that some $p$ forces  $s_{\alpha_{\mathrm{M}}}$ to be in $\mathrm{M}[G]$ is not possible to relativize to $\mathrm{M}$. These formulas are false in ZFC+ and meaningless for $\mathrm{M}[G]$ and $\mathrm{M}$, because they refer to terms not belonging to these models.} No $p$ can ever force any subset in $G$ to be specific.

This motivates the notion of a \textbf{generic} subset of $\omega$. It is a subset $s$ such that any of its properties $\phi(s)$ can be ensured by fixing a finite part of $s$.\footnote{More generally, in other applications of forcing it doesn't have to be a finite part of $s$, as long as it is a part of $s$ already present in the ground model $\mathrm{M}$. In our case, because $\mathbb{P}$ is defined to only have finite partial functions as elements, only finite conditions on $s$ are possible.}

If forcing is possible, then all subsets in $G$ are generic. This immediately does away with any worries that special subsets that create contradictions are introduced. Let's see some additional properties implied for generic subsets.

\noindent\textbf{Examples of generic properties.} 
\begin{itemize}
\item Consider two different rows $\alpha, \beta \in \aleph_2^{\mathrm{M}[G]}$ of $G$, corresponding to subsets $s_{\alpha}, s_{\beta}$ of $\omega$. Is it ever the case that $s_{\alpha} = s_{\beta}$? Could some $p\in G$ force them to be equal? Any $p$ only fixes a finite portion of the rows $\alpha$ and $\beta$, so it is always possible to extend $p$ so as to make $s_{\alpha} \not = s_{\beta}$. Notice the asymmetry here: whereas it is impossible for any $p$ to fix $s_{\alpha} = s_{\beta}$, many different $p$ make $s_{\alpha} \not = s_{\beta}$. No matter which $p$ we start with, we can always extend it to ensure $s_{\alpha} \not = s_{\beta}$. Therefore, because we insist that all properties of $s_{\alpha}, s_{\beta}$ are fixed by some $p$, all rows of a generic $G$ are distinct! 
\item Could $s_{\alpha}$ be missing some positions? Could we ``neglect'' to introduce some position $k$ to $s_{\alpha}$ and leave $s_{\alpha}[k]$ is undefined? Well, some $p$ would have to force $s_{\alpha}$ to be missing position $k$. However, it is always possible to extend $p$ by setting $(\alpha, k)$ to $0$ or $1$. Therefore, a generic $G$ is a maximal ideal and $\bigcup G$ is a complete function!
\item Could $s_{\alpha}$ equal some subset $t\in \mathrm{M}$? Again, no finite $p\in \mathbb{P}$ can force that because $p$ can always be extended with a position ensuring $s_{\alpha}$ differs from $t$. Therefore, a generic $G$ contains all new subsets not in $M$!
\end{itemize}

We are getting lots of freebies. All these properties are immediately implied by the yet to be demonstrated fact that every formula is forced by some element $p \in \mathbb{P}$. And all of these generic properties are highly desirable for $G$. Here is another freebie that we get:

\begin{itemize}
\item Consider the complete works of Shakespeare encoded in binary digits $\mathrm{Shake}[0] \cdots \mathrm{Shake}[N]$. Could a subset $s_{\alpha}$ fail to have some position $i$ such that $s_{\alpha}[i \::\: i+N] = \mathrm{Shake}[0 \::\: N]$? Well, some $p$ would have to force $s_{\alpha}$ to neglect Shakespeare's complete works. However, it is always possible to extend $p$ by finding some $i$ greater than any position of $s_{\alpha}$ set by $p$ and setting the next $N+1$ positions to $\mathrm{Shake}[0 \: : \: N]$. Therefore every row of a generic $G$ contains the complete works of Shakespeare (as well as every other finite piece of information in the universe) infinitely many times.\footnote{For any binary string $b$ and for any $i$ and any $\alpha \leq \aleph_2^{\mathrm{M}}$ there is a ZFC+ proof that $b$ appears in $|b|$ consecutive positions of $s_{\alpha}$ starting at a position $\geq i$, and this statement is forced by $p=\{\}$. Infinity is big and generic infinity contains everything infinitely many times.}

\end{itemize}
This is all fantastic news. Now it remains to show that forcing is possible.

\subsubsection{Necessity and possibility}

Recapping, the desire to leverage Comprehension$^\mathrm{M}$ to prove Comprehension$^{\mathrm{M}[G]}$ leads naturally to wishing that every sentence $\phi^{\mathrm{M}[G]}(\sigma^G)$ becomes true or false depending on a sentence $\psi_{\phi}^{\mathrm{M}}(\sigma, p)$ that says ``$p$ forces $\phi^{\mathrm{M}[G]}(\sigma^G)$ to be true''. Assuming this can be done, because any $p\in \mathrm{M}$ has very limited information, $G$ needs to be generic. Next, we will define genericity more precisely.

Let's fix some $\phi$ and $p$. We take the perspective of partial construction $p$ that will eventually be extended to a full $G$. Each possible extension  $q\geq p$ introduces incremental information about $G$ and $\mathrm{M}[G]$. Sometimes $p$ already has enough information to decide (to force) $\phi$ to be true or false, and other times $\phi$'s truth value may need to wait until some $q \geq p$ introduces more information.

What conditions should we impose on $p$ that force $\phi$ to be true? Some $q\geq p$ may make $\phi$ true, some $q \geq p$ may make $\phi$ false, and others may have uncertain status. Take for example $\phi :=$ ``$s_{\alpha}[n] = 0$''. Assuming $(\alpha, n) \not \in \mathrm{Dom}(p)$, some $q\geq p$ make $\phi$ true and some make $\phi$ false, while others leave it undecided. It seems intuitively clear that $p$ lacks sufficient information to decide $\phi$. As another example, $\psi := s_{\alpha} \not = s_{\beta}$. Given a $p$ that does not already have any entries $(\alpha, n), (\beta,n)$ set to different values for any $n$, some $q\geq p$ make $\psi$ true but none make $\psi$ false. Here, perhaps $p$ already has all the needed information to decide that $\psi$ is true in $\mathrm{M}[G]$.

Our first requirement is that \textit{every sentence $\phi$ of $\mathrm{M}[G]$ is decided by some $p \in G$}. Additionally, we want \textit{stability}: if $p$ decides $\phi$ to be true, then every $q\geq p$ should also decide $\phi$ to be true otherwise a contradiction could ensue.

With these two requirements in mind, let's ask what are sufficient conditions for $p$ to decide that a sentence  $\phi$ of ${\mathrm{M}[G]}$ is true. What if for every $q$ extending $p$, there is some $r$ extending $q$ such that $r$ makes $\phi$ true? In this scenario, could any $q\geq p$ decide $\lnot \phi$? No, because then some $r\geq q$ would decide $\phi$ and introduce a contradiction.

So this is it: $p$ should \textit{force} $\phi$ if for every $q\geq p$ there is $r\geq q$ such that $r$ makes $\phi$ true. If this condition is met, there is no way to find a $q\geq p$ that forces $\lnot \phi$. Because $\phi$ has to be true or false in $\mathrm{M}[G]$, $p$ is \textit{forced} (pun intended) to force $\phi$ to be true.

This is  the basic mechanism of forcing. The truth value of every sentence $\phi^{\mathrm{M}[G]}$ has to be decided by some $p \in \mathbb{P}$. For certain sentences this is easy: if $\phi$ asserts that $s_{\alpha}[n] = 1$, it is obvious how $\phi$ is decided when a $p$ is reached with $(\alpha, n)$ in its domain. However, certain other sentences may need to ``wait'' for potentially infinitely many choices to be made: for example, when $\phi$ asserts $s_{\alpha} = s_{\beta}$.\footnote{Or even when $\phi$ asserts that $(\alpha, n)$ is in the domain of $\mathbb{P}$ as there are limits $p_1 \leq p_2 \leq \cdots \leq p_n \leq \cdots$ that do not include given entries in the domain of $\bigcup_i p_i$.} If we want every sentence to be decided by a single $p$ so that Comprehension is transferred to $\mathrm{M}[G]$, we need to find a mechanism to decide such sentences. \textit{Forcing is precisely made for this purpose.}\footnote{The reader may ask: what if there is a sentence $\phi$ and a $p$ such that for every $q\geq p$ there is an $r \geq q$ making $\phi$ true, and at the same time there is an $r' \geq q$ making $\phi$ false? A quick reflection will reveal that this is impossible: if $r\geq q$ makes $\phi$ true, then there is no $r' \geq r$ making $\phi$ false. That's why forcing works.}

We define a forcing relation from ground formulas $x=y$, $x\in y$ all the way up to composite formulas, whereby $p\in \mathbb{P}$ \textit{forces} a formula $\phi$ to be true through structural induction on $\phi$ whenever for every $q\geq p$ there is $r \geq q$ that \textit{forces} $\phi$ to be true. The latter condition unrolls the structure of $\phi$ as will become clear.

\textbf{Necessity and possibility.} In modal logic there is the concept of \textbf{necessity} and its dual, \textbf{possibility}. An exposition of the independence of $\mathrm{CH}$ based on modal logic is covered in the excellent book by Smullyan and Fitting (2010). We will not need to cover modal logic here except for a basic concept at the heart of forcing. In our context specifically, necessity and possibility are taken from the point of view of a partial construction $p$. If \textit{every} $q$ extending $p$ has a given property $\phi$, then $\phi$ is \textbf{necessarily true} at $p$. For example, if $p$ already fixed $(\alpha, 0) \to 1$, then for any $q \geq p$ this value is still fixed, so $\phi := s_{\alpha}[0] = 1$ is necessarily true. If \textit{some} $q$ extending $p$ has a given property $\phi$, then $\phi$ is \textbf{possibly true} at $p$. For example, if $p$ has not fixed yet positions $1,\ldots, 2022$ of row $\omega + 2022$, it is possibly true that $\phi := \bigwedge_{i=1}^{2022} s_{\omega + 2022}[i] = 1$. The notion we actually care about is that of \textbf{necessarily possibly true}: no matter how we extend $p$ to some $q$, we can always extend $q$ to some $r$ that makes $\phi$ true. For example, at any $p$, it is necessarily possibly true that $s_{\alpha} \not = s_{\beta}$ for $0 \leq \alpha \not = \beta < \aleph_2^{\mathrm{M}}$. The importance of this notion was already given above: if $\phi$ is necessarily possibly true at $p$, then there is no way to find $q\geq p$ that forces $\lnot \phi$ and therefore $p$ is forced (pun intended) to force $\phi$.

In other words: if some sentence $\phi^{\mathrm{M}[G]}$ is necessarily possibly true at $p$, then  $p$ should force $\phi^{\mathrm{M}[G]}$. Applying this concept iteratively, if $\phi^{\mathrm{M}[G]}$ is necessarily possibly \textit{forced} at $p$, then $p$  should force $\phi^{\mathrm{M}[G]}$. Given our wishful thinking that every sentence $\phi^{\mathrm{M}[G]}$ is decided by some $p \in \mathbb{P}$, we are \textit{forced} to consider every $\phi^{\mathrm{M}[G]}$ that is necessarily possibly forced at $p$ to be forced by $p$.

Figure 4 summarizes this main idea. The way this is implemented in the proof below is with the concept of \textit{genericity} of $G$.

\begin{figure}[h!]
\begin{center}
\includegraphics[width=360pt]{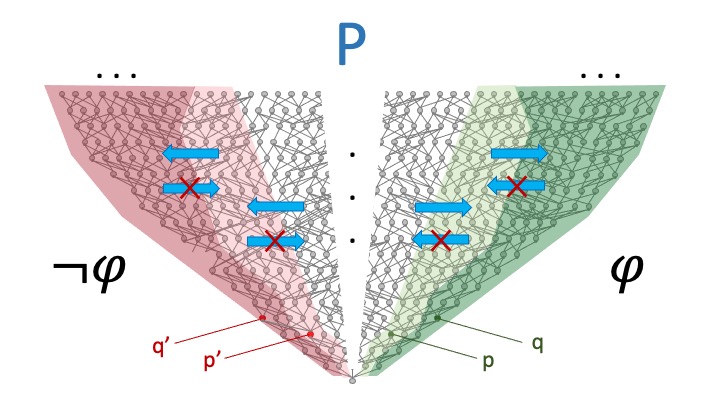}
\end{center}
\caption{\textbf{Deciding the truth value of a sentence $\phi^{\mathrm{M}[G]}$ with a condition $p\in \mathbb{P}$ .} $\mathbb{P}$ is displayed from $\{\}$ (bottom) to larger partial functions at higher levels. ($\mathbb{P}$ has infinitely many elements at each level and infinitely many edges per element.) Given a sentence $\phi$, $\mathbb{P}$ is color coded: green for true, light green for necessarily possibly true, ruby red for false, pink for necessarily possibly false, and gray for undecided. Take for example $\phi :=$ ``$(s_0[0] = 1 \land$ the value of $s_0[1]$ is set$)\lor$ the value of $s_0[1]$ is not set'', referring to row $s_0$ of $G$. Let $p$ contain $(0, 0)\to 0$. Then ``$s_0[0] = 1$'' is true for any $G\ni p$. Moreover, any $q \geq p$ that contains $(0,1)\to 0$ satisfies $\phi$. Such a $q$ is green and $p$ is light green: no matter how $p$ is extended, it is possible to extend further to some $q$ with $(0,1)$ in its domain to  satisfy $\phi$; $p$ forces $\phi$ to be true. Consider now $p'$ that contains $(0, 0)\to 1$. Now, the only way for $\phi$ to be true would be that ``the value of $s_0[1]$ is not set''. However, no matter how $p'$ is extended, it is always possible to extend it further and include $(0,1)$ in its domain, falsifying $\phi$. Therefore, $p'$ falls in the pink region. For a general $\phi$, any $p\in \mathbb{P}$ is in one of these five categories.  A $p$ that makes $\phi$ true or false can only be extended to $q$ with the same property. A $p$ that makes $\phi$ necessarily possibly true (false) can only be extended to either a $q$ with the same property or to one that makes $\phi$ true (false). $\mathbb{P}$ is colored by structural recursion on $\phi$. For theorems of ZFC, all $\mathbb{P}$ is green. Ground sentences like ``$s_0[0] = 1$'' admit only green, red and gray. Necessarily possibly true sentence like `` the value of $s_0[1]$ is set'' admit only green and light green. If $G$ is generic, every sentence $\phi^{\mathrm{M}[G]}$ has some $p\in G$ that lies in the non-gray region, as we will see.}
\end{figure}

\subsubsection{Genericity and forcing}

The condition of ``necessarily possible'', stating that for every $q\geq p$ there is some $r \geq q$ with a given property, can be defined elegantly as a graph property of $\mathbb{P}$: density of a subset $D\subseteq \mathbb{P}$, $D \in \mathrm{M}$.

\begin{quote}
\textbf{Definition of dense sets over $\mathbb{P}$.} A set $D \subseteq \mathbb{P}$ is \textbf{dense} if for all $q \in \mathbb{P}$ there is $r \geq q$ such that $r \in D$. $D$ is \textbf{dense above} $p \in \mathbb{P}$ if for all $q \geq p$ there is $r\geq q$ such that $r\in D$.
\end{quote}

To understand how the concept of density will be employed, let $D \subseteq \mathbb{P}$ be any subset, and imagine that $\phi^{\mathrm{M}[G]}(\sigma^G)$ is true for \textit{every single} $G$ that contains some element of $D$. This is a strong condition: if $p\in G$ for some $p \in D$, then $\phi^{\mathrm{M}[G]}(\sigma^G)$ is true. What is the implication if $D$ is dense? Then $\phi$ is necessarily possibly true at any $p$. If in addition $D \in \mathrm{M}$ is a set in $\mathrm{M}$ definable given parameter $\sigma \in \mathrm{M}^{\mathbb{P}}$, then we can construct formula $\psi^{\mathrm{M}}(x, \sigma)$ that forces $\phi^{\mathrm{M}[G]}(\sigma^G)$ to be true whenever $x$ is some $p \in D$.

Any $\phi^{\mathrm{M}[G]}$ that is necessarily possibly true at some $p$ (in the sense that for any $q\geq p$ there is an $r \geq q$ such that $\phi^{\mathrm{M}[G]}$ is true if $r \in G$) corresponds to a dense set $D_{p, \phi} := \{r \geq p \:|\: \phi$ is true at $r\}$. We will only care about dense sets $D\in \mathrm{M}$, because we need to be able to apply Comprehension$^\mathrm{M}$. Properties $\phi$ that correspond to dense sets $D\not\in \mathrm{M}$ will simply not relativize to $\mathrm{M}$ so we don't need to worry about them: such sentences are true or false in ZFC+ but do not have a truth value$^{\mathrm{M}[G]}$.

Some examples of dense sets $D\subseteq \mathbb{P}, D\in M$: (1) the set of functions $p$ that fix element $(\alpha, n)$ is a dense subset of $\mathbb{P}$ and is in $\mathrm{M}$ by Comprehension$^{\mathrm{M}}$;\footnote{For instance $\{p\in \mathbb{P} \:|\: \exists x \: x\in p \land x[0] = (\alpha, n)\}$ is an instance of Comprehension$^{\mathrm{M}}$ with obvious abbreviations.} (2) the set of $p$ that make rows $\alpha$ and $\beta$ different;\footnote{For instance $\{p \in \mathbb{P} \:|\: \exists n \: \exists x \: \exists y \: x \in p \land y \in p \land x[0] = (\alpha, n) \land y[0] = (\beta, n) \land x[1] \not = y[1]\}$ with obvious abbreviations.} (3) the set of $p$ that make row $\alpha$ different from a specific subset $s \subseteq \omega$ in $\mathrm{M}$; (4) the set of $p$ that make row $\alpha$ contain the binary encoding of Shakespeare's complete works in some consecutive positions; (5) the set of $p$ that are intersections of finitely many of the above sets: finitely many dense properties can be combined to form a dense property that is the $\land$ of all of them. Importantly, this allows combination of dense properties within composite formulas and definition of forcing with structural induction on formulas.

At last, we are ready to mathematically define what a generic $G$ is. A generic $G$ should never fail to intersect a set $D \in M$ that is dense over $\mathbb{P}$. The formula $\Gamma \cap \check{D} \not = \check{\emptyset}$ evaluates to $G \cap D \not = \emptyset$ in $\mathrm{M}[G]$, and is necessarily possibly true at any $p$. Therefore every $p$ should force it and consequently there should be some $q\in G\cap D$. In other words, given $D\in \mathrm{M}$ dense over $\mathbb{P}$,  the $\mathrm{M}[G]$-sentence $G \cap D \neq 0$ is either true or false. By our premise (wishful thinking) that every sentence of $\mathrm{M}[G]$ is forced by some $p\in \mathbb{P}$, either some $p$ forces $\Gamma \cap \check{D} \not = \check{\emptyset}$, or some $p$ forces $\Gamma \cap \check{D} = \check{\emptyset}$. But no $p$ can do the latter, therefore every $p$ has to do the former. \textit{To ensure that, we define genericity thus:}

\begin{quote}
\textbf{Definition of genericity.} An ideal $G$ over $\mathbb{P}$ is \textbf{generic} over $M$ if it intersects every dense set $D \in M$. That is, $G \cap D \neq \emptyset$ for every dense $D\subset \mathbb{P}, D\in \mathrm{M}$. 
\end{quote}

As we will see, this simple condition elegantly enables defining a forcing formula $\psi_{\phi}^{\mathrm{M}}(p, \ldots)$ for every sentence $\phi^{\mathrm{M}[G]}$. There is a big question though, does such a $G$ exist? Generally no, but if the ground model is countable then an easy proof guarantees the existence of a generic $G$. This is precisely why a \textit{countable} transitive model $\mathrm{M}$ is postulated.
\begin{quote}
\textbf{Proposition: existence of a generic ideal.} Let $\mathbb{P}$ be a poset over countable model $M$ and let $p\in \mathbb{P}$. A generic ideal $G$ with $p\in G$ exists. 
\end{quote}

\textbf{Proof.} Let $D_1, D_2, D_3, \ldots$ be an enumeration of all dense subsets of $\mathbb{P}$ that are members of countable model $\mathrm{M}$. We will create a sequence $p_0, p_1, p_2, \ldots$ hitting each one of them. Let $p_0 = p$. Since $D_1$ is dense, there is $p_1\in D_1$ with $p_1\geq p_0$. Similarly, $p_{i+1}\in D_{i+1}$ with $p_{i+1} \geq p_i$ can always be found.

Now, we will create an ideal $G$ that contains all elements in that sequence. $G$ is simply $\{p \in \mathbb{P} \:|\: p \leq p_i$ for some $p_i$ in the sequence$\}$. Clearly $p_i \in G\cap D_i \not = 0$, so it remains to show that $G$ is an ideal. Indeed it is downward closed by definition, and for any two $p, q \in G$, $p \leq p_i, q\leq p_j$ for some $i, j$ and then pick $k \geq i, j$ and we have $p, q \leq p_k$. Therefore every pair of elements of $G$ are compatible. $\square$

\begin{quote}
\textit{The existence of generic ideals is precisely the reason $\mathrm{M}$ needs to be countable. $G$ needs to hit dense subsets that are members of $\mathrm{M}$, of which there are only countably many.} What about dense subsets not in $\mathrm{M}$? We don't care about them because they correspond to $\phi$ with constants outside $\mathrm{M}[G]$ and therefore whose truth is settled outside $\mathrm{M}[G]$. 
\end{quote}

There is one additional property that $\mathbb{P}$ needs to have to ensure that $G\not \in \mathrm{M}$: no condition $p$ should be such that it ``finalizes'' $G$. Formally, $p$ is an \textbf{atom} if every pair $q, r \geq p$ are compatible. Once an atom is hit, $G$ is finalized: no more ``choices'' need to be made. $\mathbb{P}$ is \textbf{atomless}, if it has no atoms. In addition, if $p \leq q \leq p$ implies $p = q$, then $\mathbb{P}$ is called \textbf{separative}. For forcing to go through and generate $\mathrm{M}[G] \supsetneq \mathrm{M}$, $\mathbb{P}$ has to be separative. The $\mathbb{P}$ we have been discussing so far is clearly separative.

\textbf{Proposition: $G$ is not in $\mathrm{M}$.} If $\mathbb{P}$ is a separative poset and $G$ is generic over $\mathbb{P}$, then $G\not \in \mathrm{M}$.

\textbf{Proof.} Let $H \subseteq \mathbb{P}$, $H \in \mathrm{M}$. Consider the set $D_H := \{p\in \mathbb{P} \:|\: \exists q \in H \: p \perp q \}$. $D_H \in \mathrm{M}$ by Comprehension$^{\mathrm{M}}$, and it is easy to see that $D_H$ is dense. Therefore $G \cap D_H \not = \emptyset$. $\square$

With the definition of genericity in place, we are finally ready to define forcing.

\begin{quote}
\textbf{Definition: forcing.} An element $p \in \mathbb{P}$ \textbf{forces} the sentence $\phi(\tau_1, \ldots, \tau_k)$ if $\phi^{\mathrm{M}[G]}(\tau_1^G, \ldots, \tau_k^G)$  is true for every generic ideal $G$ with $p\in G$. Forcing establishes a separate relation $\subset \mathbb{P}\times (\mathrm{M}^{\mathbb{P}})^k$ for every $k$-ary formula $\phi$ in the language, denoted by $p \Vdash \phi(\tau_1, \ldots, \tau_k)$.
\end{quote}

Forcing is a generalization of truth in $\mathrm{M}$, or the notion of satisfaction: if $p\Vdash\phi$ and $\phi \rightarrow \psi$, then $p \Vdash \psi$. Notice that forcing is only defined on the \textit{sentences}. Forcing establishes a separate relation among elements of the ground model $\mathrm{M}$ for every sentence $\phi$ of ZFC. Given a formula of $k$ variables and a $k$-tuple of names, the relation defines the set of $p$ that make the resulting sentence true in ${\mathrm{M}[G]}$.

\begin{quote}
\textbf{Notation.} The standard notation $p \Vdash \phi(\tau_1, \ldots, \tau_k)$ when $\phi^{\mathrm{M}[G]}(\tau_1^G, \ldots, \tau_k^G)$ is true for every generic $G \ni p$ can be somewhat confusing. Nobody actually cares about the truth value of $\phi(\tau_1, \ldots, \tau_k)$, only of $\phi^{\mathrm{M}[G]}(\tau_1^G, \ldots, \tau_k^G)$. For example, $\phi := \check{0} \in \check{1}$ is forced by every $p$ because $\phi^{\mathrm{M}[G]} := \check{0}^G \in \check{1}^G :=  0 \in 1$ is true for any generic $G$. However, $\check{0} \in \check{1}$ is obviously false and its truth value is irrelevant. The reason for the notation is that the forcing relation for $\phi$ is between $p$, $\tau_1, \ldots, \tau_k$, and $G$ is \textit{not} part of the relation. In fact, the relation is definable within $\mathrm{M}$.
\end{quote}

The next task will be to prove that given any $\phi$ the forcing relation is definable in $\mathrm{M}$, so that the collection of $p$ that force $\phi$ is a set in $\mathrm{M}$. Once this is established,  Comprehension$^{\mathrm{M}}$ can be leveraged to assign truth within $\mathrm{M}[G]$. The key to showing definability of forcing is that $G$ is generic. 

\subsection{Proof that $\mathrm{M}[G]$ is a model}

The goal now is to show that $\mathrm{M}[G]$ is a model of ZFC. This is a metatheorem: for every axiom $\phi$ of ZFC, an algorithm produces a derivation of $\phi^{\mathrm{M}[G]}$ from  $\phi^{\mathrm{M}}$. Going from $\phi^{\mathrm{M}}$ to $\phi^{\mathrm{M}[G]}$ is not a simple task. It requires careful control of what statements are true in $\mathrm{M}[G]$, and the key to that is the fundamental theorem of forcing, which is a metatheorem \textit{about} ZFC rather than a theorem \textit{of} ZFC.

\subsubsection{Statement of the ``fundamental theorem of forcing'' metatheorem}

Given a formula $\phi(x_1, \ldots, x_k)$, the fundamental theorem of forcing (FTF) says that the collection of $p \in \mathbb{P}$ and $\mathbb{P}$-names $\tau_1, \ldots, \tau_k$ that force $\phi$ is a set in $\mathrm{M}$, and moreover that $\phi^{\mathrm{M}[G]}(\tau_1^G, \ldots, \tau_k^G)$ is true iff some $p$ forces $\phi(\tau_1, \ldots, \tau_k)$. 

\begin{quote}
\textbf{The fundamental theorem of forcing (FTF).} Let $\phi(x_1, \ldots, x_k)$ be a formula of $k$ free variables, $\mathbb{P}$ a separative poset in countable transitive set model $\mathrm{M}$, and $G$ a generic ideal of $\mathbb{P}$. Let $\tau_i, 1\leq i\leq k$ be $\mathbb{P}$-names of maximum name rank $\alpha$.
\begin{enumerate}
\item \textbf{Definability.} There is a set $\mathcal{F}^{\phi}_{\alpha} \in \mathrm{M}$ such that $p \Vdash \phi(\tau_1, \ldots, \tau_k)$ iff $(p, \tau_1, \ldots, \tau_k) \in \mathcal{F}^{\phi}_{\alpha}$. That is to say, \textbf{the forcing relation is definable} in $\mathrm{M}$.
\item \textbf{Truth.} For any generic ideal $G$ over $\mathbb{P}$, $\phi^{\mathrm{M}[G]}(\tau_1^G, \ldots, \tau_k^G)$ iff there is a $p\in G$ such that $p \Vdash \phi(\tau_1, \ldots, \tau_k)$. In other words, in $\mathrm{M}[G]$ \textbf{any true sentence is forced by some} $p\in G$.
\item \textbf{Coherence.} If $p \Vdash \phi$ and $q \geq p$ then $q \Vdash \phi$.
\end{enumerate}
\end{quote}

The FTF is an \textit{algorithm on ZFC syntax}: it takes a formula $\phi$ as input and produces a ZFC proof of the three assertions for $\phi$. 

\textbf{How FTF is used in the proof of $\mathrm{CH}^{\mathrm{M}[G]}$.} In the algorithm that transforms a proof of $\mathrm{CH}$ into a proof of $\mathrm{CH}^{\mathrm{M}[G]}$, as a first step all the needed axioms $\phi_1, \ldots, \phi_k$ are collected. The Reflection Principle algorithm provides a ZFC proof of existence of a countable transitive set $M$ in which every $\phi_i^{M}$ holds. Then the FTF algorithm generates a separate ZFC proof of FTF for each $\phi_i^M$, including recursively for its subformulas, to obtain $\phi_i^{M[G]}$. Then, the proof of $\mathrm{CH}$ is copied to a proof of $\mathrm{CH}^{M[G]}$.

Before proving FTF, a few words about how it will be applied to prove Comprehension$^{\mathrm{M}[G]}$ using Comprehension$^{\mathrm{M}}$. Given a formula $\phi$, whenever $\phi^{\mathrm{M}[G]}(\sigma^G)$ is true there is a $p$ that forces $\phi(\sigma)$ (FTF.Truth). Looking at the elements $(\sigma,p)$ of $\tau$, the maximum name rank $\mathrm{nr}(\sigma)$ is some $\alpha$, and now $\mathcal{F}_{\alpha}^{\phi}$ contains all possible $(\sigma, p)$ in $\tau$ such that $p$ forces $\phi(\sigma)$ (FTF.Definability). Then, by Comprehension$^{\mathrm{M}}$ the elements $(\sigma, p) \in \tau$ such that $(p,\sigma) \in \mathcal{F}_{\alpha}^{\phi}$ form a set in $\mathrm{M}$, call it $\upsilon := \{(\sigma, p) \in \tau \: | \: (p, \sigma) \in \mathcal{F}_{\alpha}^{\phi}\}$. This is a $\mathbb{P}$-name, and it is easy to show that in fact $\upsilon^G$ is the desired set in $\mathrm{M}[G]$ of all $\{ \sigma^G \in \tau^G \: | \: \phi^{\mathrm{M}[G]}(\sigma^G)\}$. Therefore Comprehension works in $\mathrm{M}[G]$ by virtue of the FTF.

The above argument, which omits some  technical details of why $\upsilon^G$ works and is presented in more detail below, is at the heart of the proof that $\mathrm{M}[G]$ is a model. The FTF's main job is to transfer Comprehension as well as Power and Replacement from $\mathrm{M}$ to $\mathrm{M}[G]$.\footnote{Another important job comes later, in showing preservation of cardinals in $\mathrm{M}[G]$.}

The next task is to prove FTF. For starters, FTF.Coherence is trivial and follows from the downward closure of an ideal $G$. The other two clauses of FTF are quite technical to prove.

\subsubsection{Architecture of the proof of FTF}

The FTF for a formula $\phi$ has a ZFC+ proof by structural induction on $\phi$. The ground case ``$=$'' is the hardest. Each case follows the same proof layout.

\noindent\textbf{FTF proof layout.}
\begin{enumerate}

\item A general lemma establishes that for any $Z\subset \mathbb{P}, Z\in \mathrm{M}$, if $G\subset Z$ then there is some $p\in G$ all of whose extensions are in $Z$. Think of $Z$ as a net definable in $\mathrm{M}$; $Z$ will in fact be defined separately for every $\phi$ and contain all $q\in \mathbb{P}$ for which $\phi$ is possibly true. If $G$ is caught inside $Z$, the $p$ in the lemma makes $\phi$ necessarily possibly true and thereby forces $\phi$. The existence of a $p$ all of whose extensions are in $Z$ is implied because $G$ is generic---\textit{this is precisely how genericity of $G$ is used in the proof of FTF.}
\item The set $\mathcal{F}_{\alpha}^{\phi}$ contains all $(p, \tau_1, \ldots, \tau_k)$ such that  $\tau_i$'s ranks are $\leq \alpha$ and $\phi^{\mathrm{M}[G]}(\tau_1^G, \ldots, \tau_k^G)$ if $p\in \mathbb{P}$. These sets are defined with a single formula $Force^{\phi}(\alpha, p, \tau_1, \ldots, \tau_k)$. If $\phi$ is a composite sentence, then $Force^{\phi}$ is composed using existing formulas for the components of $\phi$ by structural induction. 
\item A net $Z^{\phi} \subset \mathbb{P}$ is defined that contains all $q$ for which $\phi$ is possibly true. Importantly, $Z^{\phi}$ is proven to be in $\mathrm{M}$: by induction on name rank in the ground cases $\phi:=$ ``$\tau_1 = \tau_2$'' and ``$\sigma \in \tau$", and by structural induction in the composite cases.
\item Given $\phi$ true in $\mathrm{M}[G]$, it is shown that $G\subset Z^{\phi}$, therefore some $p\in G$ has all its extensions in $Z^{\phi}$, which implies that $\phi$ is necessarily possibly true at $p$, which means $p$ forces $\phi$. The converse is also shown.
\item The FTF for $\phi$ easily follows.
\end{enumerate}

The proof is long but not difficult. Its main difficulty lies in  keeping the proof structure in mind when working on the details, and paying attention to what is definable in $\mathrm{M}$ and what is not. Let's prove the important lemma at the top. 

\textbf{Lemma: the net $Z$.} If generic $G \subset Z \in \mathrm{M}$, then there is $p \in G$ all of whose extensions are in $Z$. 

\begin{quote} \textit{This lemma is the place where genericity of $G$ is utilized.} $Z\in \mathrm{M}$ will contain all $q$ for which a condition that ensures $\phi$ is possibly satisfied. The lemma asserts a $p\in G$ for which the condition is necessarily possible, and consequently a $p$ that forces the condition. Figure 5 displays the lemma.
\end{quote}

\begin{figure}
\centering
\includegraphics[width=360pt]{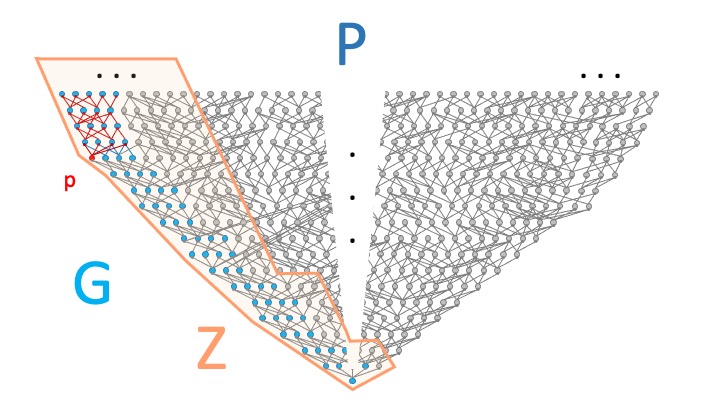}
\caption{\textbf{The net $Z \in \mathrm{M}$ that envelopes $G$.} If $G$ (green elements is a subset of  $Z \in \mathrm{M}$, then there is $p\in G$ (pictured in red here) all of whose extensions are in $Z$. The FTF will be proven by structural induction on $\phi$; each case will involve defining a net $Z$ that includes all $q$ for which $\phi$ is possibly true or possibly forced. This lemma establishes that if all of $G$ is a subset of $Z$, then there is a $p$ for which $\phi$ is necessarily possibly true / forced, and therefore that $p$ will force $\phi$ in each case.}
\end{figure}

\textbf{Proof.} Since $Z\in M$, $\bar{Z} := \mathbb{P} \setminus Z \in \mathrm{M}$. Let $D$ be the set of all elements $p\in \mathbb{P}$ that are either $\geq$ some  $q\in \bar{Z}$, or are incompatible with every  $q \in \bar{Z}$. That is, $D := \{p\in \mathbb{P} \:|\: (\exists q \in \bar{Z}\: p\geq q) \lor (\forall q\in \bar{Z} \: p \perp q) \}$. $D \in \mathrm{M}$ by Comprehension$^\mathrm{M}$ and is dense: let $p\in \mathbb{P}$; if $p$ is compatible with some $q\in \bar{Z}$ then a common extension exists in $D$, whereas if $p$ is incompatible with every $q\in \bar{Z}$ then it is in $D$ already. Therefore by genericity $G$ intersects $D$, say $r \in G\cap D$. Either $r$ extends some $q \in \bar{Z}$, which is impossible because  downward closure would imply $q\in G$ whereas  $G\subset Z$, or $r$ is incompatible with every $q\in \bar{Z}$, therefore all of $r$'s extensions are in $Z$. $\square$

\subsubsection{Proof of FTF for equality}

Following the structure of the FTF proof laid above, let's define $\mathcal{F}_{\alpha}^{x=y}$.  This is the hardest case. Equality of $\tau_1$ and $\tau_2$ in $\mathrm{M}[G]$ depends on whether each of the elements of one are in the other, and those elements in turn may become equal or unequal or even disappear depending on $G$. We saw a simple example above of how this happens. As it turns out, equality has to be settled in a transfinite induction on name rank. Same for membership, $\in$, but FTF for equality can be used in this case. The rest of the cases are settled with structural induction on $\phi$ and are much easier.

We want all triplets $(p, \tau_1, \tau_2)$ in $\mathrm{M}$ where $\mathrm{nr}(\tau_1), \mathrm{nr}(\tau_2) \leq \alpha$, such that $p$ forces $\tau_1^G = \tau_2^G$. Therefore, $p$ must be such that it is necessarily possible that  every element $\sigma^G \in \tau_1^G$ is also $\in \tau_2^G$, and vice versa. Let's start with the forward direction: whenever $q \geq p$ and $(\sigma, q) \in \tau_1$, there must be $r \geq q$ making $\sigma$ a member of $\tau_2$ (recall that $(\sigma, q) \in \tau_1$ implies that $(\sigma, r) \in \tau_1$ for any $r \geq q$). Here we have a difficulty. We may be tempted to say that $(\sigma,r )$ has to be in $\tau_2$, but as explained above on ``how $\in$ works within $\mathrm{M}[G]$'', this is too strong a condition. Perhaps $(\sigma, r)$ is not in $\tau_2$ for any $r \geq q$, but instead $(\sigma', r) \in \tau_2$ and simultaneously $r$ forces $\sigma = \sigma'$.  

So now we seem to arrive at circular reasoning. We need to define equality, and the only way to do that is through equality. Fortunately, $(\sigma, p) \in \tau$ always implies that the name rank of $\sigma$ is lower than the name rank of $\tau$. Therefore we can use transfinite induction on name ranks. Without further ado, the correct condition for equality triplet $(p, \tau_1, \tau_2)$ of name ranks $\leq \alpha$ is:
\begin{itemize}

\item For every $(\sigma, q) \in \tau_1$ with $q \geq p$ there exists $(\sigma', r) \in \tau_2$ with $r \geq q$ and $(r, \sigma, \sigma') \in \mathcal{F}_{\max\{\mathrm{nr}(\sigma), \mathrm{nr}(\sigma')\}}^{x=y}$.
\item Conversely, for every $(\sigma, q) \in \tau_2$ with $q \geq p$ there exists $(\sigma', r) \in \tau_1$ with $r \geq q$ and $(r, \sigma, \sigma') \in \mathcal{F}_{\max\{\mathrm{nr}(\sigma), \mathrm{nr}(\sigma')\}}^{x=y}$.
\end{itemize}

\textbf{A ZFC formula for forcing equality.} Similar to the definitions of formula $V(a,x)$ above, a single formula $Force^=(\alpha, p, \tau_1, \tau_2)$ is defined that is true whenever $p \in \mathbb{P}$, $\tau_1, \tau_2 \in N_{\alpha}$, and $p \Vdash \tau_1 = \tau_2$. Importantly, \textit{this formula is already relativized to $\mathrm{M}$}. Treating $\mathrm{M}$ for a moment as a model and not as mere syntax, $\mathrm{M}$ ``knows'' whether a given $p$ forces $\tau_1^G = \tau_2^G$ in $\mathrm{M}[G]$: it is whenever in every extension $q$ of $p$, every member $\sigma$ of $\tau_1$ (i.e., $(\sigma, q) \in \tau_1$) is possibly equal to a member $\sigma'$ of $\tau_2$ (i.e., $\exists r \geq q \: (\sigma',r)\in \tau_2 \land (r, \sigma, \sigma') \in \mathcal{F}_{\beta}^{x=y}$ for some $\beta < \alpha$), and vise versa. To bring this important point forward, I include an ``implementation'' of $\mathrm{Force}^=(\alpha, p, \tau_1, \tau_2)$ here:\footnote{A few obvious abbreviations are used: $\forall (\sigma, q) \in \tau_2 \ldots$ is not proper syntax, and it stands for $\forall x (x \in \tau_2 \rightarrow \exists \sigma \: \exists q \: x = (\sigma, q)\land...)$.}
$$\mathrm{Force}^=(\alpha, p, \tau_1, \tau_2) := \mathrm{Ord}(\alpha) \land p\in \mathbb{P} \land N(\alpha, \tau_1) \land N(\alpha, \tau_2) \land \exists \mathcal{F}_{aux} \in \mathrm{M} \: \mathrm{Fun}(\mathcal{F}_{aux}) \land \mathrm{Dom}(\mathcal{F}_{aux}, \alpha+1) \land $$$$\forall \beta \in \alpha + 1 \: (\forall p' \: \forall \tau_1' \: \forall \tau_2' \: (p' \in \mathbb{P} \land N(\beta, \tau_1') \land N(\beta, \tau_2')) \rightarrow ((q, \tau_1', \tau_2') \in \mathcal{F}_{aux}[\beta] \leftrightarrow$$
$$(\forall (\sigma, q) \in \tau_1' \: (q \geq p' \rightarrow \exists \sigma' \: \exists r \: (\sigma', r) \in \tau_2' \land r \geq q \land \exists \delta \in \beta \:  (r, \sigma, \sigma') \in \mathcal{F}_{aux}[\delta]) \land$$
$$(\forall (\sigma, q) \in \tau_2' \: (q \geq p' \rightarrow \exists \sigma' \: \exists r \: (\sigma', r) \in \tau_1' \land r \geq q \land \exists \delta \in \beta \: (r, \sigma, \sigma') \in \mathcal{F}_{aux}[\delta])) \land$$
$$(p, \tau_1, \tau_2) \in \mathcal{F}_{aux}[\alpha]$$
In this formula, $\mathbb{P}$ is assumed to be the poset in $\mathrm{M}$ with respect to which $N()$ is the name space within $\mathrm{M}$. Assuming $\alpha \leq \alpha_{\mathrm{M}}$ and $p, \tau_1, \tau_2 \in \mathrm{M}$, as mentioned above the entire formula is relativized to $\mathrm{M}$.\footnote{Syntactically it is not, but replacing syntax such as $\forall \beta \in \alpha+1...$ with $\forall \beta \in \mathrm{M} \: \beta \in \alpha+1 \rightarrow ...$ makes no difference as $\beta \in \alpha+1$ implies $\beta\in \mathrm{M}$ and similarly for all the other quantifiers.} 

\begin{quote}
\textbf{How Force$^=$ works.} This formula resembles computer code. For any $\alpha$, it defines the set $\mathcal{F}_{\alpha}^{x=y}\in \mathrm{M}$ of triplets $(p, \tau_1, \tau_2)$ of name ranks $\leq \alpha$ whereby $\tau_1, \tau_2$ are equal if $p \in G$. It uses an auxiliary transfinite array $\mathcal{F}_{aux}$ whose job is to store all intermediate sets $\mathcal{F}_{aux}[\beta] := \mathcal{F}_{\beta}^{x=y}$. Each of these sets uses the information in previous sets (i.e., $\mathcal{F}_{\delta}^{x=y}$ for $\delta < \beta$) to define equality for name ranks up to $\beta$. At each level $\beta$, for any $\tau_1', \tau_2' \in N_{\beta}$, $p'$ forces them equal if among the extensions of $p'$ each member of $\tau_1'$ is necessarily possibly a member of $\tau_2'$ and conversely. Here, ``member of $\tau_1'$'' means $\sigma$ such that $(\sigma, q)\in \tau_1$ for $q \geq p$,  and ``$\sigma$ is possibly a member of $\tau_2$'' means that there exists $r\geq q$ such that $r$ forces $\sigma = \sigma'$ and $(\sigma', r) \in \tau_2$. Having filled the entire array $\mathcal{F}_{aux}$, the last level $\alpha$ defines the elements $p$ that force equality of input names $\tau_1, \tau_2$. 

Because there are no descending chains of membership $\cdots \in \sigma_3^G \in \sigma_2^G \in \sigma_1^G \in \tau$ in ZFC, all the chains of forcing equality of $\tau_1, \tau_2$ at $p$ by requiring all members $\sigma_1$ of one to be necessarily possibly equal to some $\sigma_1'$ of the other are resolved in a finite chain of ``necessarily possibly": $\sigma_1 = \sigma_1'$ will be resolved by requiring in turn for all members $\sigma_2$ of one to be necessarily possibly equal to some $\sigma_2'$ of the other and so on. Such a chain creates a chain of membership $\sigma_i^G \in \sigma_{i-1}^G \in \cdots \in \sigma_1^G \in \tau_1$ and $i$ cannot go to infinity because of Foundation. Each of the successive levels $i$ down also requires successively increasing $p \leq q_1 \leq r_1 \leq q_2 \leq r_2 \leq \cdots \leq q_i \leq r_i$. Because the chain is finite, the final forcing condition (say, $r_i$) is a member of $\mathbb{P}$ rather than some kind of infinite limit.
\end{quote}

Now, using Comprehension$^\mathrm{M}$ each $\mathcal{F}_{\alpha}^{x=y}$ can be defined as set of the triplets $(p, \tau_1, \tau_2) \in \mathbb{P} \times N_{\alpha} \times N_{\alpha} \in \mathrm{M}$ that satisfy $Force^=(\alpha, p, \tau_1, \tau_2)$. 
$$\mathcal{F}_{\alpha}^{x=y} := \{x\in \mathbb{P}\times N_{\alpha}\times N_{\alpha} \: |\: \exists p \exists \tau_1 \exists \tau_2 \: x = (p, \tau_1, \tau_2) \land Force^=(\alpha, p, \tau_1, \tau_2) \}\in \mathrm{M}$$
Therefore, each $\mathcal{F}_{\alpha}^{x=y}$ is a set in $\mathrm{M}$.

Once $\mathcal{F}_{\alpha}^{x=y}$ is defined, the proof below establishes FTF for $\phi := x=y$. Specifically, two elements $\tau_1^G$ and $\tau_2^G$ in $\mathrm{M}[G]$ are equal iff some $p\Vdash \tau_1 = \tau_2$ iff $(p, \tau_1, \tau_2) \in \mathcal{F}_{\alpha}^{x=y}$. This is shown by transfinite induction on name rank $\alpha$, assuming both parts of the theorem are true for all $\beta < \alpha$, specifically for all $(\sigma, p)$ in $\tau_1$ and $\tau_2$. Because this is the hardest part of the proof of FTF, it is included here. Much of this proof is mere bookkeeping. 

\textbf{Proof of FTF for equality.} Following the structure of the FTF proof the next step is to define the net $Z$. For equality, let $Z^=:= \{q\in \mathbb{P} \:|\:$ for all $(\sigma_1, q) \in \tau_1$ there is $(\sigma_2, r) \in \tau_2$ with $r\geq q$ such that $r \Vdash \sigma_1 = \sigma_2\}$. Notice that $Z^{=}$ is the set of $q$ for which $\tau_1$ is \textit{possibly a subset of or equal to} $\tau_2$. 

By inductive assumption on name rank, FTF holds for $\sigma_1, \sigma_2$, and it follows that $Z^= \in \mathrm{M}$. Now let's go ahead and perform step 4 of the FTF proof layout above:

\begin{quote}
\textbf{Claim}. $\tau_1^G \subseteq \tau_2^G$ implies $G \subset Z^=$ implies $\exists p \in G$ $\forall q \geq p$ $q \in Z^=$ implies  $\tau_1^G \subseteq \tau_2^G$. Therefore all three conditions are equivalent.
\end{quote}

We will prove that $\tau_1^G \subseteq \tau_2^G$ implies $G \subset Z^=$, which by the net lemma implies there is $p\in G$ all of whose extensions are in $Z^=$. Then we will prove that if there is such a $p$, then $\tau_1^G \subseteq \tau_2^G$, establishing that all these three conditions are equivalent. Conversely the same argument holds for $\tau_2^G \subseteq \tau_1^G$. Then, FTF will follow easily.

\textbf{Proof of claim.} Let $\tau_1^G\subseteq \tau_2^G$. Let $p\in G$, to show that $p \in Z^=$. If there is no $(\sigma_1, p) \in \tau_1$, then $p$ vacuously satisfies the condition to be in $Z^=$. So let $(\sigma_1, p) \in \tau_1$, therefore $\sigma_1^G \in \tau_1^G$. Then there is $(\sigma_2, q) \in \tau_2$ with $q \in G$ and $\sigma_1^G = \sigma_2^G$. The name ranks of $\sigma_1, \sigma_2$ are smaller than those of $\tau_1, \tau_2$ so by  induction we can assume FTF for them, so there is $r \in G$ such that $r \Vdash \sigma_1 = \sigma_2$. Now a common extension $s \geq p, q, r$ satisfies $(\sigma_1, s) \in \tau_1$, $(\sigma_2, s) \in \tau_2$, and $s \Vdash \sigma_1 = \sigma_2$. Therefore $p \in Z^=$.

Since $p$ was arbitrary, $G \subset Z^=$, so by the lemma there is $p\in G$ such that every $q\geq p$ is in $Z^=$. Now to show that this in turn implies $\tau_1^G \subseteq \tau_2^G$.

Let $p$ be such. Let $\sigma_1^G \in \tau_1^G$, to show $\sigma_1^G \in \tau_2^G$.  There is some $q\in G$ with $(\sigma_1, q) \in \tau_1$ and we can assume $q\geq p$. Therefore every extension $q'\geq q$ is in $Z^=$ and has $(\sigma_1, q') \in \tau_1$, therefore there is $r\geq q'$ with $(\sigma_2, r) \in \tau_2$ and $r \Vdash \sigma_1 = \sigma_2$. Therefore the set $D_q$ of $r$ with $(\sigma_2, r) \in \tau_2$ and $r \Vdash \sigma_1 = \sigma_2$ is dense above $q$. By assuming FTF.Definability inductively on name rank, $D_q\in \mathrm{M}$ (because $\mathcal{F}_{max\{\mathrm{nr}(\sigma_1), \mathrm{nr}(\sigma_2)\}}^{x=y} \in \mathrm{M}$ by inductive assumption). Since $D_q$ is dense above $q\in G$ there is some $r\in G\cap D_q$ for which $(\sigma_1, r) \in \tau_1, (\sigma_2, r) \in \tau_2$, and $r \Vdash \sigma_1 = \sigma_2$ therefore $\sigma_1^G = \sigma_2^G \in \tau_2^G$. Since $\sigma_1^G \in \tau_1^G$ was arbitrary, the conclusion is that $\tau_1^G \subseteq \tau_2^G$ and moreover that $p \Vdash \tau_1 \subseteq \tau_2$. This concludes the proof of the claim. $\square$

Now it remains to prove two parts of FTF: (1) $p \Vdash \tau_1 = \tau_2$ iff $(p, \tau_1, \ldots, \tau_k) \in \mathcal{F}^{x=y}_{\alpha}$; (2) $\tau_1^G = \tau_2^G$ iff there is $p \in G$ such that $p \Vdash \tau_1 = \tau_2$.

For (1), let $p \Vdash \tau_1 = \tau_2$. We need to show that $(p, \tau_1, \tau_2) \in \mathcal{F}_{\alpha}^{=}$. The conditions of membership $\mathcal{F}_{\alpha}^=$ say that every extension $q\geq p$ must belong to $Z^=$ and conversely with $\tau_1$ and $\tau_2$ flipped. For any $q\geq p$ we can find $G$ containing $q$, and we have just shown that $G\subset Z^=$ therefore $q \in Z^=$. The direction with $\tau_1$ and $\tau_2$ flipped is exactly the same, so  $(p, \tau_1, \tau_2) \in \mathcal{F}_{\alpha}^{=}$. Conversely, let  $(p, \tau_1, \tau_2) \in \mathcal{F}_{\alpha}^{=}$ and let $G$ be a generic ideal that contains $p$. By definition of $\mathcal{F}_{\alpha}^{=}$ every extension $q\geq p$ is in $Z^=$, therefore as we just showed $\tau_1^G \subseteq \tau_2^G$. The same argument with $\tau_1$ and $\tau_2$ flipped, establishes $\tau_1^G = \tau_2^G$.

For (2), the reverse direction is trivial. So let $\tau_1^G = \tau_2^G$, therefore $G \subset Z^=\in \mathrm{M}$ therefore by the above lemma we can find $p\in G$ such that all $q\geq p$ are in $Z^=$, therefore as we showed $p \Vdash \tau_1 \subseteq \tau_2$ and by flipping $\tau_1$ and $\tau_2$ and repeating the argument including the definition of the corresponding $Z^=$, $p \Vdash \tau_1 = \tau_2$. $\square$

\subsubsection{Proof of the fundamental theorem of forcing for the remaining cases}

Equality was the hardest case of FTF to prove. The rest of the cases are easier and follow the same proof layout. Let's first define  the sets $\mathcal{F}_{\alpha}^{\phi}$ for each structural case of $\phi$. We do not need to be elaborate---they are all similar to $Force^=$ so we will just sketch them. The reader can write down some of these formulas for fun.

$\mathcal{F}^{x \in y}_{\alpha}$ contains the triplets $(p, \tau_1, \tau_2)$ such that for every $q \geq p$ there is an $r \geq q$ and $(\sigma, r) \in \tau_2$ with $r \Vdash \tau_1 = \sigma$. In other words, such that $\tau_1$ is necessarily possibly a member of $\tau_2$. (Notice that $\mathrm{nr}(\tau_1) < \alpha$ otherwise $(p, \tau_1, \tau_2)$ is never in $\mathcal{F}_{\alpha}^{x\in y}$.) The formula  $Force^{\in}$ uses a transfinite induction with $Force^=$ as a subroutine, and details are omitted. Alternatively, a combined  formula $Force^{= \: or \: \in}(\alpha, p, \tau_1, \tau_2, T)$ can be written down, where $T = 0/1$ is the type ($0:$ $=$; $1:$ $\in$), and two auxiliary arrays are defined, $\mathcal{F}_{aux}^=, \mathcal{F}_{aux}^{\in}$ holding the respective sets of $(p, \tau_1, \tau_2)$. The resulting formula is quite long but not hard to write down and details are omitted.

$\mathcal{F}^{\lnot \psi}_{\alpha}$ contains the tuples $(p, \tau_1, \ldots, \tau_k)$ such that no $q\geq p$ forces $\psi(\tau_1, \ldots, \tau_k)$:  no $q\geq p$ exists with $(q, \tau_1, \ldots, \tau_k)$ a member of  $\mathcal{F}_{\alpha}^{\psi}$, a set assumed to be already defined by structural induction on $\phi$.

$\mathcal{F}^{\psi \rightarrow \chi}_{\alpha}$ contains the tuples $(p, \tau_1, \ldots, \tau_k)$ such that it is necessarily possible that $\psi(\tau_1, \ldots, \tau_k)$ implies $\chi(\tau_1, \ldots, \tau_k)$, i.e., such that for every $q\geq p$ that forces $\psi$ there is $r\geq q$ that forces $\chi$:  for every $(q, \tau_1, \ldots, \tau_k) \in \mathcal{F}^{\psi}_{\alpha}$ with $q\geq p$ there is $(r,\tau_1, \ldots, \tau_k) \in \mathcal{F}^{\chi}_{\alpha}$ with $r\geq q$, where those latter sets are assumed defined by structural induction on $\phi$.

Finally, $\mathcal{F}_{\alpha}^{\forall x \psi}$ contains the tuples $(p, \tau_1, \ldots, \tau_k)$ such that $\forall \tau \in \mathrm{M}^{\mathbb{P}}$, it is necessarily possible that $\psi(\tau, \tau_1, \ldots, \tau_k)$ is forced: $\forall \tau \: \forall \beta$ such that $N(\beta, \tau)$ and $\forall q\geq p$ there exists $r\geq q$ that forces $\psi(\tau, \tau_1, \ldots, \tau_k)$, i.e., $(r, \tau, \tau_1, \ldots, \tau_k) \in \mathcal{F}^{\psi}_{\max(\alpha,\beta)}$. Here notice that $N(\beta, \tau)$ states that $\tau \in N_{\beta}$ is a $\mathbb{P}$-name of rank $\leq \beta$. The reader may ask: ``wait, is there no constraint on $\beta$, and can it be $>\alpha$?'' Yes: by structural inductive assumption, $\mathcal{F}^{\psi}_{\max(\alpha,\beta)}$ is defined and is a set in $\mathrm{M}$. Then, by Comprehension$^\mathrm{M}$ we can form the set of all $(p, \tau_1, \ldots, \tau_k) \in \mathbb{P}\times N_{\alpha}^k \in \mathrm{M}$ such that $\forall \tau \: \forall \beta \: (N(\beta, \tau) \rightarrow \forall q\geq p \: \exists r \geq q \: (r, \tau, \tau_1, \ldots, \tau_k) \in \mathcal{F}_{\max(\alpha, \beta)}^{\psi})$.

\begin{quote}
\textbf{The forcing chain that decides a sentence in $\mathrm{M}[G]$.} Once the sets $\mathcal{F}_{\alpha}^{\phi}$ are defined, proof of the FTF can proceed. At this point, it is instructive to ask how a sentence $\phi^{\mathrm{M}[G]}$ is decided as true or false. There is a forcing chain that goes down the structural components of $\phi$ with ``necessarily possibly forces'' conditions, and then down the name ranks in atomic formulas $=, \in$ again with ``necessarily possibly $=, \in$'' conditions, respectively. The chain is constructed straight from the definitions of the $\mathcal{F}$ sets. 

Let's look at a simple example: $\phi(y) := y = 0 \rightarrow \forall x \: \lnot (x \in y)$. At the bottom we have $\mathrm{Force}^{x \in y}(\alpha, p, x, y)$ and $\mathrm{Force}^{y=0}(\alpha, p, y, \check{0})$ defined for all $\alpha$ and $p$ as above. Then, $\mathrm{Force}^{\lnot (x \in y)}(\alpha, p, x, y) := \forall q \geq p \: \lnot \mathrm{Force}^{x\in y}(\alpha, q, x, y)$. Continuing up, $\mathrm{Force}^{\forall x \lnot (x \in y)}(\alpha, p, y) :=$ $\forall x \: \forall \beta \: (N(\beta, x) \rightarrow \forall q \geq p \: \exists r \geq q$ $\mathrm{Force}^{\lnot (x \in y)}(\max(\alpha, \beta), r, x, y))$.\footnote{Note that $\max(\alpha, \beta)$ is shorthand for asserting the existence of some ordinal $\gamma$ with the requisite property, and that $N(\beta, x)$ "type checks" $\beta$ and $x$ so there is no further need to assert for instance $\mathrm{Ord}(\beta)$.} Finally, $\mathrm{Force}^{\phi}(\alpha, p, y) :=$ $\forall q \geq p \: (\mathrm{Force}^{y=0}(\alpha, q, y, \check{0}) \rightarrow \exists r \geq q \: \mathrm{Force}^{\forall x \: \lnot (x \in y)}(\max(\alpha, \beta), r, y))$. Unrolling this, we get:
$$\mathrm{Force}^{\phi}(\alpha, p, y) := \forall q \geq p \: (\mathrm{Force}^{y=0}(\alpha, q, y, \check{0}) \rightarrow$$ $$\exists r \geq q \: \forall x \: \forall \beta \: (N(\beta, x) \rightarrow \forall s \geq r \: \exists t \geq s \: \forall u \geq t \: \lnot \mathrm{Force}^{x\in y}(\max(\alpha, \beta), u, x, y))
$$\end{quote}

\textbf{Proof of FTF for $\in$.} It is instructive to go over the proof for $\in$ and note the similarity with the proof for equality. Following the FTF proof layout, first we define the net $Z^{\in}$. Let $Z^{\in} := \{q \in \mathbb{P} \:|\:$ there exists $(\sigma, r)\in \tau_2$ with $r\geq q$ and $r\Vdash \tau_1 = \sigma \}$. By assuming FTF.Definability for equality it follows that $Z^{\in} \in \mathrm{M}$.

$Z^{\in}$ is the set of all $q$ from whose ``perspective'' it is possible that $\tau_1 \in \tau_2$. On the other hand, $\mathcal{F}^{x \in y}_{\alpha}$ contains all the triplets $(p, \tau_1, \tau_2)$ where from $p$'s ``perspective'' it is necessarily possible that $\tau_1 \in \tau_2$. Therefore, for every $(p, \tau_1, \tau_2) \in \mathcal{F}^{x \in y}_{\alpha}$ all extensions $q\geq p$ are caught in the net $Z^{\in}$. 

Next task is to show that given $\tau_1^G \in \tau_2^G$ we have $G\subset Z^{\in}$ therefore there is a $p$ all of whose extensions are in $Z^{\in}$ therefore $p$ forces $\tau_1 \in \tau_2$, and conversely.

Let $\tau_1^G \in \tau_2^G$. Let $p\in G$.  By definition $\tau_2^G = \{\sigma^G \:| \: (\sigma, q) \in \tau_2$ and $q \in G\}$, therefore there is a $q\in G$ such that $(\sigma, q) \in \tau_2$ for some $\sigma^G = \tau_1^G$.  By FTF.Truth for equality, there is some $r\in G$ with $r\Vdash \sigma = \tau$ and we can find $s \geq p, q, r$, $s \in G$. Then, $(\sigma, s) \in \tau_2$, and $s\Vdash \sigma = \tau$ follows from FTF.Coherence. Therefore $p \in Z^{\in}$. Since $p$ was arbitrary, $G \subset Z^{\in}$. Therefore, there is some $p\in G$ such that every extension $q\geq p$ is in $Z^{\in}$.

Let $p\in G$ such that every extension $q\geq p$ is in $Z^{\in}$. To show that $\tau_1^G \in \tau_2^G$.  First, note that $\tau_2^G$ is not empty: if it were, let $q$ force $\tau_2 = \check{\emptyset}$, by assuming FTF.Truth for $=$. Then take $q' \geq q, p$. Because $q' \in Z^{\in}$, there is an extension $r\geq q'$ with $r \Vdash \tau_1 = \sigma$ for $(\sigma, r) \in \tau_2$ so the set of $r$ with this property is dense above $q$ and therefore is a set in $\mathrm{M}$ by FTF.Definability for $=$. Therefore some $r'\in G$ is in that set, contradicting $\tau_2^G = \emptyset$. Similarly, letting $\tau_2^G$ be nonempty, every element of $\tau_2^G$ by definition has the form $\sigma^G$ for some $(\sigma, q) \in \tau_2$ with $q\in G$ so let $q$ be such. Then take $q'\geq q, p$, and again there is $r\geq q'$ with $r \Vdash \tau_1 = \sigma'$ for some $(\sigma', r) \in \tau_2$. Therefore the set of $r$ with that property is dense above $q$, it is in $\mathrm{M}$ by FTF.Definability for $=$, and therefore some $r'\in G$ is in that set. Therefore $\tau_1^G \in \tau_2^G$.

We just proved FTF.Truth for $\in$, namely $\tau_1^G \in \tau_2^G$ iff there is $p\in G$ such that $p \Vdash \tau_1 \in \tau_2$. 

For FTF.Definability, we need to show that $(p, \tau_1, \tau_2) \in \mathcal{F}_{\alpha}^{\in}$ iff $p\Vdash \tau_1 \in \tau_2$. If $(p, \tau_1, \tau_2) \in \mathcal{F}_{\alpha}^{\in}$, then every extension $q\geq p$ belongs to $Z^{\in}$ and we just showed that this implies $p \Vdash \tau_1 \in \tau_2$. Conversely, if $p \Vdash \tau_1 \in \tau_2$, let $q\geq p$. We must show that $q$ satisfies the condition ``there is $r\geq q$ and $(\sigma, r) \in \tau_2$ such that $(r, \tau_1, \sigma) \in \mathcal{F}_{\max(\mathrm{nr}(\tau_1), \mathrm{nr}(\sigma))}^{x=y}$''. Let $G$ be a generic ideal containing $q$, therefore $p\in G$ so $\tau_1^G \in \tau_2^G$. Therefore $q \in Z^{\in}$ as we just showed, so there is such an $r\geq q$ with  $(r, \tau_1, \sigma) \in \mathcal{F}_{\max(\mathrm{nr}(\tau_1), \mathrm{nr}(\sigma))}^{x=y}$. By FTF.Truth on $=$, that implies $r \Vdash \tau_1 = \sigma$. Therefore $q$ satisfies the condition for $p \in \mathcal{F}_{\alpha}^{\in}$. $\square$

The rest of the cases are easier. The nets $Z$ for each are given below:
\begin{itemize}

\item $Z^{\lnot} := \{q\in \mathbb{P} \:|\: q \not \Vdash \psi(\tau_1, \ldots, \tau_k) \}$ where $\phi := \lnot \psi$. 
\item $Z^{\rightarrow} := \{q \in \mathbb{P}\:|\:$ either $q \not \Vdash \psi(\tau_1, \ldots, \tau_k)$ or $\exists r \geq q$ such that $r \Vdash \chi(\tau_1, \ldots, \tau_k)\}$, where $\phi := \psi \rightarrow \chi$.
\item $Z^{\forall} := \{q \in \mathbb{P} \:|\:$ for every $\mathbb{P}$-name $\tau$ there is $r \geq q$ such that $r \Vdash \psi(\tau, \tau_1, \ldots, \tau_k)\}$, where $\phi := \forall x \: \psi$.
\end{itemize}

The proofs use the same layout and are left as exercises. $\square$

\subsubsection{Proof that $\mathrm{M}[G]$ is a model of ZFC}

Next comes a case analysis where every axiom of ZFC is proven for $\mathrm{M}[G]$. The  axioms that leverage FTF are Power, Comprehension and Replacement. The other axioms are simpler to prove. Most cases involves concocting the right $\mathbb{P}$-name that proves existence of the set postulated by the axiom.

\textbf{Extensionality.} Every transitive set is extensional and therefore $\mathrm{M}[G]$ is. $\square$

\textbf{Pairing.} Let $\sigma^G, \tau^G \in \mathrm{M}[G]$. To prove $\{\sigma^G, \tau^G \}\in \mathrm{M}[G]$ we need to devise a $\mathbb{P}$-name that evaluates to that pair. Well, let $\upsilon := \{(\sigma, p) \:|\: p \in \mathbb{P}\} \cup \{(\tau, p) \:|\: p \in \mathbb{P}\}$. The reader can verify that $\upsilon^G = \{\sigma^G, \tau^G \}$ regardless of the choice of $G$. $\square$

\textbf{Union.} Let $\tau^G \in \mathrm{M}[G]$, to prove $\cup \tau^G \in \mathrm{M}[G]$. We want a $\mathrm{P}$-name for members of members of $\tau$. Let $\upsilon := \{(\rho, p)\:|\: (\rho, p) \in \sigma$ for some $(\sigma, p) \in \tau \}$. The reader can verify with some routine bookkeeping that $\upsilon^G$ is precisely the set of elements of elements of $\tau^G$. $\square$

\textbf{Power.} Any $\tau^G \in \mathrm{M}[G]$ has subsets $s \subset \tau^G$ some of which may be in $\mathrm{M}[G]$ and some not---after all, $\mathrm{M}[G]$ is countable and therefore ``most'' subsets of any infinite $\tau^G$ are missing. The Power axiom relativized to $\mathrm{M}[G]$ says that the collection of subsets of $\tau^G$ that are in $\mathrm{M}[G]$ form a set in $\mathrm{M}[G]$: $\mathcal{P}(\tau^G)\cap \mathrm{M}[G] \in \mathrm{M}[G]$. We need to find a $\mathrm{P}$-name evaluating to precisely these subsets.

Starting with $\tau$, let $\pi = \{(\sigma, p) \:|\: \sigma \subseteq \tau$ is a $\mathrm{P}$-name and $p \in \mathbb{P}\}$.\footnote{More precisely, $\pi := \{\sigma \in \mathcal{P}(\tau)^{\mathrm{M}}\} \times \mathbb{P} \in \mathrm{M}$. It easily follows that $\pi$ is a $\mathbb{P}$-name.} Intuitively, $\pi$ contains as many subsets of $\tau$ as possible so it better work.

Clearly $\pi^G \subseteq \mathcal{P}(\tau^G)\cap \mathrm{M}[G]$. To show the reverse, the FTF will be useful. Let $\sigma^G \subseteq \tau^G$ be any subset of $\tau^G$ in $\mathrm{M}[G]$. To show that $\sigma^G \in \pi^G$ it suffices to find a $\mathbb{P}$-name $\sigma' \subseteq \tau$ such that $\sigma'^G = \sigma^G$. The following $\sigma'$ is a set in $\mathrm{M}$ thanks to FTF:Definability: $\sigma' :=  \{(\rho, p)$ such that $p \Vdash \rho \in \sigma\}$. 

To show $\sigma'^G = \sigma^G$. Let $\rho^G \in \sigma'^G$, then there is $p \in G$ such that $(\rho, p) \in \sigma'$ and $p \Vdash \rho \in \sigma$ so $\rho^G \in \sigma^G$. So $\sigma' \subseteq \sigma$. Conversely, let $\rho^G \in \sigma^G$. Using FTF:Truth, there is $p\in G$ that forces this to be true, i.e., $p\Vdash \rho \in \sigma$, therefore by definition of $\sigma'$ $(\rho, p) \in \sigma'$ which implies $\rho^G \in \sigma'^G$. $\square$

\textbf{Infinity.} It follows from $\omega \in \mathrm{M}\subset \mathrm{M}[G]$. $\square$

\textbf{Foundation.} It holds in every set. $\square$

\textbf{Comprehension Schema.} This axiom is the best example of using FTF. For any formula $\phi(x)$ relativized to $\mathrm{M}[G]$ and any set $\tau^G \in \mathrm{M}[G]$, the set $z := \{\sigma^G \in \tau^G \:|\: \phi(\sigma^G)\}$must be shown to belong to $\mathrm{M}[G]$.\footnote{The formula $\phi$ can have additional variables, $\phi(x, x_1, \ldots, x_k)$ and for any $\tau_1^G, \ldots, \tau_k^G$, the set $\{\sigma^G\in \tau^G \:|\: \phi(\sigma^G, \sigma_1^G, \ldots, \sigma_k^G)\}$ should be in $\mathrm{M}[G]$. To avoid the clutter, the extra variables are suppressed.}  It suffices to devise a $\mathbb{P}$-name $\upsilon$ evaluating to $z$. With forcing in place, the obvious $\mathbb{P}$-name works:
$$ \upsilon := \{(\sigma, p)\in \tau \:|\: p \Vdash \phi(\sigma)\} $$Importantly, $\upsilon$ is a $\mathrm{P}$-name due to FTF:Definability, specifically $\{(p, \sigma) \in \mathcal{F}_{\mathrm{nr}(\tau)}^{\phi} \:|\: (\sigma, p) \in \tau\}$ is a set in $\mathrm{M}$ due to Comprehension$^{\mathrm{M}}$ since we already know that $\mathcal{F}_{\mathrm{nr}(\tau)}^{\phi}$ is a set by Definability and that $\mathrm{nr}(\sigma) < \mathrm{nr}(\tau)$.\footnote{More explicitly, $\{x \in \mathrm{P} \times \mathrm{M}^{\mathrm{P}} \:|\: \exists \alpha\: \exists p \: \exists \sigma \: \alpha = \mathrm{nr}(\tau) \land x = (p, \sigma) \land (\sigma, p) \in \tau \land \mathrm{Force}^{\phi}(\alpha, p, \sigma)\}$.}

To show $\upsilon^G = z$ we show $\upsilon^G \subseteq z$ and $z \subseteq \upsilon^G$. Let $\sigma^G \in \upsilon^G$, therefore $(\sigma, p) \in \tau$ and $p \Vdash \phi(\sigma)$, therefore $\phi(\sigma^G)$, which implies that $\sigma^G \in z$.

Conversely, let $\sigma^G \in z$ therefore $\sigma^G \in \tau^G$ and $\phi(\sigma^G)$. Then, $(\sigma, p) \in \tau$ for some $p \in G$. Moreover, by FTF:Truth there is some $q$ that forces $\phi(\sigma^G)$. Let $r\in G$ be a common extension of $p, q$. Then $(\sigma, r) \in \tau$ and by FTF:Coherence $\tau \Vdash \phi(\sigma^G)$ so by definition of $\upsilon$, $(\sigma, r) \in \upsilon$  and therefore $\sigma^G \in \upsilon^G$. $\square$

\textbf{Replacement Schema.} Given any formula $\phi(x,y, z)$ relativized to ${\mathrm{M}[G]}$ and set $\pi^G \in \mathrm{M}[G]$ (the \textit{domain}) such that for any $\sigma^G \in \pi^G$ there is a unique $\tau^G$ such that $\phi(\sigma^G, \tau^G, \pi^G)$, we need to show that there is a set $\rho^G \in \mathrm{M}[G]$ (the \textit{image} of the domain under the class function $\phi$) that contains all elements $\tau^G$ that are making $\phi(\sigma^G, \tau^G, \pi^G)$ true for some $\sigma^G \in \pi^G$.

Fix $\mathbb{P}$-name $\pi$ in $\mathrm{M}$. For every $(\sigma, p) \in \pi$, let $\alpha_{\sigma, p}$ be the least ordinal in $\mathrm{M}$ such that there exists $\tau$ with $(p, \sigma, \tau, \pi) \in \mathcal{F}_{\alpha_{\sigma, p}}^{\phi}$, i.e., with $p \Vdash \phi(\sigma, \tau, \pi)$, if any such $p$ exists, and let $\alpha_{\sigma, p} = \emptyset$ if not. By Replacement in $\mathrm{M}$, the set of these ordinals is in $\mathrm{M}$. To see this, consider the formula $Ord(\alpha) \land (\exists p \: \exists \tau \: (p, \sigma, \tau, \pi) \in \mathcal{F}_{\alpha}^{\phi} \lor \alpha = \emptyset) \land \forall \beta \in \alpha \: \lnot (\exists p \: \exists \tau \: (p, \sigma, \tau, \pi) \in \mathcal{F}_{\beta}^{\phi})$. For any set $\pi$ and any $(\sigma, p) \in \pi$ there is a unique $\alpha$ that makes this formula true, therefore by Replacement we can collect into a set all ordinals $\alpha$ that make the formula true for elements of $\pi$'s domain.

Therefore, the supremum of $\alpha_{\sigma, p}$ for any $(\sigma, p)$, call it $\alpha^*$ is also in $\mathrm{M}$.

Now define $V_{\alpha} := \{\sigma^G \:|\: \sigma \in N_{\alpha}\}\in \mathrm{M}[G]$ for any $\alpha \in \mathrm{M}$. If we can show that any $\tau^G$ such that $\phi(\sigma^G, \tau^G, \pi^G)$ for $\sigma^G \in \pi^G$ has to be in $V_{\alpha^*}$, then we can define the set of all of such $\tau^G$ using Comprehension$^{\mathrm{M}[G]}$, which we just proved. That is to say, we can define $\{y \in V_{\alpha^*} \:|\: \exists x \in \pi^G \: \phi(x,  y, \pi^G) \}\in \mathrm{M}[G]$ by Comprehension$^{\mathrm{M}[G]}$. Therefore, all we need to show is that $\mathrm{nr}(\tau) \leq \alpha^*$ for any $\tau$ such that $\phi(\sigma^G, \tau^G, \pi^G)$.

Let $\sigma^G \in \pi^G$ be arbitrary. We know $(\sigma, p) \in \pi$ for some $p\in G$. Also, there is some $\tau^G$ making $\phi$ true and by FTF.Truth there is some $q\in G$ that forces $\phi(\sigma, \tau, \pi)$. Taking $r\geq p, q$ as a common extension, we know $(\sigma, r) \in \pi$ as well as $r \Vdash \phi(\sigma, \tau, \pi)$. Because $\alpha^*$ is the supremum of the least ordinals $\alpha_{\sigma, p}$ where $p$ forces $\phi(\sigma, \tau', \pi)$, there must be some $\tau'$ such that $\mathrm{nr}(\tau') \leq \alpha^*$ as well as $r \Vdash \phi(\sigma, \tau', \pi)$. However, $\phi$ maps unique $\sigma$ to unique $\tau$, therefore $\tau = \tau'$ and consequently $\mathrm{nr}(\tau) \leq \alpha^*$. $\square$

\textbf{Axiom of Choice.} (\textit{Sketch}.) AC is equivalent to the well-ordering principle. Let $\tau^G \in \mathrm{M}[G]$. We need to show that $\tau^G$ can be well ordered. Starting with $\tau \subset N_{\alpha} \times \mathbb{P}$ for $\alpha = \mathrm{nr}(\tau)$, we can well order $\tau$ using AC in $\mathrm{M}$: $\tau = \{(\sigma_{\beta}, p_{\beta}) \:|\: 0 \leq \beta < \gamma\}$, whereby the function $f: \tau \to \gamma$ mapping $(\sigma_{\beta}, p_{\beta})$ to $\beta$ is in $\mathrm{M}$. We can lift it up to $\mathrm{M}[G]$ with its canonical $\mathbb{P}$-name $\check{f}$ and it will easily do the job of well ordering $\tau^G$. Details are left to the reader. $\square$

Therefore, given any axiom $\phi$ of ZFC and any separative poset $\mathbb{P}$ defined in $\mathrm{M}$, it is a ZFC+ theorem that a generic $G$ over $\mathbb{P}$ exists and $\mathrm{M}[G]$ is a set containing all elements of $\mathrm{M}$ as well as $G$ and moreover $\phi^{\mathrm{M}[G]}$. Given any finite list of axioms $\phi_1, \ldots, \phi_k$, $(\phi_1 \land \ldots \land \phi_k)^{\mathrm{M}[G]}$ can be derived in ZFC+. Those axioms can then be used as explained in the next section to prove in ZFC+ that $(\lnot \mathrm{CH})^{\mathrm{M}[G]}$. Then, proof of $\mathrm{CH}$ in ZFC would lead to contradiction by relativizing the proof within $\mathrm{M}[G]$ (adding more axioms to the $\phi_i$'s as needed). We just finished outlining the proof of the following theorem, which is a numerical statement about syntactic properties of ZFC+ and can be formalized in Peano Arithmetic:

\textbf{Metatheorem.} If $\mathbb{P}$ is a separative poset in $\mathrm{M}$ and $\phi$ is any theorem of ZFC, then there is a ZFC+ proof of "$(G$ is a generic ideal over $\mathbb{P})$ and $($the set of $\mathbb{P}$-names $\mathrm{M}^{\mathbb{P}}\in \mathrm{M})$ and $($the set $\mathrm{M}[G]$ of $G$-evaluations of elements of $\mathrm{M}^{\mathbb{P}}$ is countable and transitive$)$ and $\mathrm{M} \subsetneq \mathrm{M}[G]$ and $G \in \mathrm{M}[G]$ and $\phi^{\mathrm{M}[G]}$".

\subsection{Forcing $\lnot \mathrm{CH}$ to hold in $\mathrm{M}[G]$}

To recap, $\mathrm{M}[G]$ contains $G$ together with $\aleph_2^{\mathrm{M}}$-many new subsets of $\aleph_0$ and simultaneously the proof of $CH^{\mathrm{M}[G]}$ goes through. Are we done?

Not yet: we know that $\bigcup G: \aleph_2^{\mathrm{M}} \times \aleph_0 \to \{0,1\}$ is in $\mathrm{M}[G]$ and consequently $|(2^{\aleph_0})^{\mathrm{M}[G]}|\geq \aleph_2^{\mathrm{M}}$, but we haven't shown that $\aleph_2^{\mathrm{M}[G]} = \aleph_2^{\mathrm{M}}$. Perhaps among its many new sets $\mathrm{M}[G]$ also contains a function $\aleph_1^{\mathrm{M}[G]} \to \aleph_2^{\mathrm{M}}$. Then, no contradiction with $CH^{\mathrm{M}[G]}$ is shown.

\textbf{Cardinal collapse.} We need to show that introducing $G$ through forcing does not create new functions from a smaller infinite cardinal to a larger infinite cardinal. Such a function would collapse the larger cardinal to be equal to the smaller one. Here is a simple way to achieve that through forcing. At the ground countable model $\mathrm{M}$ define the forcing notion $\mathbb{P}'$ to consist of all finite partial bijections between $\aleph_{\alpha}^{\mathrm{M}}$ and $\aleph_{\beta}^{\mathrm{M}}$ with $\beta > \alpha$. Then $\mathrm{M}[G]$ will contain $\bigcup G$, which is a complete bijection between $\aleph_{\alpha}^{\mathrm{M}}$ and $\aleph_{\beta}^{\mathrm{M}}$. Recall that both of these are countable ordinals. $\bigcup G$ would be a complete bijection because the set $D_{\gamma} := \{p\in \mathbb{P}' \:|\: \gamma$ is in the range of $p\}$ is a dense subset of $\mathbb{P}'$ for all $\gamma \in \aleph_{\beta}^{\mathrm{M}}$ and $D_{\gamma} \in \mathrm{M}$, therefore $G$ has to intersect it, i.e., $\bigcup G$ has to include $\gamma$ in its range. Similarly for its domain. Therefore  $\aleph_{\beta}^{\mathrm{M}}$ is of the same cardinality as $\aleph_{\alpha}^{\mathrm{M}}$ in $\mathrm{M}[G]$. One may ask now, what about $\aleph_{\beta}^{\mathrm{M}[G]}$? Well, this is some countable ordinal larger than $\aleph_{\beta}^{\mathrm{M}}$ that $\mathrm{M}[G]$ ``thinks'' is the $\beta$th infinite cardinal. Could new cardinals be introduced in $\mathrm{M}[G]$? No: if a bijection already exists in $\mathrm{M}$ between two ordinals, the same function still exists in $\mathrm{M}[G]$. New cardinals can never appear, but cardinals of the ground model $\mathrm{M}$ can collapse. Therefore, to conclude the proof we need to show that forcing with $\mathbb{P}$  does not collapse cardinals.

The main idea is to examine the new functions $f^G\in \mathrm{M}[G]$ created from the ``perspective'' of $\mathrm{M}$ and observe that their range cannot be ``too big'' because each $\mathbb{P}$-name $\sigma$ of the domain has at most countably$^{\mathrm{M}}$ many potential $\mathbb{P}$-name partners $\tau$ in the range. A potential partner $\tau$ in the range is one for which there is any $q_{\tau}$ forcing $f(\sigma) = \tau$. If each $\sigma$ in the domain has at most countably$^{\mathrm{M}}$ many potential partners in the range, then the range is at most $\aleph_0$ times the domain in size from the perspective of $\mathrm{M}$, and for infinite domain that implies same cardinality from the perspective of $\mathrm{M}$. 

More specifically, any function $f^G\in \mathrm{M}[G]$ between cardinals $\kappa^{\mathrm{M}} < \lambda^{\mathrm{M}}$ has an associated $\mathbb{P}$-name $f$ and some $p\in G$ that forces ``$f$ is a function $\kappa ^{\mathrm{M}}\to \lambda^{\mathrm{M}}$''. Then, for every $(\sigma^G, \tau^G) \in f^G$ (i.e., $f(\sigma^G) = \tau^G$) there has to be some $q\in G, q\geq p$ that forces $(\sigma, \tau) \in f$. Fixing $p$ and $\sigma$, the potential distinct partners $\tau$ of $\sigma$ correspond to distinct $q_{\tau}\geq p$ that are pairwise incompatible. The key property to prove now is that there are at most countably$^{\mathrm{M}}$ many pairwise incompatible elements in $\mathbb{P}$. If so, then there are at most countably$^{\mathrm{M}}$ many potential partners $\tau$ for each $\sigma$.

\begin{quote}
\textbf{Definition.} An \textbf{antichain} over a poset $\mathbb{P}$ is a set of pairwise incompatible elements. $\mathbb{P}$ satisfies the \textbf{countable chain condition (ccc)} if there is no uncountable antichain over $\mathbb{P}$. 
\end{quote}

Because $\mathrm{M}$ is countable, $\mathbb{P}$ is trivially ccc. The key will be to prove $(\mathbb{P}$ is ccc$)^{\mathrm{M}}$. Once this is done, the following lemma establishes that cardinals in $\mathrm{M}$ remain cardinals in $\mathrm{M}[G]$. The idea is  simple: given a function $f^G$, any assignment $f^G(\sigma^G) = \tau^G$ is forced by some $q\in G$. If ($\mathbb{P}$ is ccc)$^\mathrm{M}$, then $\mathrm{M}$ ``knows'' that for any given $\sigma$, there are at most countably many distinct $\tau$ forced by distinct and pairwise incompatible $q_{\tau} \Vdash f(\sigma) = \tau$. Therefore, the cardinality of the range of $f$ is forced to be at most that of the domain of $f$ times $\aleph_0$, therefore $f$ cannot map an infinite cardinal to a larger infinite cardinal. All of these statements are relativized to $\mathrm{M}$.

\begin{quote}
\textbf{Lemma.} Let $\mathbb{P}$ be a poset in $\mathrm{M}$ such that ($\mathbb{P}$ is ccc$)^{\mathrm{M}}$, $G$ a generic ideal, and $f^G: \eta \to \theta$ be a function in $\mathrm{M}[G]$ between ordinals $\eta < \theta$. Then:
\begin{enumerate}
\item There exists a function $g \in \mathrm{M}$ from $\eta$ to $2^{\theta}$ such that $\forall \alpha < \eta \: f^G(\alpha) \in g(\alpha)$, and moreover $(g(\alpha)$ is countable$)^\mathrm{M}$.
\item Infinite ordinals that are cardinals in $\mathrm{M}$ remain cardinals in $\mathrm{M}[G]$.
\end{enumerate}
\end{quote}

\textbf{Proof.} (1) Let $f$ be the $\mathbb{P}$-name for $f^G$ and $p\in G$ force ``$f: \check{\eta} \to \check{\theta}$ is a function''. For each $\alpha \in \eta$ let $g(\alpha) = \{\beta \in \theta \:|\: \exists q \geq p$ such that	$q \Vdash f(\check{\alpha}) = \check{\beta}\}$. This is a set in $\mathrm{M}$ by Comprehension and FTF.Definability. Moreover, $f^G(\alpha) \in g(\alpha)$. Finally, if for each $\beta$ we pick one element $q_{\beta}$ that forces $f(\check{\alpha}) = \check{\beta}$, then $\{ q_{\beta} \:|\: \beta \in g(\alpha)\}$ is a set of pairwise incompatible elements and therefore at most countable$^\mathrm{M}$ by ($\mathbb{P}$ is ccc$)^{\mathrm{M}}$, and is in bijection with $g(\alpha)$ by $\beta \to q_{\beta}$. Finally, because the latter function is in $\mathrm{M}$, we have $(g(\alpha)$ is countable$)^\mathrm{M}$.

(2) Let $\theta$ be an infinite ordinal  in $\mathrm{M}$ such that $(\theta$ is not a cardinal$)^{\mathrm{M}[G]}$. Then let $\eta < \theta$ be an infinite ordinal and $f^G: \eta \to \theta$ a surjection in $\mathrm{M}[G]$. Let $g: \eta \to 2^{\theta}$ be a function in $\mathrm{M}$ as in part (1). Notice that $\theta = \bigcup_{\alpha \in \eta} g(\alpha)$\footnote{Because $f^G$ is a surjection $\eta \to \theta$ and $f^G(\alpha) \in g(\alpha)$ for all $\alpha \in \eta$.} and now we get $(|\theta| = |\bigcup_{\alpha \in \eta} g(\alpha)| \leq |\eta \cdot \aleph_0| = |\eta|)^\mathrm{M}$. Therefore $\theta$ is not a cardinal in $\mathrm{M}$ either. In conclusion, a cardinal in $\mathrm{M}$ remains a cardinal in $\mathrm{M}[G]$. $\square$

The key is now to prove that $(P$ is ccc$)^{\mathrm{M}}$. Equivalently, that any uncountable$^{\mathrm{M}}$ subset of $\mathbb{P}$ contains some pair of compatible elements $p, q$. Recalling that all $p\in \mathbb{P}$ are finite partial functions, this will boil down to a basic result in infinite combinatorics , the $\Delta$-system lemma.

\begin{quote}
\textbf{$\Delta$-system lemma.} Let $X$ be an uncountable family of finite sets. Then there is an uncountable $Y \subseteq X$ and a finite set $r$, the \textit{root}, such that every pair of sets in $Y$ intersect in precisely $r$: $\forall x, y \in Y \: x\cap y = r$.
\end{quote}

\textbf{Proof sketch.} This is a basic result with many excellent expositions of the proof available online. Briefly, let $X_n$ be the sets of size $n$ among sets in $X$. At least one of these, say $X_k$ has to be uncountable otherwise $X$ would be a countable union of countable sets.  Let $m$ be the maximum size of a set $r$ such that $r\subseteq x$ for uncountably many $x \in X_k$.  $\emptyset \subset x$ for every $x\in X_k$, therefore $0 \leq m \leq k$. So let $r$ be a set of size $m$ and $Y \subseteq X_k$ the uncountable collection of all $x\supseteq r$. Any element $a$ of any $x \in Y$ is either $\in r$ or in at most countably many sets in $Y$.

Form a graph whose vertices are the sets $x \in Y$ and edges connect any sets $x, y$ whose intersection contains elements not in $r$. Any vertex in the graph has at most a finite number of elements not in $r$, and consequently at most countably many neighbors. So the graph has uncountably many vertices, each of which has at most countably many neighbors. We show that there is an uncountable \textit{independent} set of vertices, where no two vertices are neighbors. Such a set can be constructed with transfinite induction. At stage $0$ pick any vertex $x_0$. at stage $\alpha \leq \aleph_1$, we have selected a countable family of vertices $\{x_{\beta} \:|\: \beta<\alpha\}$ each of which has at most countably many neighbors outside this collection, so there are still uncountably many vertices from which to pick some $x_{\alpha}$ that is not a neighbor of any $x_{\beta}, \beta< \alpha$.

Therefore, pick an uncountable independent subset of $Y$. This is an uncountable collection in which any $x\cap y = r$. $\square$

This lemma will prove useful in showing that $(P$ is ccc$)^{\mathrm{M}}$. Given an uncountable$^{\mathrm{M}}$ collection of finite partial functions in $\mathbb{P}$, we will find an uncountable$^{\mathrm{M}}$  subcollection whose domains pairwise intersect precisely in some finite root domain $r$. Then, because each finite domain can only generate finitely many distinct functions whose range is $\{0,1\}$, there has to be some functions that are identical within $r$ and whose remaining distinct domains can be patched together. That means they are compatible, showing that there is no uncountable$^{\mathrm{M}}$ incompatible collection of partial functions in $\mathbb{P}$. More formally:

\textbf{Lemma.} $(P$ is ccc$)^{\mathrm{M}}$.

\textbf{Proof.} Let $S$ be an uncountable$^{\mathrm{M}}$ collection of $p\in \mathbb{P}$ and let $X = \{\mathrm{Dom}(p) \: |\: p \in S\}$. Any specific $\mathrm{Dom}(p)$ is a finite set, which can support at most finitely many distinct functions to $\{0,1\}$. Therefore $X$ is uncountable$^{\mathrm{M}}$. Now prove the $\Delta$-system lemma relativized to $\mathrm{M}$ and find an uncountable$^{\mathrm{M}}$ $Y \subseteq X$ and a finite root $r$ such that $x\cap y = r$ for every $x, y \in Y$. For every $x\in Y$ pick some $p_x\in S$ such that $x = \mathrm{Dom}(p_x)$. The domains of functions $\{p_x \:|\: x\in Y\}$ all intersect in precisely $r$. For any two functions to be incompatible, they should differ in how they map $r$ to $\{0,1\}$. Because there are finitely many ways to map $r$ to $\{0,1\}$, there are distinct $p_x, p_y$ such that $p_x|_r = p_y|_r$, i.e., which are compatible. $\square$

\textbf{Theorem.} $(\lnot \mathrm{CH})^{\mathrm{M}[G]}$ is a ZFC+ theorem.

\textbf{Proof.} In ZFC+, because $(P$ is ccc$)^{\mathrm{M}}$ we get $\aleph_2^{\mathrm{M}[G]} = \aleph_2^{\mathrm{M}}$. Also, $(\bigcup G: \aleph_2 \times \aleph_0 \to \{0,1\}$ is a function$)^{\mathrm{M}[G]}$. Because for each $(\alpha, n) \in \aleph_2^{\mathrm{M}} \times \aleph_0$ the set of $p \in \mathbb{P}$ containing $(\alpha, n)$ in their domain is a dense set in $\mathrm{M}$, $G$ must intersect that set and therefore all of $\aleph_2^{\mathrm{M}} \times \aleph_0$ is in the domain of $\bigcup G$, i.e., $(\bigcup G$ is a complete function$)^{\mathrm{M}[G]}$. Moreover, for any distinct $\alpha, \beta \in \aleph_2^{\mathrm{M}}$ the set of all $p \in \mathbb{P}$ that contain $(\alpha, n), (\beta, n)$ mapping to distinct values for some $n$ is a dense set in $\mathrm{M}$ and therefore $G$ must intersect that set. Therefore each ``row'' of $\bigcup G$ is distinct. By a similar argument each row of $\bigcup G$ is a subset of $\aleph_0$ not in $\mathrm{M}$. In conclusion, (there are $\aleph_2$ distinct subsets of $\aleph_0)^{\mathrm{M}[G]}$ is a ZFC+ theorem. $\square$

Therefore, if ZFC is consistent, then so is ZFC + $\lnot \mathrm{CH}$. 

\subsubsection{What values can $2^{\aleph_0}$ take?}

The above proof shows that $|(2^{\aleph_0})^{\mathrm{M}}| \geq \aleph_2^{\mathrm{M}}$ but no exact value for the cardinality of $2^{\aleph_0}$ is given. Forcing $|2^{\aleph_0}| =\aleph_2$ is possible but will not be covered here. However, notice that the value $\aleph_2$ plays no role in the proof. Can the continuum be forced to take other values?

In fact, it can take almost any value!  It can be $\aleph_3, \aleph_{2022}, \aleph_{\aleph_2}, \aleph_{\aleph_{\aleph_{2022}}}, \ldots$. There is only one restriction that comes from a well known and useful result, K\"{o}nig's theorem:

\begin{quote}
\textbf{K\"{o}nig's theorem.} Let $\kappa_i, \lambda_i$, $i \in I$ be two infinite sequences of cardinals with $\kappa_i \leq \lambda_i$, $\forall i \in I$, where $I$ is an index set of any cardinality. Then $\Sigma \kappa_i < \Pi \lambda_i$. 
\end{quote}

Note the strict inequality in the above theorem, as well as the fact that the $\kappa_i$ and $\lambda_i$ can be finite or infinite. Proof of this theorem requires the Axiom of Choice.

A constraint on the value of the continuum follows. To describe it, the notion of \textit{cofinality} needs to be introduced.

\begin{quote}
\textbf{Definition: cofinality.} Let $\kappa$ be any ordinal. A sequence of ordinals $\langle \lambda_{\xi} < \kappa \:|\: \xi < \beta\rangle$ is \textbf{cofinal} to $\kappa$ if for every $\alpha < \kappa$ there is a $\xi < \beta$ such that $\lambda_{\xi} > \alpha$. The cofinality of $\kappa$, denoted $\mathrm{cf}(\kappa)$, is the least $\beta$ for which such a sequence exists.

\textbf{Examples.} The cofinality of the ordinal $2022$ is $1$: the sequence $\langle 2021 \rangle$ suffices. The cofinality of $\omega$ is $\omega$: no finite sequence of finite numbers has $\omega$ as its limit. A cardinal $\kappa$ is called \textbf{regular} when it has this precise property of $\mathrm{cf}(\kappa) = \kappa$; otherwise it is called \textbf{singular}. The cofinality of $\aleph_{\aleph_0}$ is $\omega$:  the sequence $\langle \aleph_i \:|\: i < \aleph_0 \rangle$ is cofinal; therefore $\aleph_{\aleph_0}$ is singular. 
\end{quote}

An important corollary of K\"{o}nig's theorem involves the cofinality of an infinite cardinal:

\begin{quote}
\textbf{Corollary.} If $\lambda \geq \aleph_0$ then $\lambda < \lambda^{cf(\lambda)}$.

\textbf{Proof.} Let $\langle \kappa_{\xi} < \lambda \:|\: \xi < cf(\lambda) \rangle$ be a cofinal sequence to $\lambda$. Use K\"{o}nig's theorem with $\lambda_i = \lambda$ for $i \in \mathrm{cf}(\lambda)$ to obtain $\Sigma_{\xi < cf(\kappa)} \kappa_{\xi} < \lambda^{\mathrm{cf}(\lambda)}$. Then, observe that $\Sigma_{\xi < \mathrm{cf}(\kappa)} \kappa_{\xi} = \lambda$ so the result follows. $\square$
\end{quote}

It follows easily that the cofinality of continuum cannot be $\omega$: notice that $(2^{\omega})^{\omega} = 2^{2\omega} = 2^{\omega}$. Assume now that $\mathrm{cf}(2^{\omega}) = \omega$. Let $\langle \kappa_i < 2^{\omega} \:|\: i < \omega \rangle$ be cofinal to $2^{\omega}$. Then $\Sigma \kappa_i < (2^{\omega})^{\omega}$ by K\"{o}nig's theorem. However $\Sigma \kappa_i = 2^{\omega}$, contradiction.

Therefore, the cardinality of continuum cannot be $\aleph_{\omega}$ or $\aleph_{\omega^2}$ or $\aleph_{\aleph_{\aleph_{\aleph_{\ldots}}}}$ because all these have cofinality $\omega$. However, it can be bigger or smaller than any of these values.

\section{Bullet-point outline of the proof}

\subsection{Part I: logical argument}

\begin{enumerate}
\item Any sentence $\phi$ is equivalent to $\phi^M$ where $M$ is a countable transitive set. In particular, given any finite list of axioms of ZFC $\phi_1, \ldots, \phi_k$, $\exists M \: \mathrm{CntTrans}(M) \land \phi_1^M \land \cdots \land \phi_k^M$ is a ZFC theorem. This is established by the Reflection Principle and Mostowski collapse lemma.
\item Therefore, for any theorem $\phi$ of ZFC, the proof can be copied to a proof $\exists M \: \mathrm{CntTrans}(M) \land \phi^M$.
\item Extend ZFC to ZFC+ by adding the new constant symbol $\mathrm{M}$, the axiom $\mathrm{CntTrans}(\mathrm{M})$, and the axiom $\phi^{\mathrm{M}}$ for each axiom $\phi$ of ZFC . If ZFC is consistent, so is ZFC+.
\item In ZFC+, $\mathrm{M}$ is a symbol that represents a countable transitive set \textit{by fiat} and $\aleph_0^{\mathrm{M}} < \aleph_1^{\mathrm{M}} < \aleph_2^{\mathrm{M}} < \ldots$ are countable ordinals that are uncountable cardinals in $\mathrm{M}$.
\item Assume there is a ZFC proof of $\mathrm{CH}$ for contradiction. The proof is a finite list of axioms $\phi_1, ..., \phi_k$ and their logical consequences $\phi_{k+1}, ..., \phi_n = \mathrm{CH}$. The same proof will work in ZFC+ by relativizing each $\phi_i$ to $\mathrm{M}$ to show $CH^{\mathrm{M}}$, that is, $(2^{\aleph_0})^{\mathrm{M}} = \aleph_1^{\mathrm{M}}$. 
\item Extend $\mathrm{M}$ to a bigger set $\mathrm{M}[G]$ with a table $G$ that encodes $\aleph_2^{\mathrm{M}}$-many new subsets of $\aleph_0$, so that (a) for any axiom $\phi$ there is a ZFC+ proof that $\phi^{\mathrm{M}[G]}$; (b) there is a ZFC+ proof that $\lnot \mathrm{CH}^{\mathrm{M}[G]}$. This is accomplished in Part II, Forcing.
\item By (6a), the proof of $CH^\mathrm{M}$ can be syntactically copied to a proof $CH^{\mathrm{M}[G]}$ together with the ZFC+ proofs of $\phi^{\mathrm{M}[G]}$ for each axiom $\phi$ needed. Those axioms include $\phi_1, \ldots, \phi_i$ as well as axioms required to complete the forcing construction in (6).
\item By 7 and (6b), ZFC+ is inconsistent, therefore by (4) ZFC is inconsistent. Assuming ZFC is consistent, the conclusion is that there is no proof of $\mathrm{CH}$. 
\end{enumerate}

How to extend $M$ to $\mathrm{M}[G]$ is Part II, Forcing. This part is done inside ZFC+.

\subsection{Part II: forcing}

Every statement below, unless otherwise indicated, is preceded with ``it is derivable in ZFC+ that...'' In particular, $\mathrm{M}$ is postulated by fiat and its existence is impossible to derive in ZFC due to G\"{o}del's incompleteness theorem. 

\begin{enumerate}
\item Since $\mathrm{M}$ is countable, it contains countably many subsets of $\aleph_0^{\mathrm{M}} = \aleph_0$, and $\aleph_2^{\mathrm{M}}$ is countable. Therefore, there exist $\aleph_2^{\mathrm{M}}$-many new  subsets of $\aleph_0$ not in $\mathrm{M}$. Represent this collection by an object $G$ such that $\bigcup G: \aleph_2^{\mathrm{M}}\times \aleph_0 \to \{0,1\}$ is a complete function whose ``rows'' are functions $\aleph_0 \to \{0,1\}$ representing pairwise distinct subsets of $\aleph_0$ that are not in $\mathrm{M}$. $G$ will be added to $\mathrm{M}$ to produce a new set $\mathrm{M}[G]$ that is a model of ZFC+ with $\lnot \mathrm{CH}^{\mathrm{M}[G]}$.
\item First, define the poset $\mathbb{P}$ consisting of all finite partial functions $p: \aleph_2^{\mathrm{M}}\times \aleph_0 \to \{0,1\}$ ordered by inclusion; each $p$ is called a forcing condition. By Pairing, each $p \in \mathbb{P}$ is in $\mathrm{M}$ and $\mathbb{P} \in \mathrm{M}$. Each $p, q$ are compatible if $\exists r \: p \leq r \land q \leq r$ and can be patched together to a bigger partial function.
\item Any complete function $\bigcup G: \aleph_2^{\mathrm{M}}\times \aleph_0 \to \{0,1\}$ is the union of a maximal \textbf{ideal} of $\mathbb{P}$: a maximal downward closed subset of $\mathbb{P}$ in which members are pairwise compatible.
\item For $\mathrm{M}[G]$ to be a model and satisfy Comprehension$^{\mathrm{M}[G]}$, it has to include the preimages under $G$ of sets $r \subset x \times \mathbb{P}$ for $x, r \in\mathrm{M}$. In addition, this operation of taking preimages under $G$ can be performed transfinitely. To do so, define the $\mathbb{P}$-name hierarchy $\mathrm{M}^{\mathbb{P}} \subset \mathrm{M}$ (a ZFC+ set that is a proper $\mathrm{M}$-class), by $N_{\alpha+1} := \mathcal{P}(N_{\alpha} \times \mathbb{P}) \cap \mathrm{M}$ for any successor ordinal $\alpha$, and $N_{\alpha} := \bigcup_{\beta<\alpha} N_{\beta}$ for limit ordinal $\alpha$.
\item $\mathrm{M}^{\mathbb{P}} = \bigcup N_{\alpha}$. Define the $G$-evaluation of a $\mathbb{P}$-name $\tau$ to be $\tau^G := \{\sigma^G \:|\: (\sigma, p) \in \tau \land p \in G\}$. If $\mathrm{M}[G]$ is a model satisfying Comprehension, it has to include each such $\tau^G$. Now by wishful thinking let $\mathrm{M}[G] := \{\tau^G \:|\: \tau \in \mathrm{M}^{\mathbb{P}}\}$, and aim to prove that $\mathrm{M}[G]$ thus minimally defined is indeed a model.
\item Define the canonical name of $x\in \mathrm{M}$ to be $\check{x} := \{(x,p) \:|\: p\in \mathbb{P}\}$ and $\Gamma := \{(\check{p}, p) \:|\: p \in \mathbb{P}\}$, show $\check{x}^G = x$ and $\Gamma^G = G$, concluding that $\mathrm{M}\subsetneq \mathrm{M}[G]$ and $G \in \mathrm{M}[G]$.
\item For $\mathrm{M}[G]$ to be a model, it has to satisfy all ZFC axioms. To prove Comprehension$^{\mathrm{M}[G]}$, Comprehension$^{\mathrm{M}}$ needs to be used. Specifically, $z := \{\sigma^G \in \tau^G \:|\: \phi(\sigma^G)\}$ has to be in $\mathrm{M}[G]$ for every $\sigma, \tau \in \mathrm{M}^{\mathbb{P}}$. It suffices to show that $z = \upsilon^G$ for some $\upsilon \in \mathrm{M}^{\mathbb{P}}$.
\item The natural candidate is $\upsilon := \{(\sigma, p) \in \tau \:|\: \psi_{\phi}(\sigma, p)\}$ for some formula $\psi_{\phi}$, whereby $\psi_{\phi}(\sigma, p)$ is true when $p\in G$ forces $\phi(\sigma^G)$ to be true.
\item It immediately follows that $G$ has to be generic: no new infinite "AND" information can be encoded in $G$. In particular, for every $p$ and every sentence $\phi$, if every $q \geq p$ has an extension $r \geq q$ that makes $\phi$ true in every ideal $G \ni r$, then $p$ is forced to make $\phi$ true: $\phi$ has to be decided in $\mathrm{M}[G]$, and no $q \geq p$ can possibly decide $\lnot \phi$ is true. Iterating this requirement by structural induction, $p$ forces $\phi$ to be true whenever every $q\geq p$ has $r\geq q$ forcing $\phi$ to be true.
\item Genericity of $G$ is defined by the concept of density. A subset $D \subset \mathbb{P}$ is \textbf{dense} if for all $q\in \mathbb{P}$ there is a $r\geq q$ such that $r \in D$. $G$ is \textit{generic} over $\mathrm{M}$ if $G\cap D \not = \emptyset$ for any $D\in M$ that is a dense subset of $\mathbb{P}$. It follows that if $G$ is generic, for any $p\in G$ and any sentence $\phi$, if the set $D$ of $q\in \mathbb{P}$ that force $\phi$ to be true is dense above $p$, then $p$ forces $\phi$ to be true.
\item Generic ideals $G$ exist as long as the ground model is countable: enumerate all dense $D_1, D_2, \ldots$ in $\mathrm{M}$; at each step $i+1$ pick $p_{i+1} \in D_{i+1}$ such that $p_{i+1}\geq p_i$; let $G = \{p \in \mathbb{P} \:|\: p \leq p_i$ for some $i \}$. $G$ is not in $\mathrm{M}$ if $\mathbb{P}$ is separative: for any $f \in \mathrm{M}$ of the same domain and range as $\bigcup G$, $\{p \in \mathbb{P} \: | \: p$ disagrees with $f$ in some position$\}$ is dense $\in \mathrm{M}$.
\item A condition $p$ \textbf{forces} $\phi(\tau_1, \ldots, \tau_k)$ if $p \in G$ implies $\phi^{\mathrm{M}[G]}(\tau_1^G, \ldots, \tau_k^G)$. The fundamental theorem of forcing (FTF) states that for any $\phi$, the collection of $(p, \tau_1, \ldots, \tau_k)$ where $p\Vdash \phi$ and the $\tau_i$s have max name rank $\alpha$ is a set $\mathcal{F}^{\phi}_{\alpha} \in \mathrm{M}$ (Definability), and $\phi^{\mathrm{M}[G]}$ iff there is some $p\in G$ that forces $\phi$ (Truth). Moreover, $p\Vdash \phi$ and $q \geq p$ implies $q \Vdash \phi$. FTF is a metatheorem: for any $\phi$ there is a separate ZFC+ proof of FTF for $\phi$.
\item Using the FTF, the metatheorem that $\mathrm{M}[G]$ is a model is proven. For any axiom $\phi$ of ZFC, there is a ZFC+ proof of $\phi^{\mathrm{M}[G]}$. Therefore $CH^{\mathrm{M}[G]}$ goes through.
\item It remains to show that forcing does not collapse cardinals. This is achieved by showing that for each function $f$ in $\mathrm{M[G]}$ each element $x$ in $f$'s domain is forced to map to at most countably many distinct $y$ in $f$'s range by elements $q_y$ in $\mathbb{P}$. A combinatorial argument finishes the proof.
\end{enumerate}

\section{Coming back to omitted steps}

\subsection{Reflection Principle}

Any  sentence can be reflected in a countable transitive set $M$. This is a metatheorem, or better, an algorithm that takes any $\phi$ as input and returns a ZFC proof of $\exists M \: \mathrm{CntTrans}(M)\land \phi \leftrightarrow \phi^M$. We sketch the algorithm here following the excellent exposition in Weaver (2014).

The first step is a lemma schema that takes as input a finite list of formulas $\phi_1, \ldots, \phi_n$ and produces a ZFC proof of existence of a countable set $M$ that contains all the elements needed so that the formulas have ``similar'' behavior with their variables restricted to $M$ as when their variables are unconstrained. More precisely, for any $\phi(x_1, \ldots, x_k)$ in the list with all free variables displayed and for any $a_1, \ldots, a_k \in M$, if there exists $x_j$ such that $\phi(a_1, \ldots, a_{j-1}, x_j, a_{j+1}, \ldots, a_k)$, then there exists $x_j\in M$ such that $\phi(a_1, \ldots, a_{j-1}, x_j, a_{j+1}, \ldots, a_k)$.

\textbf{Lemma Schema.} Let $\phi_1, \ldots, \phi_n$ be any formulas with free variables among $x_1, \ldots, x_k$. There is a ZFC proof that there exists a countable set $M$ such that for any formula $\phi_i$ in the list and for any variable $x_j$ and any constants $a_1, \ldots, a_k \in M$, whenever $\exists x_j \: \phi_i(a_1, \ldots, x_j, \ldots, a_k)$ holds, $\exists x_j \in M \: \phi_i(a_1, \ldots, x_j, \ldots, a_k)$ also holds.\footnote{Here, $a_j$ is a dummy constant but it is not called out as such in order to reduce notational clutter.} That is to say, whenever there exists $x_j$ that makes formula $\phi_i$ true with other variables assigned to specific elements of $M$, keeping the same assignments for other variables there exists $x_j \in M$ making $\phi_i$ true.

\textbf{Proof Sketch.} We create a countable sequence of sets $M_0 = \omega \subset M_1 \subset M_2 \subset \cdots$ whereby $M := \bigcup_l M_l$. At each step $l+1$, for every formula $\phi_i$ and variable $x_j$, and for all possible assignments of the other variables to elements from $M_l$, if there exists $x_j$ making the formula true, choose one such element and add it to $M_{l+1}$.\footnote{This requires AC, or at least that $\{x_j \:|\: \phi(\ldots x_j\ldots)\}$ can be well-ordered. With AC, we can choose an element from the set of $x_j$ of lowest rank making $\phi_i$ true, and details are omitted.}

In other words, for each $i,j$ and for each $a_1, \ldots, a_k \in M_l$ if $\exists x_j \: \phi_i(a_1, \ldots, x_j, \ldots, a_k)$ then choose one element $a$ such that $\phi_i(a_1, \ldots, x_j = a, \ldots, a_k)$ and add $a$ to $M_{l+1}$. 

The set $M_{l+1}$ is countable by induction: assuming $M_l$ is countable, there are countably many tuples $a_1, \ldots, a_k$ times finitely many formulas times finitely many variables, and each such combination generates at most one additional element to be added to $M_{l+1}$. Therefore $M$ is countable.

Now assume that for some $a_1, \ldots, a_j-1, a_j+1,\ldots, a_k \in M$, $\exists x_j \: \phi_i(a_1, \ldots, x_j, \ldots, a_k)$. In other words for some $a$ not necessarily in $M$, $\phi_i(a_1, \ldots, a_j-1, a, a_j+1, \ldots, a_k)$. Well, each one of $a_1, \ldots, a_k$ is in some $M_l$, so choose $l^*$ such that $a_1, \ldots, a_k \in M_{l^*}$. But then some $a^*$ should have been added to $M_{l^*+1}$ such that $\phi_i(a_1, \ldots, a^*, \ldots, a_k)$. Therefore some such $a^*$ is in $M$. Therefore $\exists x_j \in M \: \phi_i(a_1, \ldots, x_j, \ldots, a_k)$.

Given any input formulas $\phi_1, \ldots, \phi_n$, the above proof can be written down mechanically in full ZFC formality. $\square$

\textbf{Theorem: The Reflection Principle Algorithm.} Given any input sentence $\phi$ there is a ZFC proof that there is a countable transitive set $M$ such that $\phi \leftrightarrow \phi^M$.

\textbf{Proof.} Given $\phi$, the strategy to generate $M$ will be to break down $\phi$ into the list $\phi_1, \ldots, \phi_n$ of all its component formulas, apply the above lemma to the list, and then show by structural induction on $\phi$ that all subformulas as well as $\phi$ are reflected in the resulting set $M$. For this to work technically, the list of formulas needs to be constructed carefully. First, the axiom of extensionality will be needed (for all $x$ and $y$, $x = y \leftrightarrow (\forall z \: z\in x \leftrightarrow z \in y)$. So let's assume our starting list is $\phi_1 = \phi, \phi_2 = $ Extensionality. Now we break down to components to expand the list: if $\phi_i = \lnot \psi$ is in the list, add $\psi$ to the list; if $\phi_i = \psi \land \chi$ or $\phi_i = \psi \rightarrow \chi$ is in the list, add both $\psi$ and $\chi$ to the list; if $\phi_i = \forall x \: \psi$, add $\psi$ to the list. Continue until no more formulas can be added. Now we need another trick: if the resulting list so far is $\phi_1, \ldots, \phi_n$, apply the above lemma to the negations: $\lnot \phi_1, \ldots, \lnot \phi_n$.

Now let $M$ be the countable set resulting from the lemma. Using structural induction starting from atomic $\phi_i$ (i.e., of the form $x = y$ or $x\in y$) and moving up in formula complexity, we will show for each $\phi_i$ that for all $a_i, \ldots, a_k \in M$, $\phi_i(a_1, \ldots, a_k) \leftrightarrow \phi_i^M(a_1, \ldots, a_k)$.

Let $\phi_i$ be atomic: well, $(x_1 = x_2)^M$ is $x_1 = x_2$ and $(x_1\in x_2)^M$ is $x_1 \in x_2$ so this case is trivial.

Let $\phi_i = \lnot \psi$. By inductive assumption we have $\psi \leftrightarrow \psi^M$ for all $a_j \in M$, and therefore $\lnot \psi \leftrightarrow \lnot \psi^M$ therefore $\lnot \psi \leftrightarrow (\lnot \psi)^M$.

Similarly the cases of $\psi \rightarrow \chi$ or $\psi \land \chi$ are trivial.

The only case that is interesting is when $\psi(x, x_1, \ldots, x_k) \leftrightarrow \psi^M(x, x_1, \ldots, x_k)$ has been proven for a subformula that has the unbound variable $x$ (with other unbound variables $x_j$ displayed), and we need to show that for all $a_1, \ldots, a_k \in M$: $$\forall x \: \psi(x, a_1, \ldots, a_k) \leftrightarrow \forall x \in M \: \psi^M(x, a_1, \ldots, a_k)$$
Fix the $a_j$s in $M$. If $\forall x \: \psi(x, a_1, \ldots, a_k)$, then certainly $\forall x \in M\: \psi(x, a_1, \ldots, a_k)$. By inductive assumption, $\psi(\ldots) \leftrightarrow \psi^M(\ldots)$ therefore $\forall x \in M\: \psi^M(x, a_1, \ldots, a_k)$. Conversely, if $\forall x \: \psi(x, a_1, \ldots, a_k)$ is false, then $\exists x \: \lnot \psi(x, a_1, \ldots, a_k)$ and by the lemma some $x\in M$ should have been included (here is why we needed to apply the lemma with all \textit{negations} of the $\phi_i$s) and therefore $\exists x\in M \: \lnot \psi(x, a_1, \ldots, a_k)$, and by inductive assumption  $\exists x\in M \: \lnot \psi^M(x, a_1, \ldots, a_k)$ or equivalently $\forall x\in M \: \psi^M(x, a_1, \ldots, a_k)$ is false.

Since the sentence $\phi$ is in the list, $\phi \leftrightarrow \phi^M$ is established. Moreover, $M$ is countable and extensional: $\forall x, y \in M \: x=y \leftrightarrow (\forall z\in M \: z \in x \leftrightarrow z \in y)$. Now the Mostowski collapse theorem, discussed below, can be applied to $M$ to obtain countable transitive set $M'$ where $\phi \leftrightarrow \phi^{M'}$.

Given any input $\phi$, the above proof can be written down mechanically in full ZFC formality. $\square$

It remains to establish that an extensional set $M$ can be converted to a transitive set $N$ so that for every formula $\phi$, $\phi^M \leftrightarrow \phi^{N}$. The map $f$ from $M$ to $N$ will be an $\in$-\textbf{isomorphism}:  a bijection such that $x \in y$ iff $f(x) \in f(y)$ for all $x, y\in M$.

\textbf{Theorem.} Let $f: M \to N$ be an $\in$-isomorphism. Let $\phi(x_1, \ldots, x_k)$ be any formula with free variables displayed. Then for any $a_1, \ldots, a_k \in M$, the formula $$\phi^M(a_1, \ldots, a_k) \leftrightarrow \phi^N(f(a_1), \ldots, f(a_k))$$is provable in ZFC.

\textbf{Proof sketch.} This is an easy induction on the structure of $\phi$. The base cases $x=y$ and $x\in y$ as well as the cases of $\rightarrow$ and $\land$ are trivial. Therefore, assume we have a ZFC proof that $\phi^M(a_1, \ldots, a_k) \leftrightarrow \phi^N(f(a_1), \ldots, f(a_k))$ for all $a_1, \ldots, a_k \in M$.

We need to show $\forall x_j \: \phi^M(a_1, \ldots, x_j, \ldots, a_k) \leftrightarrow \forall x_j \: \phi^N(f(a_1), \ldots, x_j, \ldots, f(a_k))$. Because $f$ is a bijection, if $\forall x_j \: \phi^M(\ldots)$ fails, the $a_j$ that makes it fail is mapped 1-1 to an $f(a_j)$ that makes $\phi^N(\ldots f(a_j) \ldots)$ fail by inductive assumption, and therefore $\forall x_j \: \phi^N(\ldots)$ fails. Similarly if $\forall x_j \: \phi^M(\ldots)$ holds, then for any $a_j$ $\phi^M(\ldots a_j \ldots)$ holds, by inductive assumption $\phi^N(\ldots f(a_j) \ldots)$ holds and therefore by $f$ being surjective $\forall x_j \: \phi^N(\ldots)$ holds. $\square$

Finally, given a set $M$ where extensionality holds, we can find an $\in$-isomorphism $f$ to a transitive set $N$. The following theorem is provable in a much more general form, where the relationship does not need to be $\in$ as long as it is well founded. However, we won't need the general form.

\textbf{Theorem: Mostowski Collapse Lemma.} Let $M$ be a set where Extensionality$^M$ holds. $M$ is $\in$-isomorphic to a transitive set $N$.

\textbf{Proof.} Define $f: M \to V$ by $f(x) = \{f(y) \: | \: y \in x\}$. This definition is by transfinite recursion and it makes sense because of Foundation.\footnote{More generally, instead of '$\in$' a relation $R$ that is well founded, i.e., where any set of elements has an $R$-minimal element, suffices because well-foundedness enables transfinite recursion on $R$.}

The set $N:= \{y \: |\: \exists x \in M \: y = f(x) \}$ exists by Replacement. We need to show that $f$ is an $\in$-isomorphism and that $N$ is transitive.

To show $f$ is an $\in$-isomorphism, let $y\in x \in M$. Then $f(x) = \{f(y) \: |\: y \in x\}$ so by definition it follows that $f(y) \in f(x)$.

To show that $N$ is transitive, let $w \in z \in N$, to show $w \in N$. Well, $z = f(x)$ for some $x \in M$, therefore $z = \{f(y) \: |\: y \in x\}$ therefore $w = f(y)$ for some $y \in x$ therefore $w$ is in the image of $f$, i.e., $w \in N$. $\square$

Therefore, given any sentence $\phi$, a ZFC proof of $\exists M \: \mathrm{CntTrans}(M) \land \phi \leftrightarrow \phi^M$ can be written down mechanically. This is a ridiculously powerful result. Given any finite number of axioms of ZFC, they hold in a countable transitive set $M$. Equivalently, given any ZFC proof of some $\phi$, we can convert it to a ZFC proof of $\exists M\: \mathrm{CntTrans}(M) \land \phi^M$.

\textbf{Example of Mostowski Collapse.} The reader may wonder what, if anything, is the effect of the Mostowski Collapse procedure in a set. Let's give a simple example of $M$ and calculate the Mostowski collapse:
$$M := \{1, 2, 2022, \omega, \omega\cdot 2022, \{2022, \omega\cdot 2022\}\}$$By definition $f(1) = \{f(y) \: |\: y\in x\}$. Now, $\emptyset \in 1$, however because $\emptyset \not \in M$ it is not in the domain of $f$ therefore $f(1)$ becomes $\emptyset$ or if you prefer, $0$. Similarly, $f(2) = 1$ and $f(2022) = 2$.

The reader can easily verify that $f(\omega) = \{0, 1, 2\} = 3$ and similarly $f(\omega\cdot 2022) = 4$. Finally, $f(\{2022, \omega\cdot 2022 \}) = \{2, 4\}$.

Therefore the resulting $N$ is $\{0, 1, 2, 3, 4, \{2,4\}\}$, a transitive set $\in$-isomorphic to $M$ and quite a simplification from $M$.

\textbf{Existence of $M$ versus  existence of a proof of existence of $M$.} Given $\phi$, the Reflection Principle algorithm outputs a ZFC proof of $\exists M\: \mathrm{CntTrans}(M) \land \phi\leftrightarrow \phi^M$. It always terminates and can be implemented in software. The algorithm goes through all subformulas of $\phi$ and writes down the respective ZFC subproofs as sketched above. The resulting ZFC proof states the existence of an countable set $M$ with the requisite properties. This latter proof, by contrast, is profoundly non-constructive. It involves infinitely many iterative applications of the axiom of choice. There is no ontological conclusion that can be drawn about the existence of any such $M$. $M$ is neither shown to truly exist (whatever that would mean), neither is anything of this sort needed for the proof of the independence of $\mathrm{CH}$. All we need is that we can syntactically and mechanically convert any ZFC theorem $\phi$ to a ZFC theorem $\exists M\: \mathrm{CntTrans}(M) \land \phi^M$.

\subsection{Forcing consistency of $\mathrm{CH}$}

Consistency of $\mathrm{CH}$ is the second part of independence of $\mathrm{CH}$. Godel (1939) originally showed it by building an inner model of ZFC --- a formula $L(x)$ that describes a proper class of sets such that all axioms $\phi$ hold as $\phi^L$ and such that $CH^L$ holds. $L(x)$ is the class of constructible sets, a topic beyond scope here. The reader is referred to Smullyan and Fitting (2010).

The forcing machinery can be deployed to provide a much simpler proof of this part. The argument below follows Weaver (2014).

Starting from the same ground model $\mathrm{M}$ in ZFC+, the goal now is to extend through forcing to a model $\mathrm{M}[G]$ that contains a bijection between $(2^{\aleph_0})^{\mathrm{M}[G]}$ and $\aleph_1^{\mathrm{M}[G]}$ so that $(|2^{\aleph_0}| = \aleph_1)^{\mathrm{M}[G]}$. Can the reader think of what object $G$ we may want to add to $\mathrm{M}$?

The first thing that comes to mind is a bijection between the countable sets $(2^{\aleph_0})^{\mathrm{M}}$ and $\aleph_1^{\mathrm{M}}$. Let's define a separative poset $\mathbb{P}$ that contains partial constructions $p$ of this bijection such that $p\in \mathrm{M}$ for all $p\in \mathbb{P}$. Analogously to proving consistency of $\lnot \mathrm{CH}$, we may be inclined to let $p$ be finite partial bijections $(2^{\aleph_0})^{\mathrm{M}} \to \aleph_1^{\mathrm{M}}$. So let's follow this through and see what properties $\bigcup G$ and $\mathrm{M}[G]$ have: 
\begin{itemize}
\item For any subset $s \subset \aleph_0$ in $\mathrm{M}$, the set $D_{s} := \{p \in \mathbb{P} \:|\: s \in \mathrm{Dom}(p)\}$ is a dense set in $\mathrm{M}$, therefore $s \in \mathrm{Dom}(\bigcup G)$.
\item For any element $\alpha \in \aleph_1^{\mathrm{M}}$, the set $D'_{\alpha} := \{p\in \mathbb{P} \:|\: \alpha \in \mathrm{Ran}(p)\}$ is a dense set in $\mathrm{M}$, therefore $\alpha \in \mathrm{Ran}(\bigcup G)$.
\item In conclusion, $\bigcup G$ is a complete bijection and $((2^{\aleph_0})^{\mathrm{M}} = \aleph_1^{\mathrm{M}})^{\mathrm{M}[G]}$. 
\end{itemize}

Are we done? Not quite: the set $(2^{\aleph_0})^{\mathrm{M}[G]}$ could be bigger than $(2^{\aleph_0})^{\mathrm{M}}$ if forcing with $\mathbb{P}$ introduced new subsets of $\aleph_0$. Moreover, $\aleph_1^{\mathrm{M}}$ could be countable in $\mathrm{M}[G]$ if forcing with $\mathbb{P}$ introduces a new surjection $\aleph_0 \to \aleph_1^{\mathrm{M}}$. And indeed, just like in the case of forcing for consistency of $\lnot \mathrm{CH}$, forcing with finite partial functions $(2^{\aleph_0})^{\mathrm{M}} \to \aleph_1^{\mathrm{M}}$ introduces new functions with domain $\aleph_0$ and in particular new subsets of $\aleph_0$.

Here is some intuition why. Take any injective function $f\in \mathrm{M}$, $f: \: \omega \to (2^{\aleph_0})^{\mathrm{M}}$. There are many such functions that are always $\in \mathrm{M}$. For example, $f(n) = n$ is such a function because $n \in (2^{\aleph_0})^{\mathrm{M}}$ when $n$ is viewed as an ordinal. Also, $f_k(n) := kn$ is such a function, as is really any injective function $f :\omega \to 2^{\aleph_0}$ definable by a formula $\phi(x, y)$ that is true iff $x \in \omega$ and $y \in f(x)$. 

Now take $\bigcup G \circ f$ to produce a function $s_f := \omega \to \aleph_1^{\mathrm{M}}$. Can the reader see how new subsets of $\omega$ are formed?

Let's cut the range into two infinite parts within $\mathrm{M}$, and call one of them $0$ and the other $1$. For example (as in Weaver (2014)) we could say that $\alpha < \omega$ corresponds to $0$ and $\omega \leq \alpha < \aleph_1^{\mathrm{M}}$ corresponds to 1. Any partition in two infinite pieces works as long as they are in $\mathrm{M}$.

This automatically converts $s_f$ to a function $\omega \to \{0,1\}$, or equivalently a subset of $\omega$. Could this subset possibly be in $\mathrm{M}$ already? Not if $G$ is generic: given any $p \in \mathbb{P}$, and given any $\mathrm{M}$-set $t\subset \omega$, it is necessarily possible to extend $p$ to $q\geq p$ that makes $s_f$ and $t$ differ in some position. That is where we need both pieces of the partition to be infinite so that $p$ exhausts neither of them.

The problem of generating new $\omega$-subsets arises because any $p \in \mathbb{P}$ only forces a finite part of $s_f$. How can we change $\mathbb{P}$ in a way to ensure that $s_f \in \mathrm{M}$ whenever $f \in \mathrm{M}$? The solution of how to do this emerges by observing that given $\bigcup G: (2^{\aleph_0})^{\mathrm{M}}\to \aleph_1^{\mathrm{M}}$, any $f : \omega \to (2^{\aleph_0})^{\mathrm{M}}$ in $\mathrm{M}$ induces a countable partial bijection $p: (2^{\aleph_0})^{\mathrm{M}}|_{\mathrm{Ran}(f)}  \to \aleph_1^{\mathrm{M}}$, whose domain is $\mathrm{Ran}(f)$ and such that $p \subset \bigcup G$. 

So, if we ensure that $\mathbb{P}$ contains all countable partial bijections $p: (2^{\aleph_0})^{\mathrm{M}} \to \aleph_1^{\mathrm{M}}$ such that $p \in \mathrm{M}$, the hope is that $G$ won't introduce any new ones, and therefore it won't introduce any new subsets $s$ of $\omega$. Here is the rationale: consider a chain of elements of $\mathbb{P}$, $p_1 < p_2 < p_3 < \cdots$. Given $f : \omega \to (2^{\aleph_0})^{\mathrm{M}}$ in $\mathrm{M}$, and given a partition of $\aleph_1^{\mathrm{M}}$ in two parts mapping to $0$ and $1$, each $p_i$ partially defines a subset $s_f$ of $\omega$ by $p_i \circ f: \omega \to \{0,1\}$. In $\mathrm{M}[G]$, the complete subset $s_f$ defined by $\bigcup G |_{\mathrm{Ran}(f)} \circ f$ is present. However, notice that $\bigcup G |_{\mathrm{Ran}(f)} = \bigcup_i p_i \circ f$, where $p := \bigcup_i p_i$ is just a countable partial bijection $(2^{\aleph_0})^{\mathrm{M}} \to \aleph_1^{\mathrm{M}}$. Because each $p_i \in \mathrm{M}$ we have $p \in \mathrm{M}$ and consequently $p \in \mathbb{P}$. Therefore $s_f$ is already a member of $\mathrm{M}$!

We can formalize this with a lemma:

\textbf{Definition.} A poset $\mathbb{P}$ is \textbf{$\omega$-closed} if for any infinite chain $p_1 < p_2 < p_3 <\cdots$ of elements $p_i \in \mathbb{P}$, $\bigcup_i p_i \in \mathbb{P}$.

\textbf{Lemma.} Let $\mathbb{P}$ be  ($\omega$-closed)$^\mathrm{M}$ and let $X \in \mathrm{M}$. Any function $f: \omega \to X$ in $\mathrm{M}[G]$ is already in $\mathrm{M}$.

\textbf{Proof.} Let $S\in \mathrm{M}$ be the set of all functions $\omega \to X$ in $\mathrm{M}$. Let $\tau$ be a $\mathbb{P}$-name for $f$ and let $p$ force ``$\tau$ is a function $\check{\omega} \to \check{X}$''. If we show that the set of $r$ that force $\tau \in \check{S}$ is dense above $p$, i.e., that at $p$ it is necessarily possible that $f$ is in $S$, then $(f\in S)^{\mathrm{M}[G]}$ and therefore $f\in S$.

Because $p$ forces ``$\tau$ is a function $\check{\omega} \to \check{X}$'', it forces ``$\forall x \in \check{\omega} \:\exists y \in \check{X}$ such that $(x,y)$ is a pair of values in $\tau$''. Therefore for any $n\in \omega$ and  $\sigma = \check{n}$, the set of $q$ that force ``$\exists y \in \check{X}$ such that $(\sigma,y)$ is a pair of values in $\tau$'' is dense above $p$. Therefore the set of $r$ that force ``$(\check{n}, \check{c_n})$ is a pair of values in $\tau$'' for some $c_n \in X$ is dense above $q$.

Now starting with any $q \geq p$ let $p_0 \geq q$ force ``$\tau$ maps $\check{0}$ to $\check{c_0}$'' for some $\check{c_0}$, let $p_1 \geq p_0$ force ``$\tau$ maps $\check{1}$ to $\check{c_1}$'' for some $\check{c_1}$, and so on defining a sequence $p_0 \leq p_1 \leq p_2 \leq \cdots$. Since ($\mathbb{P}$ is $\omega$-closed)$^{\mathrm{M}}$, $r := \bigcup p_i \in \mathbb{P}$ and $r$ forces ``$(\check{n}, \check{c_n})$ is a pair of values in $\tau$'' for every $n\in \omega$ therefore $r\Vdash \tau \in \check{S}$. Because $q$ was arbitrary, the set of $r$ that force $\tau \in \check{S}$ is indeed dense above $p$ and therefore $f\in S$. $\square$

\textbf{Theorem.} $CH^{\mathrm{M}[G]}$ is a ZFC+ theorem.

\textbf{Proof.} By the lemma above, in $\mathrm{M}[G]$ there are no new functions $f: \omega \to X$ for any $X \in \mathrm{M}$. Therefore there are no new subsets of $\aleph_0$ and moreover there is no surjection $\aleph_0 \to \aleph_1^\mathrm{M}$, which means that $(2^\aleph_0)^{\mathrm{M}[G]} = (2^\aleph_0)^\mathrm{M}$ and $\aleph_1^{\mathrm{M}[G]} = \aleph_1^{\mathrm{M}}$.

In conclusion, $\bigcup G$ is a bijection from $(2^\aleph_0)^{\mathrm{M}[G]}$ to $\aleph_1^{\mathrm{M}[G]}$. $\square$

\section{Alternative forcing proofs of consistency of $\lnot \mathrm{CH}$w}

There are many alternative proofs of consistency of $\lnot \mathrm{CH}$. However, all of them are based on forcing. They share the underlying structure of the poset $\mathbb{P}$ that encodes a collection of new subsets of $\aleph_2$ and some nonstandard version of logic that implements the nonstandard logic of forcing---the ``necessarily possible'', or $\forall q \geq p \: \exists r \geq q$ such that ``$\ldots$''. 

\subsection{Original syntactic approach by Cohen}

The very original construction by Cohen described in two parts in his 1963 and 1964 PNAS papers was perhaps the cleanest after all. In the first part (1963) he defined forcing. Here is part of the definition, with notation changed to match this manuscript:

\begin{quote}
"\textbf{Definition 6.} By induction, we define the concept of ``$p$ forces $\phi$'' as follows:
I. If $r > 0$, $p$ forces $\forall_{\alpha} x\: \phi(x)$ if for all $q \geq p$ $q$ does not force $\lnot \phi(\mathcal{F}_{\beta})$ for $\beta < \alpha$. $p$ forces $\exists_{\alpha} x \:\phi(x)$ if for some $\beta < \alpha$, $p$ forces $\phi(\mathcal{F}_{\beta})$.
[...]
The most important part of Definition 6 is I, the other parts are merely obvious derivatives of it."
\end{quote}

Here, $\forall_{\alpha}$ and $\exists_{\alpha}$ are quantifiers limited to rank $< \alpha$. Describing Cohen's construction of sets $\mathcal{F}$ would be too big a digression.

In Cohen's original construction, the poset $\mathbb{P}$ is defined in the ground model and then truth is assigned in $\mathrm{M}[G]$ as follows. The collection of all sentences of ZFC, allowing for parameters ranging across $\mathrm{M}[G]$, is countable. Enumerate all sentences:
$$\psi_1, \psi_2, \psi_3, \ldots, \psi_n, \ldots$$
Now it is shown that if $p \not \Vdash \phi$ then there is $q \geq p$ such that $q \Vdash \lnot \phi$. Then, starting with $p_0 = \{\}$, at each step $i$ extend to $p_{i+1} \geq p_i$ by considering $\psi_i$. If $p_i \Vdash \phi$, then let $p_{i+1} = p_i$. Otherwise there is $q\geq p_i$ such that $q \Vdash \lnot \psi_i$. Let $p_{i+1} := q$. At the end of the construction we have a \textit{complete sequence} $p_0 \leq p_1 \leq p_2 \leq \cdots$ that decides the truth of every sentence $\psi_i$ in the resulting model.

There are many tricky technical details for all this to work. The ordering of the sentences as well as ranks of each element of the model are quite involved. Notice also that the ground model has to be countable from the perspective of the ``outside'' universe $V$: all possible sentences are enumerated in $V$.

Despite the technical complexities, the main idea is perhaps the simplest: Cohen wanted to add $\aleph_2$-many subsets of $\aleph_0$ and obtain a model $\mathrm{M}[G]$. To achieve that, he picked an order that works for all sentences $\psi_i$ of $\mathrm{M}[G]$ and constructed a $G$ that decides them one-by-one by brute force. The later installations of forcing achieve the same task in a more elegant and less syntactical way by using dense sets, but obscure the honest clarity of the goal.

\subsection{Standard versus nonstandard models}

It is worthwhile to comment on the role of countable models versus nonstandard models. A standard model $M$ of ZFC is one where $\in_M := \in_V \cap (M\times M)$ is true set membership and sentences follow the standard binary logic. To produce the standard model $M[G]$ where $\mathrm{CH}$ breaks, a generic $G$ is needed, which in turn necessitates that $M$ is countable. However, the consistency of $\lnot \mathrm{CH}$ does not strictly require that a standard model is produced. The forcing relation itself provides a notion of $\in$ in the ground model $M^{\mathbb{P}}$ with nonstandard rules---$\sigma$ ``in'' $\tau$ with respect to $p$ if $\sigma^G \in \tau^G$ in $M[G]$ when $p\in G$---and nonstandard logic that defines different truth values of every sentence localized to every $p\in \mathbb{P}$. It can be shown that all axioms are true in this nonstandard assignment of truth---they are true at $p = \{\}$---and therefore all theorems go through as before. Then, it can be shown that $\lnot \mathrm{CH}$ holds at $p = \{\}$ in this nonstandard logic.\footnote{In fact, even a sentence $\phi$ that is not false at $\{\}$ and is true in some $p\in \mathbb{P}$ is automatically shown to be consistent; for $\lnot \mathrm{CH}$ the stronger statement of truth at $p=\{\}$ applies.} In this way, the ground model can be any $M$ or even $V$ itself. With that idea in mind, there are proofs of consistency of $\lnot \mathrm{CH}$ that are based on nonstandard models such as Boolean valued models and S4 modal logical models. In these models, existence of a generic $G$ is not needed. Instead, the properties that such a generic $G$ would have, had it existed, are just the guiding principle in assigning truth values to sentences for every $p$ in the forcing notion $\mathbb{P}$. The forcing notion itself plays pretty much the same role in all these approaches.

\subsection{Boolean valued models}

A Boolean algebra is an algebra on the subsets of a set $S$, with operations $a \land b := a \cap b$, $a \lor b := a \cup b$, $\lnot a := S \setminus a$, where $1 := S$ and $0 := \{\}$. The partial order $\leq$ is defined by $\subseteq$. Boolean algebras generalize truth values: the algebra of $\{0,1\}$ is the simplest Boolean algebra, generated by $S = \{\emptyset\}$.

The poset $\mathbb{P}$ can give rise to a Boolean algebra, and then a Boolean valued model can be defined. The reader is referred to Bell (1977) or Jech (2002) for the full details, as well as to the highly readable exposition by Timothy Chow (2007). Risking confusion of the reader, below is a very brief introduction to Boolean algebras and a few words on how this is done. 

An ideal of a Boolean algebra $B$ is a set $I$ such that (1) $0\in I; 1 \not \in I$; (2) for every $a, b \in I$, $a \lor b \in I$; (2) if $a \in I$ and $b\leq a$ then $b \in I$. A maximal ideal of a Boolean algebra, also called a prime ideal, is one where for every $a \in B$ either $a \in I$ or $\lnot a \in I$.

Often Boolean algebras are defined in terms of the algebraic operations $\lor$ (sometimes denoted by $+$), $\land$ (sometimes denoted by $\cdot$) and $\setminus$ (sometimes denoted by $-$), and then \textbf{Stone's representation theorem} applies: \textit{every Boolean algebra is isomorphic to an algebra of sets}.

The natural partial order on a Boolean algebra allows definition of the supremum and infimum of a set of elements: $\sup\{a_1, \ldots, a_k\} = a_1 \lor \ldots \lor a_k$, which in an algebra of sets is just $\bigcup_1^k a_i$, and similarly $\inf\{a_1, \ldots, a_k\} = a_1 \land \cdots \land a_k = \bigcap_1^k a_i$. These can generalize to infinite joins and meets: $\bigwedge_{a \in X} a := \inf \: X$ and $\bigvee_{a \in X} a := \sup \: X$. However, these may not be elements $\in B$. If $\bigwedge, \bigvee$ exist for all $X \subset B$ then B is called \textbf{complete}. If they exist for any $X$ of cardinality $< \kappa$ where $\kappa$ is a regular cardinal, then $B$ is $\kappa$-\textbf{complete}.

\textit{Every Boolean algebra can be densely embedded in a complete Boolean algebra}, its \textbf{completion}. This is accomplished in a similar way that $\mathbb{R}$ completes $\mathbb{Q}$ by the method of Dedekind cuts. The elements $\bigwedge_X, \bigvee_X$ for any $X\subset B$ are defined similarly to the way $\pi$ is defined as the limit of all $a \in \mathbb{Q} \: \:a < \pi$. 

A separative poset $\mathbb{P}$ can be embedded in a complete Boolean algebra $B$: $\mathbb{P}\subseteq B$, $\leq_{\mathbb{P}}$ agrees with $\leq_{B}$, and $\mathbb{P}$ is dense in $B$. Moreover, $B$ is unique up to isomorphism. If $\mathbb{P}$ is a forcing notion defined on a ground model $M$, the resulting Boolean algebra can be used to define a Boolean-valued model for ZFC.

Each set $x\in M$ is converted to a Boolean-valued set in a structure $M^{B}$ that is analogous to the $\mathrm{M}$-class of $\mathbb{P}$-names. In particular, instead of $x$ containing members $y \in x$, it becomes a collection of values $x(y) \in B$: each membership of a set $y$ in $x$ is associated with a value in $B$. Then, sentences from the ground ones $x = y$ and $x \in y$ up to composites such as $\forall x \: \phi$ are given Boolean values in $B$ instead of true/false. 

Once these definitions are made, $M^B$ is proven to be a model of ZFC in the sense that each axiom evaluates to a Boolean value of 1 and consequently each theorem evaluates to $1$ as well.

The role of the generic ideal $G$ is played by a generic ultrafilter. Filters are the duel of ideals; in this exposition here we reversed $\leq$ of the poset $\mathbb{P}$, and proved the existence of a generic ideal. Once $\leq$ reverses back, the object needed is a generic ultrafilter. This object  defines a homomorphism of the Boolean algebra to the algebra $\{0, 1\}$, i.e., decides binary truth and falsehood of every sentence. Formally, the quotient $M^B / G$ decides truth. If $G$ is generic, $M^B/G$ is shown to be isomorphic to a standard transitive model of ZFC and to contain $G$. And again, if $M$ is countable then a generic ultrafilter $G$ is proven to exist.

An exceptionally interesting proof version by Dana Scott (1967) uses a Boolean valued model of ordinary real algebra axiomatized in just three levels: the reals, functions over reals, and \textit{functionals}---functions over functions of the reals. To quote Scott (1967) ``\textit{To those readers familiar with Cohen's original proof [...] our approach may seem at first very different. Actually it is not. [...] Cohen's definition of forcing could be viewed as a way of assigning Boolean values to formulas.}"

My personal opinion is that the Boolean algebra machinery is an overkill for the forcing methodology. As a matter of taste, I find this approach overpowered and consequently less intuitive. Little is achieved by embedding $\mathbb{P}$ to a complete Boolean algebra because $\mathbb{P}$ already has the needed information. Most of the machinery of the Boolean algebra is unused by forcing. Notably, $\lnot \mathrm{CH}$ evaluates to $1$ as opposed to some other nonzero value. What is used is just the logic of ``necessarily possibly true'', which naturally leads to the concept of genericity. This key idea is obscured in an exposition of forcing with Boolean valued models.

\subsection{Modal logical models}

Modal logic is the study of the deductive behavior of expressions that include the modal ``necessarily'' ($\Box$) and its dual ``possibly'' ($\Diamond := \lnot \Box \lnot$). The semantics of modal logic involves Saul Kripke's frames of \textit{possible worlds} formulated as a directed graph $(W, E)$ of vertices $w$ that represent worlds and edges $(w, w')$ that represent reachability of world $w'$ from world $w$. Each world $w$ has truth assignments for the atomic formulas. The common logical connectives work as usual, and $\Box \phi$ is true at world $w$ if for all $(w,w')\in E$, $\phi$ is true in $w'$. This is the meaning of $\phi$ to be necessarily true at $w$. Importantly, reachability is not transitive or reflexive in principle: all reachable worlds $w'$ from $w$ are exposed with explicit edges $(w, w')\in E$. The edge structure can be encoded in logical rules governing $\Box$ and its dual $\Diamond$. In particular, if $(w,w)$ is always present then $\Box \phi \rightarrow \phi$ is an axiom: if $\phi$ is necessarily true at $w$ and $w$ is reachable from $w$, then $\phi$ is true at $w$.

The basic rules \textbf{K} (for Kripke) apply to all modal structures: (1) \textit{Necessitation deduction rule:} if $\phi$ is a theorem, then it holds in all worlds and therefore $\Box \phi$ is also a theorem; (2) \textit{Distribution axiom}: $\Box (\phi \rightarrow \psi) \rightarrow (\Box \phi \rightarrow \Box \psi)$. Both of these rules hold in any graph structure $(W, E)$ as the reader may wish to verify.

In addition, the rule ``M'', $\Box \phi \rightarrow \phi$, corresponds to $E$ being reflexive, and the rule $\Box \phi \rightarrow \Box \Box \phi$ corresponds to $E$ being transitive. The system resulting from adding these two rules to \textbf{K} is called S4.

In forcing, a poset $\mathbb{P}$ is formed with reachability structure that is reflexive and transitive, neatly fitting in the S4 framework. The whole forcing construction can take place within an S4 modal model. Every formula $\phi$ is embedded in the model by attaching $\Box \Diamond$ in front of it and all its subformulas. For example, $\forall x \: \forall y \: \exists z \: x \in z \land y \in z$ becomes:
$$\Box \Diamond\forall x \: (\Box \Diamond\forall y \: (\Box \Diamond\exists z \:  (\Box \Diamond(\Box \Diamond(x \in z) \land \Box \Diamond(y \in z)))))$$It is shown that every axiom $\phi$ of ZFC is a valid sentence of the S4 modal model if thus embedded, and therefore every theorem of ZFC is a theorem of the modal model under the above embedding.\footnote{Notice the similarity with proving that every axiom $\phi$ is necessarily possibly forced at $p = \{\}\in \mathbb{P}$ with parameters from $\mathrm{M}^{\mathbb{P}}$ using induction on $\phi$'s structure with "necessarily possibly forces" logic.} To show that, similar machinery as the one we described above is employed. Some of the proofs are easier because of the syntactic power of modal logic, and some are in my opinion a bit more cumbersome. To prove Comprehension and Replacement, the forcing relation is again defined with first order formulas just like in this manuscript. Then, it is shown that $\mathrm{CH}$ is not valid by creating the same poset $\mathbb{P}$ of finite partial functions representing $\aleph_2$-many subsets of $\aleph_0$, and building an S4 modal model over it. 

The beauty of this approach is that it formalizes precisely the amount of nonstandard logic needed---no less and no more. The approach described in this manuscript arguably uses less: it deploys ``$\forall q\geq p \: \exists r\geq q \ldots$'' to express ``necessarily possibly''. The Boolean  approach uses too much power: it builds a  complete Boolean model just to use the ``necessarily possibly'' logic within it. The modal approach distills just what is needed.

\section{Discussion}

What does the proof of independence of $\mathrm{CH}$ tell us about the truth value of CH? Here is a list of possibilities: (1) there is a truth value for $\mathrm{CH}$ and we know it; (2) there is a truth value for $\mathrm{CH}$ and we don't know it; (3) there is no truth value for CH; (4) the question is meaningless. The positions among philosophers and set theorists today span all these possibilities. There is extensive literature on the topic, which I will not attempt to summarize other than give the reader a few references for further reading (Shelah 2000, Hamkins 2001, Shelah 2001, Woodin 2001, Woodin 2001b, Bellotti 2005, Woodin 2011, Hamkins 2015, Rittberg 2015, Barton and Friedman 2021).

The device employed in the proof holds a clue. A countable ground model $M$ that contains its versions of $2^{\mathbb{N}}$ and $\aleph_{\alpha}$ is extended with countably many distinct new subsets of $\mathbb{N}$ mapped surjectively to $M$'s countable $\aleph_{\alpha}$ impostor. Unwrapping this, there is a syntactic layer and a semantic layer. The semantic layer concerns $M$'s nature as a model and the membership of certain sets in $M$. However, the existence of $M$ is not needed. It is safe to assume $M$ exists \textit{syntactically} because every proof is finite and employs a finite list of axioms; for any such finite list $\Phi$, there is a ZFC proof that there is a model $M_{\Phi}$. So if an inconsistency proof $\phi_1, \ldots, \phi_k := \emptyset \not = \emptyset$ arises from assuming $M$, then pick the $\Phi$ used in the proof, write down the ZFC proof of existence of $M_{\Phi}$, and copy the proof to $\phi_1^{M_{\Phi}}, \ldots, \phi_k^{M_{\Phi}}$ resulting in an inconsistency proof without assuming $M$. The syntactic layer is an algorithm that shows it is no worse to assume $M$ than to assume ZFC. $M$ is never demonstrated and the ZFC proof of existence of any $M_{\Phi}$ is highly non-constructive. The question of existence of $M$ is not meaningful while the syntactic argument that it is safe to assume $M$ is a terminating algorithm.

From this point on, we are already on shaky ontological ground. Since the ground model is syntactically concrete and semantically dubious, anything we build on top of it is definitely semantically dubious, even if it is syntactically concrete. One may object: what is semantically dubious and who is the arbiter of that? 

Here is another clue. Consider the question of whether ZFC is consistent, call it $Con(ZFC)$. We know from G\"{o}del's second incompleteness theorem that $Con(ZFC)$ is independent of ZFC. Is the statement $Con(ZFC)$ semantically dubious? Not at all: it is the arithmetic statement that there is no finite string of letters $\phi_1, \ldots, \phi_k:= 0=1$ with the requisite syntactic properties. We can write a computer program that will go through all syntactically valid strings and if $\lnot Con(ZFC)$ it will terminate reporting the fact, while the longer the program takes, the more evidence we get in favor of $Con(ZFC)$. The semantic statement \textit{is} the syntactic statement in the case of $Con(ZFC)$ even though $Con(ZFC)$ is independent of ZFC. The situation is  different with $2^{\aleph_0} = \aleph_1$. We already know that the vast majority of members of $2^{\aleph_0}$ can never be listed or named by a formula. That is why any new axiom or argument to further settle the question of the size of $2^{\aleph_0}$ can only hope to do so by further specifying $2^{\aleph_0}$. \textit{Our concept of the real line is not specific, regardless how familiar we feel we are with $\mathbb{R}$.} The independence of $2^{\aleph_0} = \aleph_1$ tells us that the statement is syntactically independent of ZFC, and does not bother telling us anything about what we were hoping to answer---whether $2^{\aleph_0} = \aleph_1$.

Two views that I find especially interesting for philosophical debate are the multiverse view (Hamkins 2012), and countabilism (Barton and Friedman 2021). The multiverse position is that the concept of a set has infinitely many different formulations depending on the axioms chosen, and the solutions to $\mathrm{CH}$ vary in each formulation. According to Hamkins, the multiverse view ``holds that there are diverse distinct concepts of set, each instantiated in a corresponding set-theoretic universe, which exhibit diverse set-theoretic truths.'' This view can be taken syntactically, where each universe is a syntactic construct. However, Hamkins goes further and sees the multiverse as literally existing. ``Each such universe exists independently in the same Platonic sense that proponents of the universe view regard their universe to exist.'' According to this view, the solution to $\mathrm{CH}$ is simply the totality of theorems of how it behaves across the multiverse. Intriguingly, Hamkins proposes the ``Countability Principle'' stating that every universe is countable from the perspective of another, ``better'' universe. We saw this with $\mathrm{M}$ being countable from the perspective of $V$. If we view the universe $V$ as truly existing, an immediate question is whether ``our'' universe---what we perceive as $V$---is countable from the perspective of some outside, better universe that is not accessible to us.

Countabilism (Barton and Friedman 2021) goes in the opposite direction and explores the idea of rejecting the Power axiom and asserting that every set is countable and that the continuum is a proper class.  To quote Scott as also quoted in this paper, ``Perhaps we would be pushed in the end to say that all sets are countable (and that the continuum is not even a set) when at last all cardinals are absolutely destroyed.'' (Scott 1977). The fact that forcing can collapse any cardinal to be countable in a given model is an intuitive argument in favor of this view. As the authors observe, (a) Power implies that there are uncountable sets; (b) given a model $M$ and an $M$-cardinal $\kappa$, there is a forcing notion $\mathbb{P}$ of partial functions $\omega \to \kappa$ forcing $\kappa$ to be countable in the extension; and (c) there is another forcing notion $\mathbb{P}'$ that adds $\kappa$ new subsets of $\omega$ and pushes the continuum to be $\geq \kappa$. The mainstream resolution of this paradoxical situation is that all these impostor models $M$ are missing the ``real'' Power operation. The authors explore instead the removal of Power and addition of alternative axioms that enhance the ability of the remaining ZFC to be a foundation of mathematics while keeping all sets countable and elevating the continuum to a proper class.

Personally I find countabilism attractive. If we accept that every number exists by virtue of being the successor of a previous number that exists all the way down to $0$ and $1$, it makes sense to say that $\mathbb{N}$ exists as simply the collection of all of them. Beyond that, regardless of how syntactically convenient it is to postulate Power and AC for mathematical proofs, these axioms are not self evident. Power in particular is super convenient in allowing definition of higher structures such as arbitrary functions and relations, functionals and so on, but it does let all hell break loose. To get to ``the set of all subsets of $\mathbb{N}$'' one has to first get to a single subset of $\mathbb{N}$. If that subset is generic, it contains Shakespeare's complete works and everything else infinitely many times, so it will be quite a while before one finishes with one subset to get to the next. In perceiving Cantor's infinity of infinities, we may be fooled philosophically by what is nothing more than syntactic sugar. The collection of all possible mathematical statements is countable and therefore our ability to reason about the uncountable is a leap of faith.

\newpage
\section{Acknowledgements}

I thank David McAllester who introduced me to the result when I was a sophomore intern in his lab at MIT back in 1993, and Bob Givan, his PhD student at the time who was supervising my "research": my task was to read and understand the proof as a step towards machine-verifying it. I did read much of the proof line-by-line but didn't at the time "get" it. In this manuscript, I attempt to provide what could have been helpful for a student like me at the time. I especially thank Vassilis Papakonstantinou and Timothy Chow for providing feedback and edits. I thank Ilias Kastanas and Pekka Kauranen for edits, and Nik Weaver, Scott Aaronson, Richard Durbin and Bob Givan for positive feedback on an earlier version of the manuscript.

\newpage

\section{Formula Abbreviations}
$s[i] = x$:  $(i, x) \in s$.

\noindent$\mathcal{P}(x)$: the power set of $x$.

\noindent$\mathrm{Aleph}(\kappa, \alpha)$:  $\alpha$ is an ordinal and $\kappa$ is precisely the $\alpha$th infinite cardinal.

\noindent$\mathrm{Bij}(f)$:  $f$ is a bijection.

\noindent$\mathrm{Card}(x):$ $x$ is a cardinal.

\noindent$\mathrm{Cnt}(x)$: $x$ is countable.

\noindent$\mathrm{CntTrans}(x)$: $x$ is countable and transitive.

\noindent$\mathrm{Comp}(\phi, x)$: Comprehension axiom given formula $\phi$ and domain set $x$. 

\noindent$\mathrm{Dom}(f, A)$: The domain of function $f$ is $A$.

\noindent$\mathrm{Dom}(f)$: The domain of function $f$.

\noindent$\mathrm{Empty}(x)$: $x$ is empty.

\noindent$\mathrm{Inf}(x)$: $x$ is infinite.

\noindent$\mathrm{Inj}(f)$: $f$ is an injection.

\noindent$\mathrm{LimOrd}(\alpha)$: $\alpha$ is a limit ordinal.

\noindent$\mathrm{Ord}(\alpha)$: $\alpha$ is an ordinal.

\noindent$\mathrm{Pow}(x, y)$: $y$ is the power set of $x$.

\noindent$\mathrm{Ran}(f, B)$: The range of function $f$ is $B$.

\noindent$\mathrm{Ran}(f)$: The range of function $f$.

\noindent$\mathrm{Succ}(x, y)$: $y = x \cup \{x\}$

\noindent$\mathrm{SuccOrd}(\alpha)$: $\alpha$ is a successor ordinal.

\noindent$\mathrm{Surj}(f)$: $f$ is a surjection.

\noindent$\mathrm{Trans}(x)$: $x$ is transitive.

\noindent$\mathrm{Uncountable}(x)$: $x$ is uncountable.

\newpage

\section*{References}

\begin{enumerate}

\item Aaronson S. Shtetl-Optimized---The Complete Idiot’s Guide to the Independence of the Continuum Hypothesis: Part 1 of <=$\aleph_0$, 2020. https://scottaaronson.blog/?p=4974
\item Arjona M, Alonso E. Completeness: from Gödel to Henkin. History and Philosophy of Logic DOI: 10.1080/01445340.2013.816555, 2014.
\item Barton N and Friedman S-D. Countabilism and maximality principles. Paper in preparation/under review, https://philpapers.org/rec/BARCAM-5, 2021.
\item Batzoglou S. G\"{o}del's incompleteness theorem. arXiv:2112.06641, 2021.
\item Bell JL. Boolean-valued models and independence proofs in set theory. 2nd edn, Oxford University Press, New York, 1977.
\item Bellotti L. Woodin on the Continuum Problem: an overview and some objections. Logic and Philosophy of Science. 3(1) 2005.
\item Buldt B. The Scope of Gödel’s First Incompleteness Theorem. Logica Universalis (8):499-552, 2014.
\item Chow T. A beginner's guide to focing.	arXiv:0712.1320, 2007.
\item Cohen PJ. The independence of the continuum hypothesis. Proceedings of the National Academy of Sciences U S A. 50(6): 1143–1148, 1963.
\item Cohen PJ. The independence of the continuum hypothesis, II. Proceedings of the National Academy of Sciences U S A. 51(1): 105–110, 1964.
\item Cohen PJ. Set theory and the continuum hypothesis. New York: W.A. Benjamin, 1966.
\item G\"{o}del K. Consistency-proof for the generalized continuum hypothesis. Proceedings of the National Academy of Sciences USA 25: 220-224, 1939.
\item Hamkins JD. The set-theoretic multiverse. The Review of Symbolic Logic 5(3): 416-440, 2012.
\item Hamkins JD. Is the dream solution of the continuum hypothesis attainable?. Notre Dame Journal of Formal Logic 56(1):135-45, 2015.
\item Jech T. Set theory. Springer, 3rd edition, 2002.
\item Kim B. Complete proofs of G\"{o}del's incompleteness theorems. https://web.yonsei.ac.kr/bkim/goedel.pdf.
\item Kunen, Kenneth. Set Theory: An Introduction to Independence Proofs. Elsevier. ISBN 0-444-86839-9. 1980.
\item Mancosu P. Measuring the size of infinite collections of natural numbers: was Cantor's theory of infinite number inevitable? The Review of Symbolic Logic, 2(4), 612-646. doi:10.1017/S1755020309990128, 2009.
\item Nagel E, Newman JR. G\"{o}del's proof. NYU Press; Revised ed. 2001.
\item Rittberg CJ. How Woodin changed his mind: new thoughts on the Continuum Hypothesis. Archive for History of Exact Sciences. 69(2):125-51, 2015.
\item Scott D. A proof of the independence of the continuum hypothesis. Mathematical systems theory 1.2: 89-111, 1967.
\item Scott D. Foreword to boolean-valued models and independence proofs. In (Bell 2011), p. xiii–xviii. Oxford University 1977.
\item Shelah S. The generalized continuum hypothesis revisited. Israel Journal of Mathematics. 116(1):285-321, 2000.
\item Shelah S. You can enter Cantor's paradise. arXiv preprint math/0102056, 2001.
\item Smith P.  An Introduction to Gödel's Theorems, 2nd ed, Cambridge UK,  2013.
\item Smullyan RM and Fitting M. Set theory and the continuum problem. Dover Publications, New York, 2010.
\item Weaver N. Forcing for mathematicians. World Scientific Publishing Co, Singapore, 2014.
\item Woodin WH. The continuum hypothesis, I. Notices Amer. Math. Soc. 48(6):567-76, 2001.
\item Woodin WH. The continuum hypothesis, II. Notices Amer. Math. Soc. 48(7):681-90, 2001b.
\item Woodin WH. The continuum hypothesis, the generic multiverse of sets, and the $\omega$ conjecture. Set theory, arithmetic, and foundations of mathematics: theorems, philosophies 36, 2011.
\end{enumerate}

\end{document}